\documentclass[11pt]{article}

\usepackage{amsfonts}
\usepackage{amscd}
\usepackage{amssymb}
\usepackage{amsthm}
\usepackage{amsmath}

\input prepictex
\input pictex
\input postpictex

 \theoremstyle{plain}
\newtheorem{thm}{Theorem}[section]
\newtheorem{lemma}[thm]{Lemma}
\newtheorem{prop}[thm]{Proposition}

\theoremstyle{definition}

\newtheorem{remark}[thm]{Remark}
\theoremstyle{example}

\theoremstyle{remark}

\numberwithin{equation}{section}

\def\KK{\mathbb{K}}

\def\QQ{\mathbb{Q}} 
\def\RR{\mathbb{R}}

\def\ZZ{\mathbb{Z}}

\def\fh{\mathfrak{h}}

\def\fp{\mathfrak{p}}

\def\wt{\mathrm{wt}}

\def\mapleft#1{\smash{
   \mathop{\longleftarrow}\limits^{#1}}}

\makeatletter
\renewcommand{\@makefnmark}{\mbox{\textsuperscript{}}}
\makeatother

\title{Alcove walks, Hecke algebras, spherical functions, \\
crystals and
column strict tableaux}
\author{Arun Ram \\
Department of Mathematics\\
University of Wisconsin\\
Madison, WI 53706 \\
ram@math.wisc.edu}
\date{}

\begin{document}

\maketitle

\centerline{{\sl Dedicated to R. MacPherson on the occasion of his 60th birthday}}

\begin{section}{Introduction}

Together, Sections 2 and  5 of this paper form a self contained treatment of
the theory of crystals and the path model.
It is my hope that this will be useful to the many people
who, over the years, have told me that they wished they understood crystals
but have found the existing literature too daunting.
One goal of the presentation here is to clarify the
relationship between the general  path model and the crystal operators of Lascoux and
Sch\"utzenberger used in the type A case [LS].
More specifically, Section 2 is a basic pictorial  exposition of
Weyl groups and affine Weyl groups
and Section 5 is an exposition of the theory of (a) symmetric
functions, (b) crystals and (c) the path model
which is designed for readers whose only background is the material
in Section 2.   These two sections can be read independently of Sections 3 and 4. 

Sections 3 and 4 give an exposition of the affine Hecke algebra
and recent results regarding the combinatorics of spherical functions 
on p-adic groups (Hall-Littlewood polynomials) using only the material 
in Section 2 as background.   The $q$-analogue of the theory of crystals developed in 
Section 4 specializes to the path model version [Li1-3] of the 
``classical'' crystal theory at $q^{-1}=0$.  This specialization property is 
``to be expected'' since Macdonald's formula for the spherical function specializes 
to the Weyl character formula at $q^{-1}=0$.  It is my hope that the presentation of the 
results will illustrate the close connection between the affine Hecke algebra,
the path model, and the theory of crystals.

The motivation for this paper is the following.
The classical path model is a combinatorial tool for working with
crystal bases of integrable representations of symmetrizable
Kac-Moody algebras, a generalization of column strict Young 
tableaux.  The same combinatorics can be 
used to study the equivariant K-theory of the flag variety [PR]
and this point provides a translation between
the path model and the structure of the nil-affine Hecke algebra.
Thus, conceptually,
\begin{align*}
\hbox{crystal bases} &=~
\hbox{the path model} \\
&=~
\hbox{the nil affine Hecke algebra} \\
&=~\hbox{the $T$-equivariant $K$-theory of $G/B$.}
\end{align*}
The paper of Gaussent-Littelmann [GL] 
(see also the work of Kapovich-Millson in [KM])
indicated that there is a refined model which removes the `nil' in this chain
of equalities and models the geometry of Kac-Moody groups over a local field.
The result is the conceptual chain of equalities
\begin{align*}
\hbox{positively folded alcove walks} 
&=~\hbox{geodesics in the affine building} \\
&=~\hbox{the affine Hecke algebra} \\
&=~\hbox{MV cycles in the loop Grassmanian.}
\end{align*}
The connection to the affine Hecke algebra and the approach to spherical
functions for a $p$-adic group in [NR]
was made concrete by C. Schwer [Sc] who told me that ``the periodic Hecke
module encodes the positively folded galleries''.

The main results of this paper are obtained by viewing the affine Hecke algebra
in terms of a new construct, the alcove walk algebra.  This combinatorics provides 
\begin{enumerate}
\item[(a)] a change of basis formula in the affine Hecke algebra which is an
analogue of the change of basis formula for the nil affine Hecke algebra
given in [GR],
\item[(b)] a derivation of the 
combinatorial formula for the Hall-Littlewood polynomials in terms
of monomials given by Schwer [Sc],
\item[(c)] a derivation of the product formula for Hall-Littlewood polynomials
given by Kapovich-Millson [KM] and Schwer [Sc],
\item[(d)] a description of the restriction rule for Hall-Littlewood polynomials
\end{enumerate}
in the same way that the theory of
crystal bases gives expansions for Weyl characters (Schur functions)
in terms of monomials and produces general ``Littlewood-Richardson
rules'' for Weyl characters.
The result in (a) is a generalization of a result given in Schwer in the same spirit
as the $q^{-1}=0$ result given in the paper of Lenart-Postnikov [LP].  Thus
(a) generalizes both the result in [Sc] and the result in [LP].  Results
(b) and (c) are all either explicit or implicit in [GL], [KM] and [Sc].
The result in (d) has not appeared in the literature before but it is an
easy consequence of the combination of the techniques of [Sc] and the 
``crystal structure'' on positively folded galleries implicit
in [GL].  The $q$-crystal structure gives easy proofs of the results in (c) and (d).

Unfortunately, the limitations of time and space have not allowed me to
provide an exposition of the connection between the alcove walk combinatorics 
and the geometry of the loop group (which is implicit in the affine Hecke algebra 
and explicit in [Mac], [BD], [GL], [KM] and [Sc]).  Hopefully this will be done soon in a future
paper.

\bigskip\noindent
\textbf{Acknowledgments.}  I thank P. Littelmann for useful conversations and C. Schwer for 
teaching me the right way to think about positively folded galleries.   I thank the
National Science Foundation (grants DMS-0353038 and DMS-0097977),
the Vilas Foundation and the Max-Planck-Institut f\"ur Mathematik for support of this work.
I thank my hosts in Rome, C.\ DeConcini, C.\ Procesi and E.\ Strickland, for their hospitality
during a sabbatical when the writing of this paper was completed.
It is a pleasure to dedicate this paper to R. MacPherson.  He is an inspiration and a role model.

\end{section}

\begin{section}{The affine Weyl group}

This section is a summary of the main facts and notations that are needed for 
working with the affine Weyl group $\tilde W$.   The main point is that the 
elements of the affine Weyl group can be identified with alcoves via the bijection in (2.11).

Let $\fh_\RR^*$ be a finite dimensional vector space over $\RR$.
A \emph{reflection} is a diagonalizable element of $GL(\fh_\RR^*)$ which has exactly
one eigenvalue not equal to $1$.  A \emph{lattice} is a free $\ZZ$-module.
A \emph{Weyl group} is a finite subgroup $W$ of $GL(\fh^*_{\RR})$ which is
$$\hbox{generated by reflections}\qquad\hbox{and}\qquad
\hbox{acts on a lattice $L$ in $\fh_{\RR}^*$}
$$
such that $\fh_\RR^* = L\otimes_\ZZ \RR$.
Let $R^+$ be an index set for the reflections in $W$ so that, for 
$\alpha\in R^+$,
$$s_\alpha
\quad\hbox{is the reflection in the hyperplane}\quad
H_\alpha = (\fh_\RR^*)^{s_\alpha},
$$
the fixed point space of the transformation $s_\alpha$.
The \emph{chambers} are the connected components of the complement
$$\fh^*_\RR\backslash \Big(\bigcup_{\alpha\in R^+} H_\alpha\big)$$
of these hyperplanes in $\fh_{\RR}^*$.  These are fundamental regions for the action of $W$.

Let $\langle,\rangle$ be a nondegenerate $W$-invariant bilinear form
on $\fh_\RR^*$.
Fix a chamber $C$ and choose vectors $\alpha^\vee\in \fh_\RR^*$ such that 
\begin{equation}
C = \{ x\in \fh_\RR^*\ |\ \langle x,\alpha^\vee\rangle >0\}
\qquad\hbox{and}\qquad
P\supseteq L\supseteq Q,
\end{equation} 
where
\begin{equation}
P = \{ \lambda\in \fh_\RR^*\ |\ \langle\lambda,\alpha^\vee\rangle \in \ZZ\} 
\qquad\hbox{and}\qquad
Q = \sum_{\alpha\in R^+} \ZZ\alpha,
\quad\hbox{where}\quad \alpha = \frac{2\alpha^\vee}{\langle\alpha^\vee,\alpha^\vee\rangle}.
\end{equation}
Pictorially,
$$
\hbox{$\langle\lambda,\alpha^\vee\rangle$\quad is the distance from 
$\lambda$ to the hyperplane $H_{\alpha}$.}
$$

The \emph{alcoves} are the connected components of the complement
$$\fh_\RR^*\backslash \Big(\bigcup_{\alpha\in R^+\atop j\in \ZZ} H_{\alpha,j}\Big)
\qquad\hbox{of the (affine) hyperplanes}\qquad
H_{\alpha,j}= \{ x\in \fh_\RR^*\ |\ \langle x,\alpha^\vee\rangle = j\}
$$
in $\fh_\RR^*$.
The \emph{fundamental alcove} is the alcove
\begin{equation}
A\subseteq C
\qquad\hbox{such that}\qquad 0 \in \overline{A},\end{equation}
where $\overline{A}$ is the closure of $A$.
An example is the case of type $C_2$, where the picture is
$$
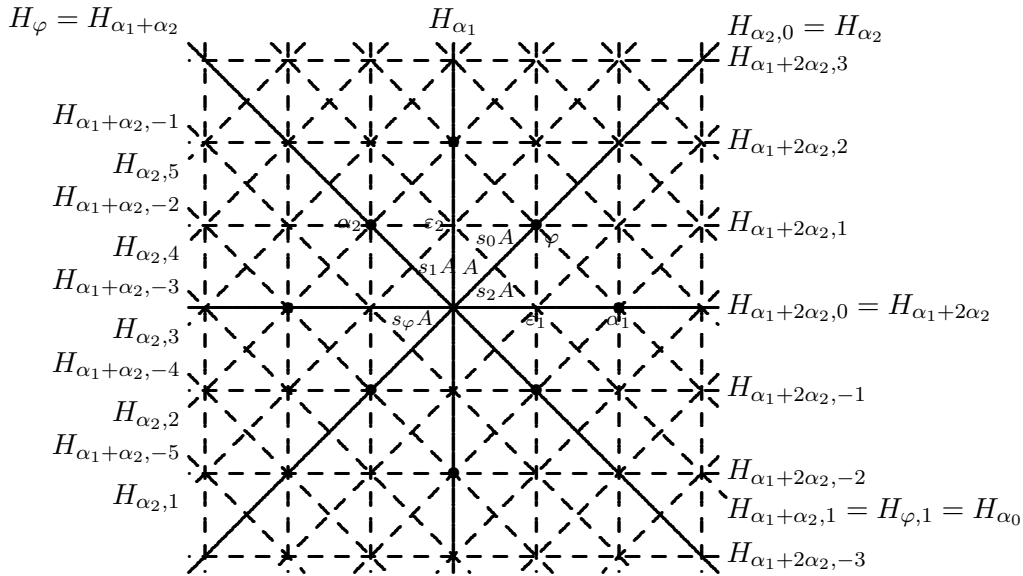
\beginpicture
\setcoordinatesystem units <1.1cm,1.1cm>         
\setplotarea x from -4.5 to 4.5, y from -3.5 to 4    
\put{$H_{\alpha_1}$}[b] at 0 3.3
\put{$H_{\alpha_2,0}=H_{\alpha_2}$}[bl] at 3.3 3.2
\put{$H_{\alpha_2,1}$}[r] at -3.3 -2.3
\put{$H_{\alpha_2,2}$}[r] at -3.3 -1.3
\put{$H_{\alpha_2,3}$}[r] at -3.3 -0.3
\put{$H_{\alpha_2,4}$}[r] at -3.3 0.7
\put{$H_{\alpha_2,5}$}[r] at -3.3 1.7
\put{$H_\varphi = H_{\alpha_1+\alpha_2}$}[br] at -3.3 3.3
\put{$H_{\alpha_1+\alpha_2,1}=H_{\varphi,1}=H_{\alpha_0}$}[tl] at 3.3 -2.3 
\put{$H_{\alpha_1+\alpha_2,-1}$}[r] at -3.3 2.3
\put{$H_{\alpha_1+\alpha_2,-2}$}[r] at -3.3 1.3
\put{$H_{\alpha_1+\alpha_2,-3}$}[r] at -3.3 0.3
\put{$H_{\alpha_1+\alpha_2,-4}$}[r] at -3.3 -0.7
\put{$H_{\alpha_1+\alpha_2,-5}$}[r] at -3.3 -1.7
\put{$H_{\alpha_1+2\alpha_2,0}=H_{\alpha_1+2\alpha_2}$}[l] at 3.3 0
\put{$H_{\alpha_1+2\alpha_2,1}$}[l] at 3.3 1
\put{$H_{\alpha_1+2\alpha_2,2}$}[l] at 3.3 2
\put{$H_{\alpha_1+2\alpha_2,3}$}[l] at 3.3 3
\put{$H_{\alpha_1+2\alpha_2,-1}$}[l] at 3.3 -1
\put{$H_{\alpha_1+2\alpha_2,-2}$}[l] at 3.3 -2
\put{$H_{\alpha_1+2\alpha_2,-3}$}[l] at 3.3 -3
\put{$\scriptstyle{A}$} at 0.2 0.5
\put{$\scriptstyle{s_1A}$} at -0.2 0.5
\put{$\scriptstyle{s_2A}$} at 0.5 0.2
\put{$\scriptstyle{s_0A}$} at 0.5 0.85
\put{$\scriptstyle{s_\varphi A}$} at -0.5 -0.15
\put{$\scriptstyle{\varepsilon_1}$}[t] at 1 -0.1
\put{$\scriptstyle{\varepsilon_2}$}[r] at -0.1 1
\put{$\scriptstyle{\alpha_1}$}[t] at   2  -0.1 
\put{$\scriptstyle{\alpha_2}$}[r] at -1.1   1 
\put{$\scriptstyle{\varphi}$}[tl] at  1.1  0.9  
\put{$\bullet$} at  0  2
\put{$\bullet$} at -1  1
\put{$\bullet$} at  1  1
\put{$\bullet$} at -2  0
\put{$\bullet$} at  2  0
\put{$\bullet$} at -1 -1
\put{$\bullet$} at  1 -1
\put{$\bullet$} at  0 -2
\plot -3.2 -3.2   3.2 3.2 /
\plot  3.2 -3.2  -3.2 3.2 /
\plot  0  3.2   0 -3.2 /
\plot  3.2  0  -3.2  0 /
\setdashes
\plot  -3 3.2   -3 -3.2 /
\plot  -2 3.2   -2 -3.2 /
\plot  -1 3.2   -1 -3.2 /
\plot   1 3.2    1 -3.2 /
\plot   2 3.2    2 -3.2 /
\plot   3 3.2    3 -3.2 /
\plot  3.2 -3   -3.2 -3 /
\plot  3.2 -2   -3.2 -2 /
\plot  3.2 -1   -3.2 -1 /
\plot  3.2  1   -3.2  1 /
\plot  3.2  2   -3.2  2 /
\plot  3.2  3   -3.2  3 /
\plot  3.2 -1.8   1.8 -3.2 /
\plot  3.2 -0.8   0.8 -3.2 /
\plot  3.2 0.2    -0.2 -3.2 /
\plot  3.2 1.2   -1.2 -3.2 /
\plot  3.2 2.2   -2.2 -3.2 /
\plot  2.2 3.2   -3.2 -2.2 /
\plot  1.2 3.2   -3.2 -1.2 /
\plot  0.2 3.2   -3.2 -0.2 /
\plot  -0.8 3.2  -3.2 0.8 /
\plot  -1.8 3.2  -3.2 1.8 /
\plot  -2.2 3.2   3.3 -2.3 /
\plot  -1.2 3.2   3.2 -1.2 /
\plot  -0.2 3.2   3.2 -0.2 /
\plot  0.8 3.2  3.2 0.8 /
\plot  1.8 3.2  3.2 1.8 /
\plot  -3.2 -1.8   -1.8 -3.2 /
\plot  -3.2 -0.8   -0.8 -3.2 /
\plot  -3.2 0.2    0.2 -3.2 /
\plot  -3.2 1.2   1.2 -3.2 /
\plot  -3.2 2.2   2.2 -3.2 /
\endpicture
$$

The \emph{translation} in $\lambda$ is the operator $t_\lambda\colon \fh_\RR^*\to \fh^*_\RR$ given
by
\begin{equation}
t_\lambda(x) = \lambda+x,
\qquad\hbox{for $\lambda\in P, x\in \fh_\RR^*$.}
\end{equation}
The reflection $s_{\alpha,k}$ in the hyperplane $H_{\alpha,k}$ is given by
\begin{equation}
s_{\alpha,k} = t_{k\alpha}s_\alpha = s_\alpha t_{-k\alpha}.
\end{equation}
The \emph{extended affine Weyl group} is 
\begin{equation}
\widetilde{W}= P\rtimes W = \{ t_\lambda w\ |\ \lambda\in P, w\in W\}
\qquad\hbox{with}\qquad
wt_{\lambda} =  t_{w\lambda}w.
\end{equation}
Denote the walls of $C$ by $H_{\alpha_1},\ldots, H_{\alpha_n}$ and extend this indexing so that 
$$H_{\alpha_0},\ldots, H_{\alpha_n}
\qquad\hbox{are the walls of}\qquad 
A,$$
the fundamental alcove.  Then the \emph{affine Weyl group},
\begin{equation}
W_{\mathrm{aff}} = Q\rtimes W
\quad\hbox{is generated by $s_0,\ldots, s_n$,}
\end{equation}
the reflections in the hyperplanes $H_{\alpha_0},\ldots, H_{\alpha_n}$.
Furthermore, $A$ is a fundamental region
for the action of $W_{\mathrm{aff}}$ on $\fh_\RR^*$ and so there is a bijection
$$\begin{matrix}
W_{\mathrm{aff}} &\longrightarrow &\{\hbox{alcoves in $\fh_\RR^*$}\} \\
w &\longmapsto &w^{-1}A.
\end{matrix}
$$
The \emph{length} of $w\in \widetilde{W}$ is 
\begin{equation}
\ell(w) = \hbox{number of hyperplanes between $A$ and $wA$.}
\end{equation}
The difference between $W_{\mathrm{aff}}$ and $\widetilde{W}$ is the group
\begin{equation}
\Omega = \widetilde{W}/W_{\mathrm{aff}}\cong P/Q.
\end{equation}
The group $\Omega$ is the set of elements of $\widetilde{W}$ of length $0$.
An element of $\Omega$ acts on the fundamental
alcove $A$ by an automorphism.  Its action on $A$ induces a permutation of the walls of $A$,
and hence a permutation of $0,1,\ldots, n$.  
If $g\in \Omega$ and $g\ne 1$ let $\omega_i$ be the image of the origin under the action of $g$ on $A$.
If $s_j$ denotes the reflection in the 
$j$th wall of $A$ and $w_i$ denotes the longest element of the stabilizer
$W_{\omega_i}$ of $\omega_i$ in $W$, then
\begin{equation}
g s_i g^{-1} = s_{g(i)}
\qquad\hbox{and}\qquad
gw_0 w_i = t_{\omega_i}.
\end{equation}

The group $\widetilde{W}$ acts freely on  $\Omega\times \fh_\RR^*$
($|\Omega|$ copies of $\RR^n$ tiled by alcoves) so that $g^{-1}A$ is in the
same spot as $A$ except on the $g$th ``sheet'' of $\Omega\times \fh_\RR^*$.
It is helpful to think of the elements of $\Omega$ as the \emph{deck transformations}
which transfer between the sheets in $\Omega\times \fh_\RR^*$.
Then
\begin{equation}
\begin{matrix}
\widetilde W &\longrightarrow &\{\hbox{alcoves in $\Omega\times \fh_\RR^*$}\} \\
w &\longmapsto &w^{-1}A 
\end{matrix}
\end{equation}
is a bijection.
In type $C_2$, the two sheets in $\Omega\times \fh_\RR^*$ look like
$$
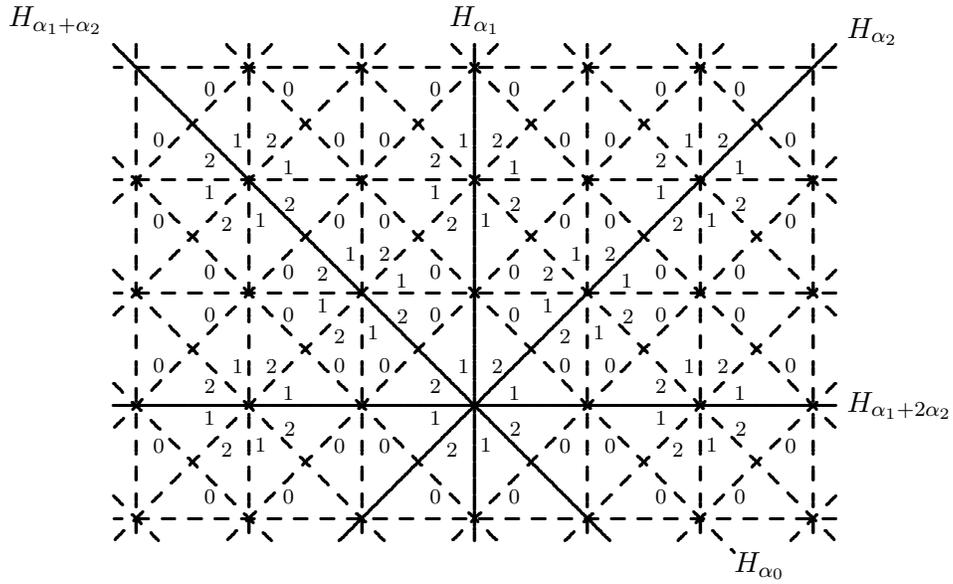
\beginpicture
\setcoordinatesystem units <1.5cm,1.5cm>         
\setplotarea x from -4.5 to 4.5, y from -1.5 to 4    
\put{$H_{\alpha_1}$}[b] at 0 3.3
\put{$H_{\alpha_2}$}[bl] at 3.3 3.2
\put{$H_{\alpha_1+\alpha_2}$}[br] at -3.3 3.3
\put{$H_{\alpha_0}$}[tl] at 2.3 -1.3 
\put{$H_{\alpha_1+2\alpha_2}$}[l] at 3.3 0
\put{$\scriptstyle{0}$}[bl] at 0.3 0.75
\put{$\scriptstyle{0}$}[br] at -0.3 0.75
\put{$\scriptstyle{0}$}[bl] at 0.75  0.3
\put{$\scriptstyle{0}$}[br] at -0.75 0.3
\put{$\scriptstyle{0}$}[tl] at 0.3 -0.75
\put{$\scriptstyle{0}$}[tr] at -0.3 -0.75
\put{$\scriptstyle{0}$}[tl] at 0.75 -0.3
\put{$\scriptstyle{0}$}[tr] at -0.75 -0.3
\put{$\scriptstyle{1}$}[r] at -0.05 0.35
\put{$\scriptstyle{2}$}[tr] at -0.3 0.23
\put{$\scriptstyle{1}$}[t] at -0.35 -0.05
\put{$\scriptstyle{2}$}[tl] at -0.25 -0.33
\put{$\scriptstyle{1}$}[l] at 0.05 -0.35
\put{$\scriptstyle{2}$}[bl] at 0.31 -0.27
\put{$\scriptstyle{1}$}[b] at 0.35 0.05
\put{$\scriptstyle{2}$}[r] at 0.25  0.35
\put{$\scriptstyle{0}$}[bl] at 2.3 0.75
\put{$\scriptstyle{0}$}[br] at 1.7 0.75
\put{$\scriptstyle{0}$}[bl] at 2.75  0.3
\put{$\scriptstyle{0}$}[br] at 1.25 0.3
\put{$\scriptstyle{0}$}[tl] at 2.3 -0.75
\put{$\scriptstyle{0}$}[tr] at 1.7 -0.75
\put{$\scriptstyle{0}$}[tl] at 2.75 -0.3
\put{$\scriptstyle{0}$}[tr] at 1.25 -0.3
\put{$\scriptstyle{1}$}[r] at 1.95 0.35
\put{$\scriptstyle{2}$}[tr] at 1.7 0.23
\put{$\scriptstyle{1}$}[t] at 1.65 -0.05
\put{$\scriptstyle{2}$}[tl] at 1.75 -0.33
\put{$\scriptstyle{1}$}[l] at 2.05 -0.35
\put{$\scriptstyle{2}$}[bl] at 2.31 -0.27
\put{$\scriptstyle{1}$}[b] at 2.35 0.05
\put{$\scriptstyle{2}$}[r] at 2.25  0.35
\put{$\scriptstyle{0}$}[bl] at -1.7 0.75
\put{$\scriptstyle{0}$}[br] at -2.3 0.75
\put{$\scriptstyle{0}$}[bl] at -1.25  0.3
\put{$\scriptstyle{0}$}[br] at -2.75 0.3
\put{$\scriptstyle{0}$}[tl] at -1.7 -0.75
\put{$\scriptstyle{0}$}[tr] at -2.3 -0.75
\put{$\scriptstyle{0}$}[tl] at -1.25 -0.3
\put{$\scriptstyle{0}$}[tr] at -2.75 -0.3
\put{$\scriptstyle{1}$}[r] at -2.05 0.35
\put{$\scriptstyle{2}$}[tr] at -2.3 0.23
\put{$\scriptstyle{1}$}[t] at -2.35 -0.05
\put{$\scriptstyle{2}$}[tl] at -2.25 -0.33
\put{$\scriptstyle{1}$}[l] at -1.95 -0.35
\put{$\scriptstyle{2}$}[bl] at -1.69 -0.27
\put{$\scriptstyle{1}$}[b] at -1.65 0.05
\put{$\scriptstyle{2}$}[r] at -1.75  0.35
\put{$\scriptstyle{0}$}[bl] at 0.3 2.75
\put{$\scriptstyle{0}$}[br] at -0.3 2.75
\put{$\scriptstyle{0}$}[bl] at 0.75  2.3
\put{$\scriptstyle{0}$}[br] at -0.75 2.3
\put{$\scriptstyle{0}$}[tl] at 0.3  1.25
\put{$\scriptstyle{0}$}[tr] at -0.3  1.25
\put{$\scriptstyle{0}$}[tl] at 0.75  1.7
\put{$\scriptstyle{0}$}[tr] at -0.75  1.7
\put{$\scriptstyle{1}$}[r] at -0.05 2.35
\put{$\scriptstyle{2}$}[tr] at -0.3 2.23
\put{$\scriptstyle{1}$}[t] at -0.35  1.95
\put{$\scriptstyle{2}$}[tl] at -0.25  1.67
\put{$\scriptstyle{1}$}[l] at 0.05   1.65
\put{$\scriptstyle{2}$}[bl] at 0.31  1.73
\put{$\scriptstyle{1}$}[b] at 0.35 2.05
\put{$\scriptstyle{2}$}[r] at 0.25  2.35
\put{$\scriptstyle{0}$}[bl] at 2.3 2.75
\put{$\scriptstyle{0}$}[br] at 1.7 2.75
\put{$\scriptstyle{0}$}[bl] at 2.75  2.3
\put{$\scriptstyle{0}$}[br] at 1.25 2.3
\put{$\scriptstyle{0}$}[tl] at 2.3 1.25
\put{$\scriptstyle{0}$}[tr] at 1.7  1.25
\put{$\scriptstyle{0}$}[tl] at 2.75  1.7
\put{$\scriptstyle{0}$}[tr] at 1.25  1.7
\put{$\scriptstyle{1}$}[r] at 1.95  2.35
\put{$\scriptstyle{2}$}[tr] at 1.7  2.23
\put{$\scriptstyle{1}$}[t] at 1.65  1.95
\put{$\scriptstyle{2}$}[tl] at 1.75  1.67
\put{$\scriptstyle{1}$}[l] at 2.05   1.65
\put{$\scriptstyle{2}$}[bl] at 2.31   1.73
\put{$\scriptstyle{1}$}[b] at 2.35  2.05
\put{$\scriptstyle{2}$}[r] at 2.25   2.35
\put{$\scriptstyle{0}$}[bl] at -1.7  2.75
\put{$\scriptstyle{0}$}[br] at -2.3  2.75
\put{$\scriptstyle{0}$}[bl] at -1.25  2.3
\put{$\scriptstyle{0}$}[br] at -2.75  2.3
\put{$\scriptstyle{0}$}[tl] at -1.7  1.25
\put{$\scriptstyle{0}$}[tr] at -2.3  1.25
\put{$\scriptstyle{0}$}[tl] at -1.25  1.7
\put{$\scriptstyle{0}$}[tr] at -2.75  1.7
\put{$\scriptstyle{1}$}[r] at -2.05  2.35
\put{$\scriptstyle{2}$}[tr] at -2.3  2.23
\put{$\scriptstyle{1}$}[t] at -2.35  1.95
\put{$\scriptstyle{2}$}[tl] at -2.25  1.67
\put{$\scriptstyle{1}$}[l] at -1.95  1.65
\put{$\scriptstyle{2}$}[bl] at -1.69  1.73
\put{$\scriptstyle{1}$}[b] at -1.65  2.05
\put{$\scriptstyle{2}$}[r] at -1.75  2.35
\put{$\scriptstyle{1}$}[r] at  0.95 1.35
\put{$\scriptstyle{2}$}[tr] at  0.7 1.23
\put{$\scriptstyle{1}$}[t] at  0.65  0.95
\put{$\scriptstyle{2}$}[tl] at  0.75  0.67
\put{$\scriptstyle{1}$}[l] at 1.05  0.65
\put{$\scriptstyle{2}$}[bl] at 1.31  0.73
\put{$\scriptstyle{1}$}[b] at 1.35 1.05
\put{$\scriptstyle{2}$}[r] at 1.25  1.35
\put{$\scriptstyle{1}$}[r] at -1.05  1.35
\put{$\scriptstyle{2}$}[tr] at -1.3  1.23
\put{$\scriptstyle{1}$}[t] at -1.35  0.95
\put{$\scriptstyle{2}$}[tl] at -1.25  0.67
\put{$\scriptstyle{1}$}[l] at  -0.95  0.65
\put{$\scriptstyle{2}$}[bl] at -0.69  0.73
\put{$\scriptstyle{1}$}[b] at -0.65  1.05
\put{$\scriptstyle{2}$}[r] at -0.75   1.35
\plot -1.2 -1.2   3.2 3.2 /
\plot  1.2 -1.2  -3.2 3.2 /
\plot  0  3.2   0 -1.2 /
\plot  3.2  0  -3.2  0 /
\setdashes
\plot  -3 3.2   -3 -1.2 /
\plot  -2 3.2   -2 -1.2 /
\plot  -1 3.2   -1 -1.2 /
\plot   1 3.2    1 -1.2 /
\plot   2 3.2    2 -1.2 /
\plot   3 3.2    3 -1.2 /
\plot  3.2 -1   -3.2 -1 /
\plot  3.2  1   -3.2  1 /
\plot  3.2  2   -3.2  2 /
\plot  3.2  3   -3.2  3 /
%
\plot  3.2 -0.8   2.8 -1.2 /
\plot  3.2 0.2    1.8 -1.2 /
\plot  3.2 1.2   0.8 -1.2 /
\plot  3.2 2.2   -0.2 -1.2 /
\plot  2.2 3.2   -2.2 -1.2 /
\plot  1.2 3.2   -3.2 -1.2 /
\plot  0.2 3.2   -3.2 -0.2 /
\plot  -0.8 3.2  -3.2 0.8 /
\plot  -1.8 3.2  -3.2 1.8 /
\plot  -2.2 3.2   2.3 -1.3 /
\plot  -1.2 3.2   3.2 -1.2 /
\plot  -0.2 3.2   3.2 -0.2 /
\plot  0.8 3.2  3.2 0.8 /
\plot  1.8 3.2  3.2 1.8 /
%
\plot  -3.2 -0.8   -2.8 -1.2 /
\plot  -3.2 0.2    -1.8 -1.2 /
\plot  -3.2 1.2   -0.8  -1.2 /
\plot  -3.2 2.2   0.2 -1.2 /
\endpicture
$$
and
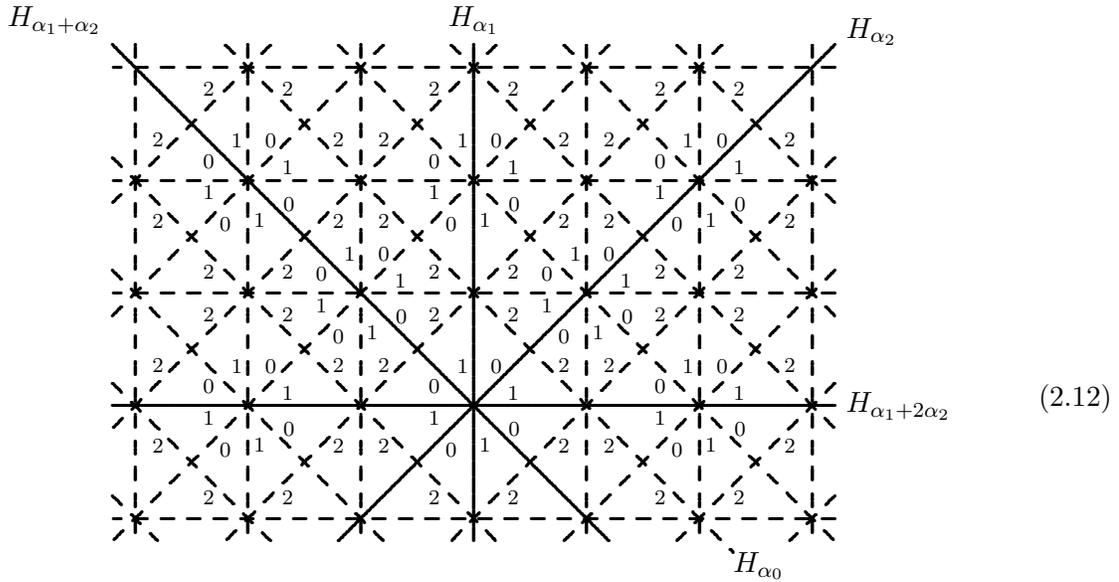
\begin{equation}
\beginpicture
\setcoordinatesystem units <1.5cm,1.5cm>         
\setplotarea x from -4.5 to 4.5, y from -1.5 to 4    
\put{$H_{\alpha_1}$}[b] at 0 3.3
\put{$H_{\alpha_2}$}[bl] at 3.3 3.2
\put{$H_{\alpha_1+\alpha_2}$}[br] at -3.3 3.3
\put{$H_{\alpha_0}$}[tl] at 2.3 -1.3 
\put{$H_{\alpha_1+2\alpha_2}$}[l] at 3.3 0
\put{$\scriptstyle{2}$}[bl] at 0.3 0.75
\put{$\scriptstyle{2}$}[br] at -0.3 0.75
\put{$\scriptstyle{2}$}[bl] at 0.75  0.3
\put{$\scriptstyle{2}$}[br] at -0.75 0.3
\put{$\scriptstyle{2}$}[tl] at 0.3 -0.75
\put{$\scriptstyle{2}$}[tr] at -0.3 -0.75
\put{$\scriptstyle{2}$}[tl] at 0.75 -0.3
\put{$\scriptstyle{2}$}[tr] at -0.75 -0.3
\put{$\scriptstyle{1}$}[r] at -0.05 0.35
\put{$\scriptstyle{0}$}[tr] at -0.3 0.23
\put{$\scriptstyle{1}$}[t] at -0.35 -0.05
\put{$\scriptstyle{0}$}[tl] at -0.25 -0.33
\put{$\scriptstyle{1}$}[l] at 0.05 -0.35
\put{$\scriptstyle{0}$}[bl] at 0.31 -0.27
\put{$\scriptstyle{1}$}[b] at 0.35 0.05
\put{$\scriptstyle{0}$}[r] at 0.25  0.35
\put{$\scriptstyle{2}$}[bl] at 2.3 0.75
\put{$\scriptstyle{2}$}[br] at 1.7 0.75
\put{$\scriptstyle{2}$}[bl] at 2.75  0.3
\put{$\scriptstyle{2}$}[br] at 1.25 0.3
\put{$\scriptstyle{2}$}[tl] at 2.3 -0.75
\put{$\scriptstyle{2}$}[tr] at 1.7 -0.75
\put{$\scriptstyle{2}$}[tl] at 2.75 -0.3
\put{$\scriptstyle{2}$}[tr] at 1.25 -0.3
\put{$\scriptstyle{1}$}[r] at 1.95 0.35
\put{$\scriptstyle{0}$}[tr] at 1.7 0.23
\put{$\scriptstyle{1}$}[t] at 1.65 -0.05
\put{$\scriptstyle{0}$}[tl] at 1.75 -0.33
\put{$\scriptstyle{1}$}[l] at 2.05 -0.35
\put{$\scriptstyle{0}$}[bl] at 2.31 -0.27
\put{$\scriptstyle{1}$}[b] at 2.35 0.05
\put{$\scriptstyle{0}$}[r] at 2.25  0.35
\put{$\scriptstyle{2}$}[bl] at -1.7 0.75
\put{$\scriptstyle{2}$}[br] at -2.3 0.75
\put{$\scriptstyle{2}$}[bl] at -1.25  0.3
\put{$\scriptstyle{2}$}[br] at -2.75 0.3
\put{$\scriptstyle{2}$}[tl] at -1.7 -0.75
\put{$\scriptstyle{2}$}[tr] at -2.3 -0.75
\put{$\scriptstyle{2}$}[tl] at -1.25 -0.3
\put{$\scriptstyle{2}$}[tr] at -2.75 -0.3
\put{$\scriptstyle{1}$}[r] at -2.05 0.35
\put{$\scriptstyle{0}$}[tr] at -2.3 0.23
\put{$\scriptstyle{1}$}[t] at -2.35 -0.05
\put{$\scriptstyle{0}$}[tl] at -2.25 -0.33
\put{$\scriptstyle{1}$}[l] at -1.95 -0.35
\put{$\scriptstyle{0}$}[bl] at -1.69 -0.27
\put{$\scriptstyle{1}$}[b] at -1.65 0.05
\put{$\scriptstyle{0}$}[r] at -1.75  0.35
\put{$\scriptstyle{2}$}[bl] at 0.3 2.75
\put{$\scriptstyle{2}$}[br] at -0.3 2.75
\put{$\scriptstyle{2}$}[bl] at 0.75  2.3
\put{$\scriptstyle{2}$}[br] at -0.75 2.3
\put{$\scriptstyle{2}$}[tl] at 0.3  1.25
\put{$\scriptstyle{2}$}[tr] at -0.3  1.25
\put{$\scriptstyle{2}$}[tl] at 0.75  1.7
\put{$\scriptstyle{2}$}[tr] at -0.75  1.7
\put{$\scriptstyle{1}$}[r] at -0.05 2.35
\put{$\scriptstyle{0}$}[tr] at -0.3 2.23
\put{$\scriptstyle{1}$}[t] at -0.35  1.95
\put{$\scriptstyle{0}$}[tl] at -0.25  1.67
\put{$\scriptstyle{1}$}[l] at 0.05   1.65
\put{$\scriptstyle{0}$}[bl] at 0.31  1.73
\put{$\scriptstyle{1}$}[b] at 0.35 2.05
\put{$\scriptstyle{0}$}[r] at 0.25  2.35
\put{$\scriptstyle{2}$}[bl] at 2.3 2.75
\put{$\scriptstyle{2}$}[br] at 1.7 2.75
\put{$\scriptstyle{2}$}[bl] at 2.75  2.3
\put{$\scriptstyle{2}$}[br] at 1.25 2.3
\put{$\scriptstyle{2}$}[tl] at 2.3 1.25
\put{$\scriptstyle{2}$}[tr] at 1.7  1.25
\put{$\scriptstyle{2}$}[tl] at 2.75  1.7
\put{$\scriptstyle{2}$}[tr] at 1.25  1.7
\put{$\scriptstyle{1}$}[r] at 1.95  2.35
\put{$\scriptstyle{0}$}[tr] at 1.7  2.23
\put{$\scriptstyle{1}$}[t] at 1.65  1.95
\put{$\scriptstyle{0}$}[tl] at 1.75  1.67
\put{$\scriptstyle{1}$}[l] at 2.05   1.65
\put{$\scriptstyle{0}$}[bl] at 2.31   1.73
\put{$\scriptstyle{1}$}[b] at 2.35  2.05
\put{$\scriptstyle{0}$}[r] at 2.25   2.35
\put{$\scriptstyle{2}$}[bl] at -1.7  2.75
\put{$\scriptstyle{2}$}[br] at -2.3  2.75
\put{$\scriptstyle{2}$}[bl] at -1.25  2.3
\put{$\scriptstyle{2}$}[br] at -2.75  2.3
\put{$\scriptstyle{2}$}[tl] at -1.7  1.25
\put{$\scriptstyle{2}$}[tr] at -2.3  1.25
\put{$\scriptstyle{2}$}[tl] at -1.25  1.7
\put{$\scriptstyle{2}$}[tr] at -2.75  1.7
\put{$\scriptstyle{1}$}[r] at -2.05  2.35
\put{$\scriptstyle{0}$}[tr] at -2.3  2.23
\put{$\scriptstyle{1}$}[t] at -2.35  1.95
\put{$\scriptstyle{0}$}[tl] at -2.25  1.67
\put{$\scriptstyle{1}$}[l] at -1.95  1.65
\put{$\scriptstyle{0}$}[bl] at -1.69  1.73
\put{$\scriptstyle{1}$}[b] at -1.65  2.05
\put{$\scriptstyle{0}$}[r] at -1.75  2.35
\put{$\scriptstyle{1}$}[r] at  0.95 1.35
\put{$\scriptstyle{0}$}[tr] at  0.7 1.23
\put{$\scriptstyle{1}$}[t] at  0.65  0.95
\put{$\scriptstyle{0}$}[tl] at  0.75  0.67
\put{$\scriptstyle{1}$}[l] at 1.05  0.65
\put{$\scriptstyle{0}$}[bl] at 1.31  0.73
\put{$\scriptstyle{1}$}[b] at 1.35 1.05
\put{$\scriptstyle{0}$}[r] at 1.25  1.35
\put{$\scriptstyle{1}$}[r] at -1.05  1.35
\put{$\scriptstyle{0}$}[tr] at -1.3  1.23
\put{$\scriptstyle{1}$}[t] at -1.35  0.95
\put{$\scriptstyle{0}$}[tl] at -1.25  0.67
\put{$\scriptstyle{1}$}[l] at  -0.95  0.65
\put{$\scriptstyle{0}$}[bl] at -0.69  0.73
\put{$\scriptstyle{1}$}[b] at -0.65  1.05
\put{$\scriptstyle{0}$}[r] at -0.75   1.35
\plot -1.2 -1.2   3.2 3.2 /
\plot  1.2 -1.2  -3.2 3.2 /
\plot  0  3.2   0 -1.2 /
\plot  3.2  0  -3.2  0 /
\setdashes
\plot  -3 3.2   -3 -1.2 /
\plot  -2 3.2   -2 -1.2 /
\plot  -1 3.2   -1 -1.2 /
\plot   1 3.2    1 -1.2 /
\plot   2 3.2    2 -1.2 /
\plot   3 3.2    3 -1.2 /
\plot  3.2 -1   -3.2 -1 /
\plot  3.2  1   -3.2  1 /
\plot  3.2  2   -3.2  2 /
\plot  3.2  3   -3.2  3 /
%
\plot  3.2 -0.8   2.8 -1.2 /
\plot  3.2 0.2    1.8 -1.2 /
\plot  3.2 1.2   0.8 -1.2 /
\plot  3.2 2.2   -0.2 -1.2 /
\plot  2.2 3.2   -2.2 -1.2 /
\plot  1.2 3.2   -3.2 -1.2 /
\plot  0.2 3.2   -3.2 -0.2 /
\plot  -0.8 3.2  -3.2 0.8 /
\plot  -1.8 3.2  -3.2 1.8 /
\plot  -2.2 3.2   2.3 -1.3 /
\plot  -1.2 3.2   3.2 -1.2 /
\plot  -0.2 3.2   3.2 -0.2 /
\plot  0.8 3.2  3.2 0.8 /
\plot  1.8 3.2  3.2 1.8 /
%
\plot  -3.2 -0.8   -2.8 -1.2 /
\plot  -3.2 0.2    -1.8 -1.2 /
\plot  -3.2 1.2   -0.8  -1.2 /
\plot  -3.2 2.2   0.2 -1.2 /
\endpicture
\end{equation}
where the numbering on the walls of the alcoves is $\widetilde{W}$ equivariant
so that, for $w\in \tilde W$,
the numbering on the walls of $wA$ is the $w$ image of the numbering on the
walls of $A$.

The \emph{0-polygon} is the $W$-orbit of $A$ in $\Omega\times \fh_\RR^*$
and for $\lambda\in P$, the 
$$\hbox{the \emph{$\lambda$-polygon} is $\lambda+WA$,}
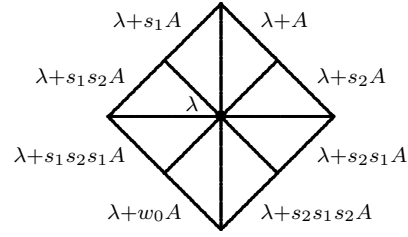
\qquad\qquad
\beginpicture
\setcoordinatesystem units <1.5cm,1.5cm>         
\setplotarea x from -1.5 to 1.5, y from -1.5 to 1.5    
\put{$\bullet$} at 0 0
\put{$\scriptstyle{\lambda+A}$}[l] at 0.35   0.85
\put{$\scriptstyle{\lambda+s_1A}$}[r] at -0.35   0.85
\put{$\scriptstyle{\lambda+s_2A}$}[l] at 0.85   0.35
\put{$\scriptstyle{\lambda+s_1s_2A}$}[r] at -0.85   0.35
\put{$\scriptstyle{\lambda+s_2s_1A}$}[l] at 0.85   -0.35
\put{$\scriptstyle{\lambda+s_1s_2s_1A}$}[r] at -0.85   -0.35
\put{$\scriptstyle{\lambda+s_2s_1s_2A}$}[l] at 0.35   -0.85
\put{$\scriptstyle{\lambda+w_0A}$}[r] at -0.35   -0.85
\put{$\scriptstyle{\lambda}$}[br] at -0.2   0.05
\plot 0 1  1 0 /
\plot 0 -1  1 0 /
\plot  -1 0   0  -1 /
\plot  -1 0   0  1 /
\plot  0 -1  0 1 /
\plot  1 0  -1 0 /
\plot  -0.5 0.5  0.5 -0.5 /
\plot  -0.5 -0.5  0.5 0.5 /
%
%
\endpicture
$$
the translate of the $W$ orbit of $A$ by $\lambda$.
The space
$\Omega\times \fh_\RR^*$ is tiled by the polygons and, via (2.11), we make identifications
between $W$, $\widetilde{W}$, $P$ and their geometric counterparts in $\Omega\times \fh_\RR^*$:
\begin{equation}
\widetilde{W} = \{\hbox{alcoves}\},
\qquad
W = \{\hbox{alcoves in the $0$-polygon}\},
\qquad
P = \{\hbox{centers of polygons}\}.
\end{equation}
Define
\begin{equation}
P^+ = P\cap \overline C
\qquad\hbox{and}\qquad
P^{++} = P\cap C
\end{equation}
so that $P^+$ is a set of representatives of the
orbits of the action of $W$ on $P$.
 The \emph{fundamental weights} are the generators 
$\omega_1,\ldots, \omega_n$ of the
$\ZZ_{\ge 0}$-module $P^+$ so that
\begin{equation}
C = \sum_{i=1}^n \RR_{\ge 0} \omega_i,
\qquad P^+ = \sum_{i=1}^n \ZZ_{\ge 0}\omega_i,
\qquad\hbox{and}\qquad
P^{++} = \sum_{i=1}^n \ZZ_{>0}\omega_i.
\end{equation}
The lattice $P$ has $\ZZ$-basis $\omega_1,\ldots,\omega_n$ and the map
\begin{equation}
\begin{matrix}
P^+ &\longrightarrow &P^{++} \cr
\lambda &\longmapsto &\rho+\lambda, \cr
\end{matrix}
\qquad\hbox{where}\quad
\rho = \omega_1+\ldots+\omega_n,
\end{equation}
is a bijection.
The \emph{simple coroots} are 
$\alpha_1^\vee,\ldots, \alpha_n^\vee$ the dual basis to
the fundamental weights,
\begin{equation}
\langle \omega_i,\alpha_j^\vee\rangle = \delta_{ij}.
\end{equation}
Define
\begin{equation}
\overline{C^\vee} = \sum_{i=1}^n \RR_{\le 0} \alpha_i^\vee
\qquad\hbox{and}\qquad
C^\vee = \sum_{i=1}^n \RR_{< 0} \alpha_i^\vee.
\end{equation}
The {\it dominance order} is the partial order on $\fh_\RR^*$ given by
\begin{equation}
\mu\le \lambda
\qquad\hbox{if}\qquad \mu\in \lambda+\overline{C^\vee}.
\end{equation}

In type $C_2$ the lattice $P= \ZZ\varepsilon_1+\ZZ\varepsilon_2$ with 
$\{\varepsilon_1,\varepsilon_2\}$ an orthonormal basis of $\fh_\RR^*\cong \RR^2$
and $W = \{ 1, s_1, s_2, s_1s_2, s_2s_1, s_1s_2s_1, s_2s_1s_2, 
s_1s_2s_1s_2\}$ is the dihedral group of order $8$ generated by the 
reflections $s_1$ and $s_2$ in the hyperplanes $H_{\alpha_1}$ 
and $H_{\alpha_2}$, respectively, where
$$
H_{\alpha_1} = \{ x\in \fh_\RR^*\ |\ \langle x,\varepsilon_1\rangle = 0\} 
\qquad\hbox{and}\qquad
H_{\alpha_2} = \{ x\in \fh_\RR^*\ |\ \langle x, \varepsilon_2-\varepsilon_1\rangle = 0\}. \\
$$
$$
\beginpicture
\setcoordinatesystem units <1cm,1cm>           
\setplotarea x from -4 to 4, y from -4 to 4    
\put{$H_{\alpha_1}$}[b] at 0 3.1
\put{$H_{\alpha_2}$}[l] at 3.1 3.1
\put{$H_{\alpha_1+\alpha_2}$}[r] at -3.1 3.1
\put{$H_{\alpha_1+2\alpha_2}$}[l] at 4.1 0
\put{$C$}             at  1.5  3
\put{$s_1C$}          at -1.5  3 
\put{$s_2C$}          at  4    1.7
\put{$s_1s_2C$}       at -4    1.7
\put{$s_2s_1C$}       at  4   -1.7 
\put{$s_2s_1s_2C$}    at  1.5 -3
\put{$s_1s_2s_1C$}    at -4   -1.7
\put{$s_1s_2s_1s_2C$} at -1.5 -3
\put{$\bullet$} at -3 -2
\put{$\bullet$} at -2 -2
\put{$\bullet$} at -1 -2
\put{$\bullet$} at  0 -2
\put{$\bullet$} at  1 -2
\put{$\bullet$} at  2 -2
\put{$\bullet$} at  3 -2
\put{$\bullet$} at -3 -1
\put{$\bullet$} at -2 -1
\put{$\bullet$} at -1 -1
\put{$\bullet$} at  0 -1
\put{$\bullet$} at  1 -1
\put{$\bullet$} at  2 -1
\put{$\bullet$} at  3 -1
\put{$\bullet$} at -3 0
\put{$\bullet$} at -2 0
\put{$\bullet$} at -1 0
\put{$\bullet$} at  0 0
\put{$\bullet$} at  1 0
\put{$\bullet$} at  2 0
\put{$\bullet$} at  3 0
\put{$\bullet$} at -3 1
\put{$\bullet$} at -2 1
\put{$\bullet$} at -1 1
\put{$\bullet$} at  0 1
\put{$\bullet$} at  1 1
\put{$\bullet$} at  2 1
\put{$\bullet$} at  3 1
\put{$\bullet$} at -3 2
\put{$\bullet$} at -2 2
\put{$\bullet$} at -1 2
\put{$\bullet$} at  0 2
\put{$\bullet$} at  1 2
\put{$\bullet$} at  2 2
\put{$\bullet$} at  3 2
\put{$\scriptstyle{\alpha_2}$}[l] at -0.9 1.1
\put{$\scriptstyle{\alpha_1}$}[b] at 2.0  0.2  
\put{$\scriptstyle{\alpha_1+\alpha_2}$}[t] at 1.2  0.9  
\put{$\scriptstyle{\alpha_1+2\alpha_2}$}[bl] at -0.3  2.2  
\put{$\scriptstyle{\varepsilon_1}$}[t] at 1 -0.1 
\put{$\scriptstyle{\varepsilon_2}$}[l] at 0.1 1 
\plot -3 -3   3 3 /
\plot  3 -3  -3 3 /
\plot  0  3   0 -3 /
\plot  4  0  -4  0 /
\endpicture
$$
$$
\begin{array}{cc}
\beginpicture
\setcoordinatesystem units <1cm,1cm>         
\setplotarea x from -2 to 3.5, y from -2 to 3.5    
\put{$H_{\alpha_1}$}[b] at 0 3.1
\put{$H_{\alpha_2}$}[l] at 3.1 3.1
\put{$\scriptstyle{0}$}[tr] at -0.05 -0.25
\put{$\scriptstyle{\omega_1}$}[l] at 0.6 0.5
\put{$\scriptstyle{\varepsilon_2=\omega_2}$}[r] at -0.1 0.5 
\put{$\scriptstyle{C}$}             at   1.3   3.2
\put{$\scriptstyle{s_1C}$}          at  -0.9   2.1
\put{$\scriptstyle{s_2C}$}          at   2.3   0.9
\put{$\scriptstyle{s_1s_2C}$}       at  -2.3   0.9
\put{$\scriptstyle{s_2s_1C}$}       at   2.3  -0.9 
\put{$\scriptstyle{s_2s_1s_2C}$}    at   0.9  -2.0
\put{$\scriptstyle{s_1s_2s_1C}$}    at  -2.3  -0.9
\put{$\scriptstyle{s_1s_2s_1s_2C}$} at  -0.9  -2.0
\put{$\bullet$} at  0 0
\put{$\bullet$} at  0 0.5
\put{$\bullet$} at  0.5 0.5
\put{$\bullet$} at  0 1
\put{$\bullet$} at  0.5 1
\put{$\bullet$} at  1 1
\put{$\bullet$} at  0 1.5
\put{$\bullet$} at  0.5 1.5
\put{$\bullet$} at  1 1.5
\put{$\bullet$} at  1.5 1.5
\put{$\bullet$} at  0 2
\put{$\bullet$} at  0.5 2
\put{$\bullet$} at  1 2
\put{$\bullet$} at  1.5 2
\put{$\bullet$} at  2 2
\put{$\bullet$} at  0 2.5
\put{$\bullet$} at  0.5 2.5
\put{$\bullet$} at  1 2.5
\put{$\bullet$} at  1.5 2.5
\put{$\bullet$} at  2 2.5
\put{$\bullet$} at  2.5 2.5
\plot  0.5 0  0.5 0.1 /
\put{$\scriptstyle{\varepsilon_1}$}[t] at 0.5 -0.1
\plot -2 -2   3 3 /
\plot  2 -2  -1.7 1.7 /
\plot  0  3   0 -2.3 /
\plot  3  0  -3  0 /
\endpicture
&\quad
\beginpicture
\setcoordinatesystem units <1cm,1cm>               
\setplotarea x from -2 to 3.5, y from -2 to 3.5    
\put{$H_{\alpha_1}$}[b] at 0 3.1
\put{$H_{\alpha_2}$}[l] at 3.1 3.1
\put{$\scriptstyle{\rho}$}[r] at 0.4 1.0
\put{$\scriptstyle{C}$}             at   1.3   3.4
\put{$\scriptstyle{s_1C}$}          at  -0.9   2.1
\put{$\scriptstyle{s_2C}$}          at   2.3   0.9
\put{$\scriptstyle{s_1s_2C}$}       at  -2.3   0.9
\put{$\scriptstyle{s_2s_1C}$}       at   2.3  -0.9
\put{$\scriptstyle{s_2s_1s_2C}$}    at   0.9  -2.0
\put{$\scriptstyle{s_1s_2s_1C}$}    at  -2.3  -0.9
\put{$\scriptstyle{s_1s_2s_1s_2C}$} at  -0.9  -2.0
\put{$\bullet$} at  0.5 1
\put{$\bullet$} at  0.5 1.5
\put{$\bullet$} at  1 1.5
\put{$\bullet$} at  0.5 2
\put{$\bullet$} at  1 2
\put{$\bullet$} at  1.5 2
\put{$\bullet$} at  0.5 2.5
\put{$\bullet$} at  1 2.5
\put{$\bullet$} at  1.5 2.5
\put{$\bullet$} at  2 2.5
\put{$\bullet$} at  0.5 3
\put{$\bullet$} at  1 3
\put{$\bullet$} at  1.5 3
\put{$\bullet$} at  2 3
\put{$\bullet$} at  2.5 3
\plot    0 0.5  0.1 0.5 /
\plot  0.5 0    0.5 0.1 /
\put{$\scriptstyle{\varepsilon_1}$}[t] at 0.5 -0.1
\put{$\scriptstyle{\varepsilon_2}$}[r] at -0.1 0.5
\plot -2 -2   3 3 /
\plot  2 -2  -1.7 1.7 /
\plot  0  3   0 -2.3 /
\plot  3  0  -3  0 /
\setdots
\plot  0.5  3   0.5 -2 /
\plot -2.5 -2   2.5 3 /
\endpicture
\cr
\cr
\cr
\hbox{The set $P^+$}
&\hbox{The set $P^{++}$}\end{array}
$$
In this case
$$
\begin{array}{lll}
\omega_1 = \varepsilon_1+\varepsilon_2, \qquad\qquad
&\alpha_1 = 2\varepsilon_1,\qquad\qquad 
&\alpha_1^\vee = \varepsilon_1,\\
\omega_2 = \varepsilon_2, \qquad\qquad
&\alpha_2 = \varepsilon_2-\varepsilon_1, \qquad\qquad
&\alpha_2^\vee = \alpha_2,
\end{array}
$$
and
$$
R = \{\pm \alpha_1, \pm \alpha_2, \pm(\alpha_1+\alpha_2),
\pm(\alpha_1+2\alpha_2)\}.
$$

\end{section}

\begin{section}{The affine Hecke algebra}

\begin{subsection}{The alcove walk algebra}

Fix notations for the Weyl group $W$, the extended affine Weyl group $\widetilde{W}$,
and their action on $\Omega\times \fh_\RR^*$ as in Section 2.  Label the walls of the alcoves
so that the fundamental alcove has walls labeled $0,1,\ldots, n$ and the labeling
is $\widetilde{W}$-equivariant (see the picture in (2.12)).

The \emph{periodic orientation} is the orientation of the walls of the alcoves 
given by 
\begin{equation}
\hbox{setting the positive side of $H_{\alpha,j}$}
\qquad\hbox{ to be }\qquad
\{x\in \fh_\RR^*\ |\ \langle x,\alpha^\vee\rangle > j\}.
\end{equation}  
This is an orientation
of the walls of the alcoves such that if $\triangle$ is an alcove and $\lambda\in P$ then
$$\hbox{the walls of $\lambda+\triangle$ have the same orientation as the walls of $\triangle$.}$$

Let $\KK$ be a field.
Use the notations for elements of $\Omega$ as in (2.10).
The \emph{alcove walk algebra} is the algebra over $\KK$ given by generators $g\in \Omega$ and
$$
\begin{matrix}
\beginpicture
\setcoordinatesystem units <1cm,1cm>         
\setplotarea x from -1.5 to 1.5, y from -0.5 to 0.5  
\put{$i$}[b] at 0 0.6
\put{$\scriptstyle{-}$}[b] at -0.4 0.25
\put{$\scriptstyle{+}$}[b] at 0.4 0.25
\plot  0 -0.4  0 0.5 /
\arrow <5pt> [.2,.67] from -0.5 0 to 0.5 0   %
\endpicture
&\beginpicture
\setcoordinatesystem units <1cm,1cm>         
\setplotarea x from -1.5 to 1.5, y from -0.5 to 0.5  
\put{$i$}[b] at 0 0.6
\put{$\scriptstyle{-}$}[b] at -0.4 0.25
\put{$\scriptstyle{+}$}[b] at 0.4 0.25
\plot  0 -0.4  0 0.5 /
\arrow <5pt> [.2,.67] from 0.5 0 to -0.5 0   %
\endpicture
&
\beginpicture
\setcoordinatesystem units <1cm,1cm>         
\setplotarea x from -1.5 to 0.5, y from -0.5 to 0.5  
\put{$i$}[b] at 0 0.7
\put{$\scriptstyle{-}$}[b] at -0.4 0.35
\put{$\scriptstyle{+}$}[b] at 0.4 0.35
\plot  0 -0.4  0 0.6 /
\plot 0.5 0  0.05 0 /
\arrow <5pt> [.2,.67] from 0.05 0.1 to 0.5 0.1   %
\plot 0.05 0 0.05 0.1 /
\endpicture 
&
\beginpicture
\setcoordinatesystem units <1cm,1cm>         
\setplotarea x from -1.5 to 0.5, y from -0.5 to 0.5  
\put{$i$}[b] at 0 0.7
\put{$\scriptstyle{-}$}[b] at -0.4 0.35
\put{$\scriptstyle{+}$}[b] at 0.4 0.35
\plot  0 -0.4  0 0.6 /
\plot -0.5 0  -0.05 0 /
\plot -0.05 0 -0.05 0.1 /
\arrow <5pt> [.2,.67] from -0.05 0.1 to -0.5 0.1   %
\endpicture 
\\
\hbox{positive $i$-crossing} &\hbox{negative $i$-crossing}
&\hbox{positive $i$-fold} &\hbox{negative $i$-fold} 
\end{matrix}
\qquad (1\le i\le n)
$$
with relations (straightening laws)
\begin{equation}
\begin{array}{ccc}
\beginpicture
\setcoordinatesystem units <1cm,1cm>         
\setplotarea x from -0.7 to 0.7, y from -0.5 to 0.5  
\put{$i$}[b] at 0 0.6
\put{$\scriptstyle{-}$}[b] at -0.4 0.25
\put{$\scriptstyle{+}$}[b] at 0.4 0.25
\plot  0 -0.4  0 0.5 /
\arrow <5pt> [.2,.67] from -0.5 0 to 0.5 0   %
\endpicture
=
\beginpicture
\setcoordinatesystem units <1cm,1cm>         
\setplotarea x from -0.7 to 0.7, y from -0.5 to 0.5  
\put{$i$}[b] at 0 0.6
\put{$\scriptstyle{-}$}[b] at -0.4 0.25
\put{$\scriptstyle{+}$}[b] at 0.4 0.25
\plot  0 -0.4  0 0.5 /
\arrow <5pt> [.2,.67] from 0.5 0 to -0.5 0   %
\endpicture
+
\beginpicture
\setcoordinatesystem units <1cm,1cm>         
\setplotarea x from -0.7 to 0.7, y from -0.5 to 0.5  
\put{$i$}[b] at 0 0.7
\put{$\scriptstyle{-}$}[b] at -0.4 0.35
\put{$\scriptstyle{+}$}[b] at 0.4 0.35
\plot  0 -0.4  0 0.6 /
\plot 0.5 0  0.05 0 /
\arrow <5pt> [.2,.67] from 0.05 0.1 to 0.5 0.1   %
\plot 0.05 0 0.05 0.1 /
\endpicture  
&\qquad\hbox{and}\qquad
&\beginpicture
\setcoordinatesystem units <1cm,1cm>         
\setplotarea x from -0.7 to 0.7, y from -0.5 to 0.5  
\put{$i$}[b] at 0 0.6
\put{$\scriptstyle{-}$}[b] at -0.4 0.25
\put{$\scriptstyle{+}$}[b] at 0.4 0.25
\plot  0 -0.4  0 0.5 /
\arrow <5pt> [.2,.67] from 0.5 0 to -0.5 0   %
\endpicture
=
\beginpicture
\setcoordinatesystem units <1cm,1cm>         
\setplotarea x from -0.7 to 0.7, y from -0.5 to 0.5  
\put{$i$}[b] at 0 0.6
\put{$\scriptstyle{-}$}[b] at -0.4 0.25
\put{$\scriptstyle{+}$}[b] at 0.4 0.25
\plot  0 -0.4  0 0.5 /
\arrow <5pt> [.2,.67] from -0.5 0 to 0.5 0   %
\endpicture
+
\beginpicture
\setcoordinatesystem units <1cm,1cm>         
\setplotarea x from -0.7 to 0.7, y from -0.5 to 0.5  
\put{$i$}[b] at 0 0.7
\put{$\scriptstyle{-}$}[b] at -0.4 0.35
\put{$\scriptstyle{+}$}[b] at 0.4 0.35
\plot  0 -0.4  0 0.6 /
\plot -0.5 0  -0.05 0 /
\plot -0.05 0 -0.05 0.1 /
\arrow <5pt> [.2,.67] from -0.05 0.1 to -0.5 0.1   %
\endpicture 
\end{array}
\end{equation}
and
$$
g\left(\beginpicture
\setcoordinatesystem units <1cm,1cm>         
\setplotarea x from -0.6 to 0.6, y from -0.5 to 0.5  
\put{$i$}[b] at 0 0.6
\put{$\scriptstyle{-}$}[b] at -0.4 0.25
\put{$\scriptstyle{+}$}[b] at 0.4 0.25
\plot  0 -0.4  0 0.5 /
\arrow <5pt> [.2,.67] from -0.5 0 to 0.5 0   %
\endpicture\right)
= 
\left(\beginpicture
\setcoordinatesystem units <1cm,1cm>         
\setplotarea x from -0.6 to 0.6, y from -0.5 to 0.5  
\put{$g(i)$}[b] at 0 0.6
\put{$\scriptstyle{-}$}[b] at -0.4 0.25
\put{$\scriptstyle{+}$}[b] at 0.4 0.25
\plot  0 -0.4  0 0.5 /
\arrow <5pt> [.2,.67] from -0.5 0 to 0.5 0   %
\endpicture
\right) g,
\qquad
g\left(\beginpicture
\setcoordinatesystem units <1cm,1cm>         
\setplotarea x from -0.6 to 0.6, y from -0.5 to 0.5  
\put{$i$}[b] at 0 0.6
\put{$\scriptstyle{-}$}[b] at -0.4 0.25
\put{$\scriptstyle{+}$}[b] at 0.4 0.25
\plot  0 -0.4  0 0.5 /
\arrow <5pt> [.2,.67] from 0.5 0 to -0.5 0   %
\endpicture\right)
=\left(\beginpicture
\setcoordinatesystem units <1cm,1cm>         
\setplotarea x from -0.6 to 0.6, y from -0.5 to 0.5  
\put{$g(i)$}[b] at 0 0.6
\put{$\scriptstyle{-}$}[b] at -0.4 0.25
\put{$\scriptstyle{+}$}[b] at 0.4 0.25
\plot  0 -0.4  0 0.5 /
\arrow <5pt> [.2,.67] from 0.5 0 to -0.5 0   %
\endpicture
\right) g,
$$
$$
g\left(\beginpicture
\setcoordinatesystem units <1cm,1cm>         
\setplotarea x from -0.6 to 0.6, y from -0.5 to 0.5  
\put{$i$}[b] at 0 0.7
\put{$\scriptstyle{-}$}[b] at -0.4 0.35
\put{$\scriptstyle{+}$}[b] at 0.4 0.35
\plot  0 -0.4  0 0.6 /
\plot -0.5 0  -0.05 0 /
\plot -0.05 0 -0.05 0.1 /
\arrow <5pt> [.2,.67] from -0.05 0.1 to -0.5 0.1   %
\endpicture\right)
=\left(\beginpicture
\setcoordinatesystem units <1cm,1cm>         
\setplotarea x from -0.6 to 0.6, y from -0.5 to 0.5  
\put{$g(i)$}[b] at 0 0.7
\put{$\scriptstyle{-}$}[b] at -0.4 0.35
\put{$\scriptstyle{+}$}[b] at 0.4 0.35
\plot  0 -0.4  0 0.6 /
\plot -0.5 0  -0.05 0 /
\plot -0.05 0 -0.05 0.1 /
\arrow <5pt> [.2,.67] from -0.05 0.1 to -0.5 0.1   %
\endpicture 
\right) g,
\qquad
g\left(
\beginpicture
\setcoordinatesystem units <1cm,1cm>         
\setplotarea x from -0.6 to 0.6, y from -0.5 to 0.5  
\put{$i$}[b] at 0 0.7
\put{$\scriptstyle{-}$}[b] at -0.4 0.35
\put{$\scriptstyle{+}$}[b] at 0.4 0.35
\plot  0 -0.4  0 0.6 /
\plot 0.5 0  0.05 0 /
\arrow <5pt> [.2,.67] from 0.05 0.1 to 0.5 0.1   %
\plot 0.05 0 0.05 0.1 /
\endpicture\right) 
=
\left(\beginpicture
\setcoordinatesystem units <1cm,1cm>         
\setplotarea x from -0.6 to 0.6, y from -0.5 to 0.5  
\put{$g(i)$}[b] at 0 0.7
\put{$\scriptstyle{-}$}[b] at -0.4 0.35
\put{$\scriptstyle{+}$}[b] at 0.4 0.35
\plot  0 -0.4  0 0.6 /
\plot 0.5 0  0.05 0 /
\arrow <5pt> [.2,.67] from 0.05 0.1 to 0.5 0.1   %
\plot 0.05 0 0.05 0.1 /
\endpicture 
\right) g.$$
Viewing the product as concatenation
each word in the generators can be represented as a sequence of
arrows, with the first arrow having its head or its tail in the fundamental
alcove.
An \emph{alcove walk} is a word in the generators such that, 
\begin{enumerate}
\item[(a)] the tail of the first step is in the fundamental alcove $A$,
\item[(b)] at every step, the head of each arrow is in the same alcove as the tail of the next arrow. 
\end{enumerate}
The \emph{type} of a walk $p$ is the sequence of labels on the arrows.  Note that, 
if $w\in \widetilde{W}$ then
\begin{equation}
\ell(w) = \hbox{length of a minimal length walk from $A$ to $wA$.}
\end{equation}
For example, in type $C_2$,
$$
\beginpicture
\setcoordinatesystem units <1.5cm,1.5cm>         
\setplotarea x from -3.5 to 3.5, y from -1.5 to 4    
\put{$H_{\alpha_1}$}[b] at 0 3.3
\put{$H_{\alpha_2}$}[bl] at 3.3 3.2
\put{$H_{\alpha_1+\alpha_2}$}[br] at -3.3 3.3
\put{$H_{\alpha_0}$}[tl] at 2.3 -1.3 
\put{$H_{\alpha_1+2\alpha_2}$}[l] at 3.3 0
\put{$\scriptstyle{0}$}[bl] at 0.3 0.75
\put{$\scriptstyle{0}$}[br] at -0.3 0.75
\put{$\scriptstyle{0}$}[bl] at 0.75  0.3
\put{$\scriptstyle{0}$}[br] at -0.75 0.3
\put{$\scriptstyle{0}$}[tl] at 0.3 -0.75
\put{$\scriptstyle{0}$}[tr] at -0.3 -0.75
\put{$\scriptstyle{0}$}[tl] at 0.75 -0.3
\put{$\scriptstyle{0}$}[tr] at -0.75 -0.3
\put{$\scriptstyle{1}$}[r] at -0.05 0.35
\put{$\scriptstyle{2}$}[tr] at -0.3 0.23
\put{$\scriptstyle{1}$}[t] at -0.35 -0.05
\put{$\scriptstyle{2}$}[tl] at -0.25 -0.33
\put{$\scriptstyle{1}$}[l] at 0.05 -0.35
\put{$\scriptstyle{2}$}[bl] at 0.31 -0.27
\put{$\scriptstyle{1}$}[b] at 0.35 0.05
\put{$\scriptstyle{2}$}[r] at 0.25  0.35
\put{$\scriptstyle{0}$}[bl] at 2.3 0.75
\put{$\scriptstyle{0}$}[br] at 1.7 0.75
\put{$\scriptstyle{0}$}[bl] at 2.75  0.3
\put{$\scriptstyle{0}$}[br] at 1.25 0.3
\put{$\scriptstyle{0}$}[tl] at 2.3 -0.75
\put{$\scriptstyle{0}$}[tr] at 1.7 -0.75
\put{$\scriptstyle{0}$}[tl] at 2.75 -0.3
\put{$\scriptstyle{0}$}[tr] at 1.25 -0.3
\put{$\scriptstyle{1}$}[r] at 1.95 0.35
\put{$\scriptstyle{2}$}[tr] at 1.7 0.23
\put{$\scriptstyle{1}$}[t] at 1.65 -0.05
\put{$\scriptstyle{2}$}[tl] at 1.75 -0.33
\put{$\scriptstyle{1}$}[l] at 2.05 -0.35
\put{$\scriptstyle{2}$}[bl] at 2.31 -0.27
\put{$\scriptstyle{1}$}[b] at 2.35 0.05
\put{$\scriptstyle{2}$}[r] at 2.25  0.35
\put{$\scriptstyle{0}$}[bl] at -1.7 0.75
\put{$\scriptstyle{0}$}[br] at -2.3 0.75
\put{$\scriptstyle{0}$}[bl] at -1.25  0.3
\put{$\scriptstyle{0}$}[br] at -2.75 0.3
\put{$\scriptstyle{0}$}[tl] at -1.7 -0.75
\put{$\scriptstyle{0}$}[tr] at -2.3 -0.75
\put{$\scriptstyle{0}$}[tl] at -1.25 -0.3
\put{$\scriptstyle{0}$}[tr] at -2.75 -0.3
\put{$\scriptstyle{1}$}[r] at -2.05 0.35
\put{$\scriptstyle{2}$}[tr] at -2.3 0.23
\put{$\scriptstyle{1}$}[t] at -2.35 -0.05
\put{$\scriptstyle{2}$}[tl] at -2.25 -0.33
\put{$\scriptstyle{1}$}[l] at -1.95 -0.35
\put{$\scriptstyle{2}$}[bl] at -1.69 -0.27
\put{$\scriptstyle{1}$}[b] at -1.65 0.05
\put{$\scriptstyle{2}$}[r] at -1.75  0.35
\put{$\scriptstyle{0}$}[bl] at 0.3 2.75
\put{$\scriptstyle{0}$}[br] at -0.3 2.75
\put{$\scriptstyle{0}$}[bl] at 0.75  2.3
\put{$\scriptstyle{0}$}[br] at -0.75 2.3
\put{$\scriptstyle{0}$}[tl] at 0.3  1.25
\put{$\scriptstyle{0}$}[tr] at -0.3  1.25
\put{$\scriptstyle{0}$}[tl] at 0.75  1.7
\put{$\scriptstyle{0}$}[tr] at -0.75  1.7
\put{$\scriptstyle{1}$}[r] at -0.05 2.35
\put{$\scriptstyle{2}$}[tr] at -0.3 2.23
\put{$\scriptstyle{1}$}[t] at -0.35  1.95
\put{$\scriptstyle{2}$}[tl] at -0.25  1.67
\put{$\scriptstyle{1}$}[l] at 0.05   1.65
\put{$\scriptstyle{2}$}[bl] at 0.31  1.73
\put{$\scriptstyle{1}$}[b] at 0.35 2.05
\put{$\scriptstyle{2}$}[r] at 0.25  2.35
\put{$\scriptstyle{0}$}[bl] at 2.3 2.75
\put{$\scriptstyle{0}$}[br] at 1.7 2.75
\put{$\scriptstyle{0}$}[bl] at 2.75  2.3
\put{$\scriptstyle{0}$}[br] at 1.25 2.3
\put{$\scriptstyle{0}$}[tl] at 2.3 1.25
\put{$\scriptstyle{0}$}[tr] at 1.7  1.25
\put{$\scriptstyle{0}$}[tl] at 2.75  1.7
\put{$\scriptstyle{0}$}[tr] at 1.25  1.7
\put{$\scriptstyle{1}$}[r] at 1.95  2.35
\put{$\scriptstyle{2}$}[tr] at 1.7  2.23
\put{$\scriptstyle{1}$}[t] at 1.65  1.95
\put{$\scriptstyle{2}$}[tl] at 1.75  1.67
\put{$\scriptstyle{1}$}[l] at 2.05   1.65
\put{$\scriptstyle{2}$}[bl] at 2.31   1.73
\put{$\scriptstyle{1}$}[b] at 2.35  2.05
\put{$\scriptstyle{2}$}[r] at 2.25   2.35
\put{$\scriptstyle{0}$}[bl] at -1.7  2.75
\put{$\scriptstyle{0}$}[br] at -2.3  2.75
\put{$\scriptstyle{0}$}[bl] at -1.25  2.3
\put{$\scriptstyle{0}$}[br] at -2.75  2.3
\put{$\scriptstyle{0}$}[tl] at -1.7  1.25
\put{$\scriptstyle{0}$}[tr] at -2.3  1.25
\put{$\scriptstyle{0}$}[tl] at -1.25  1.7
\put{$\scriptstyle{0}$}[tr] at -2.75  1.7
\put{$\scriptstyle{1}$}[r] at -2.05  2.35
\put{$\scriptstyle{2}$}[tr] at -2.3  2.23
\put{$\scriptstyle{1}$}[t] at -2.35  1.95
\put{$\scriptstyle{2}$}[tl] at -2.25  1.67
\put{$\scriptstyle{1}$}[l] at -1.95  1.65
\put{$\scriptstyle{2}$}[bl] at -1.69  1.73
\put{$\scriptstyle{1}$}[b] at -1.65  2.05
\put{$\scriptstyle{2}$}[r] at -1.75  2.35
\put{$\scriptstyle{1}$}[r] at  0.95 1.35
\put{$\scriptstyle{2}$}[tr] at  0.7 1.23
\put{$\scriptstyle{1}$}[t] at  0.65  0.95
\put{$\scriptstyle{2}$}[tl] at  0.75  0.67
\put{$\scriptstyle{1}$}[l] at 1.05  0.65
\put{$\scriptstyle{2}$}[bl] at 1.31  0.73
\put{$\scriptstyle{1}$}[b] at 1.35 1.05
\put{$\scriptstyle{2}$}[r] at 1.25  1.35
\put{$\scriptstyle{1}$}[r] at -1.05  1.35
\put{$\scriptstyle{2}$}[tr] at -1.3  1.23
\put{$\scriptstyle{1}$}[t] at -1.35  0.95
\put{$\scriptstyle{2}$}[tl] at -1.25  0.67
\put{$\scriptstyle{1}$}[l] at  -0.95  0.65
\put{$\scriptstyle{2}$}[bl] at -0.69  0.73
\put{$\scriptstyle{1}$}[b] at -0.65  1.05
\put{$\scriptstyle{2}$}[r] at -0.75   1.35
\plot -1.2 -1.2   3.2 3.2 /
\plot  1.2 -1.2  -3.2 3.2 /
\plot  0  3.2   0 -1.2 /
\plot  3.2  0  -3.2  0 /
\arrow <5pt> [.2,.67] from 0.175 0.5 to -0.175 0.5   %
\arrow <5pt> [.2,.67] from -0.175 0.5 to -0.5 0.175   %
\arrow <5pt> [.2,.67] from -0.5 0.175 to -0.825 0.5   %
\arrow <5pt> [.2,.67] from -0.825 0.5 to -1.175 0.5   %
\plot -1.175 0.5   -1.316 0.375  /  %
\plot -1.316  0.375  -1.35 0.42  / %
\arrow <5pt> [.2,.67] from -1.35 0.42 to -1.2 0.55   %
\arrow <5pt> [.2,.67] from -1.2 0.55 to -1.5 0.825   %
\arrow <5pt> [.2,.67] from -1.5 0.825 to -1.5 1.125   %
\arrow <5pt> [.2,.67] from -1.5 1.125 to -1.2 1.5   %
\plot -1.2 1.5  -1.05 1.5 /   %
\plot -1.05 1.5 -1.05 1.55 /  %
\arrow <5pt> [.2,.67] from -1.05 1.55 to -1.22 1.55   %
\arrow <5pt> [.2,.67] from -1.22 1.55 to -1.5 1.82   %
\arrow <5pt> [.2,.67] from -1.5 1.82 to -1.5 2.2   %
\arrow <5pt> [.2,.67] from -1.5 2.2 to -1.825 2.5   %
\setdashes
\plot  -3 3.2   -3 -1.2 /
\plot  -2 3.2   -2 -1.2 /
\plot  -1 3.2   -1 -1.2 /
\plot   1 3.2    1 -1.2 /
\plot   2 3.2    2 -1.2 /
\plot   3 3.2    3 -1.2 /
\plot  3.2 -1   -3.2 -1 /
\plot  3.2  1   -3.2  1 /
\plot  3.2  2   -3.2  2 /
\plot  3.2  3   -3.2  3 /
%
\plot  3.2 -0.8   2.8 -1.2 /
\plot  3.2 0.2    1.8 -1.2 /
\plot  3.2 1.2   0.8 -1.2 /
\plot  3.2 2.2   -0.2 -1.2 /
\plot  2.2 3.2   -2.2 -1.2 /
\plot  1.2 3.2   -3.2 -1.2 /
\plot  0.2 3.2   -3.2 -0.2 /
\plot  -0.8 3.2  -3.2 0.8 /
\plot  -1.8 3.2  -3.2 1.8 /
\plot  -2.2 3.2   2.3 -1.3 /
\plot  -1.2 3.2   3.2 -1.2 /
\plot  -0.2 3.2   3.2 -0.2 /
\plot  0.8 3.2  3.2 0.8 /
\plot  1.8 3.2  3.2 1.8 /
%
\plot  -3.2 -0.8   -2.8 -1.2 /
\plot  -3.2 0.2    -1.8 -1.2 /
\plot  -3.2 1.2   -0.8  -1.2 /
\plot  -3.2 2.2   0.2 -1.2 /
\endpicture
$$
is an alcove walk $p$ of type $(1,2,0,1,0,2,1,2,1,0,1,2)$ with two folds.  Using the notation
\begin{equation}
\begin{array}{llll}
c_i^+\ &\hbox{for a positive $i$-crossing}, \qquad
&f_i^+\ &\hbox{for a positive $i$-fold}, \\
c_i^-\ &\hbox{for a negative $i$-crossing}, 
&f_i^-\ &\hbox{for a negative $i$-fold}, 
\end{array}
\end{equation}
the walk in the picture is
$c_1^-c_2^-c_0^+c_1^-f_0^+c_2^+c_1^+c_2^+f_1^-c_0^+c_1^+c_2^+$.

The proof of the following lemma is straightforward following the 
scheme indicated by the example which follows.

\begin{lemma} The set of alcove walks is a basis of the alcove walk algebra.
\end{lemma}

For example, in type $C_2$, a product of the generators which is not a walk is 
$$c_1^-c_2^+c_0^+c_1^-f_0^-c_2^+c_1^-c_2^+f_1^-c_0^+c_1^+c_2^-,$$
$$
\beginpicture
\setcoordinatesystem units <1.5cm,1.5cm>         
\setplotarea x from -3.5 to 3.5, y from -1.5 to 4    
\put{$H_{\alpha_1}$}[b] at 0 3.3
\put{$H_{\alpha_2}$}[bl] at 3.3 3.2
\put{$H_{\alpha_1+\alpha_2}$}[br] at -3.3 3.3
\put{$H_{\alpha_0}$}[tl] at 2.3 -1.3 
\put{$H_{\alpha_1+2\alpha_2}$}[l] at 3.3 0
\put{$\scriptstyle{0}$}[bl] at 0.3 0.75
\put{$\scriptstyle{0}$}[br] at -0.3 0.75
\put{$\scriptstyle{0}$}[bl] at 0.75  0.3
\put{$\scriptstyle{0}$}[br] at -0.75 0.3
\put{$\scriptstyle{0}$}[tl] at 0.3 -0.75
\put{$\scriptstyle{0}$}[tr] at -0.3 -0.75
\put{$\scriptstyle{0}$}[tl] at 0.75 -0.3
\put{$\scriptstyle{0}$}[tr] at -0.75 -0.3
\put{$\scriptstyle{1}$}[r] at -0.05 0.35
\put{$\scriptstyle{2}$}[tr] at -0.3 0.23
\put{$\scriptstyle{1}$}[t] at -0.35 -0.05
\put{$\scriptstyle{2}$}[tl] at -0.25 -0.33
\put{$\scriptstyle{1}$}[l] at 0.05 -0.35
\put{$\scriptstyle{2}$}[bl] at 0.31 -0.27
\put{$\scriptstyle{1}$}[b] at 0.35 0.05
\put{$\scriptstyle{2}$}[r] at 0.25  0.35
\put{$\scriptstyle{0}$}[bl] at 2.3 0.75
\put{$\scriptstyle{0}$}[br] at 1.7 0.75
\put{$\scriptstyle{0}$}[bl] at 2.75  0.3
\put{$\scriptstyle{0}$}[br] at 1.25 0.3
\put{$\scriptstyle{0}$}[tl] at 2.3 -0.75
\put{$\scriptstyle{0}$}[tr] at 1.7 -0.75
\put{$\scriptstyle{0}$}[tl] at 2.75 -0.3
\put{$\scriptstyle{0}$}[tr] at 1.25 -0.3
\put{$\scriptstyle{1}$}[r] at 1.95 0.35
\put{$\scriptstyle{2}$}[tr] at 1.7 0.23
\put{$\scriptstyle{1}$}[t] at 1.65 -0.05
\put{$\scriptstyle{2}$}[tl] at 1.75 -0.33
\put{$\scriptstyle{1}$}[l] at 2.05 -0.35
\put{$\scriptstyle{2}$}[bl] at 2.31 -0.27
\put{$\scriptstyle{1}$}[b] at 2.35 0.05
\put{$\scriptstyle{2}$}[r] at 2.25  0.35
\put{$\scriptstyle{0}$}[bl] at -1.7 0.75
\put{$\scriptstyle{0}$}[br] at -2.3 0.75
\put{$\scriptstyle{0}$}[bl] at -1.25  0.3
\put{$\scriptstyle{0}$}[br] at -2.75 0.3
\put{$\scriptstyle{0}$}[tl] at -1.7 -0.75
\put{$\scriptstyle{0}$}[tr] at -2.3 -0.75
\put{$\scriptstyle{0}$}[tl] at -1.25 -0.3
\put{$\scriptstyle{0}$}[tr] at -2.75 -0.3
\put{$\scriptstyle{1}$}[r] at -2.05 0.35
\put{$\scriptstyle{2}$}[tr] at -2.3 0.23
\put{$\scriptstyle{1}$}[t] at -2.35 -0.05
\put{$\scriptstyle{2}$}[tl] at -2.25 -0.33
\put{$\scriptstyle{1}$}[l] at -1.95 -0.35
\put{$\scriptstyle{2}$}[bl] at -1.69 -0.27
\put{$\scriptstyle{1}$}[b] at -1.65 0.05
\put{$\scriptstyle{2}$}[r] at -1.75  0.35
\put{$\scriptstyle{0}$}[bl] at 0.3 2.75
\put{$\scriptstyle{0}$}[br] at -0.3 2.75
\put{$\scriptstyle{0}$}[bl] at 0.75  2.3
\put{$\scriptstyle{0}$}[br] at -0.75 2.3
\put{$\scriptstyle{0}$}[tl] at 0.3  1.25
\put{$\scriptstyle{0}$}[tr] at -0.3  1.25
\put{$\scriptstyle{0}$}[tl] at 0.75  1.7
\put{$\scriptstyle{0}$}[tr] at -0.75  1.7
\put{$\scriptstyle{1}$}[r] at -0.05 2.35
\put{$\scriptstyle{2}$}[tr] at -0.3 2.23
\put{$\scriptstyle{1}$}[t] at -0.35  1.95
\put{$\scriptstyle{2}$}[tl] at -0.25  1.67
\put{$\scriptstyle{1}$}[l] at 0.05   1.65
\put{$\scriptstyle{2}$}[bl] at 0.31  1.73
\put{$\scriptstyle{1}$}[b] at 0.35 2.05
\put{$\scriptstyle{2}$}[r] at 0.25  2.35
\put{$\scriptstyle{0}$}[bl] at 2.3 2.75
\put{$\scriptstyle{0}$}[br] at 1.7 2.75
\put{$\scriptstyle{0}$}[bl] at 2.75  2.3
\put{$\scriptstyle{0}$}[br] at 1.25 2.3
\put{$\scriptstyle{0}$}[tl] at 2.3 1.25
\put{$\scriptstyle{0}$}[tr] at 1.7  1.25
\put{$\scriptstyle{0}$}[tl] at 2.75  1.7
\put{$\scriptstyle{0}$}[tr] at 1.25  1.7
\put{$\scriptstyle{1}$}[r] at 1.95  2.35
\put{$\scriptstyle{2}$}[tr] at 1.7  2.23
\put{$\scriptstyle{1}$}[t] at 1.65  1.95
\put{$\scriptstyle{2}$}[tl] at 1.75  1.67
\put{$\scriptstyle{1}$}[l] at 2.05   1.65
\put{$\scriptstyle{2}$}[bl] at 2.31   1.73
\put{$\scriptstyle{1}$}[b] at 2.35  2.05
\put{$\scriptstyle{2}$}[r] at 2.25   2.35
\put{$\scriptstyle{0}$}[bl] at -1.7  2.75
\put{$\scriptstyle{0}$}[br] at -2.3  2.75
\put{$\scriptstyle{0}$}[bl] at -1.25  2.3
\put{$\scriptstyle{0}$}[br] at -2.75  2.3
\put{$\scriptstyle{0}$}[tl] at -1.7  1.25
\put{$\scriptstyle{0}$}[tr] at -2.3  1.25
\put{$\scriptstyle{0}$}[tl] at -1.25  1.7
\put{$\scriptstyle{0}$}[tr] at -2.75  1.7
\put{$\scriptstyle{1}$}[r] at -2.05  2.35
\put{$\scriptstyle{2}$}[tr] at -2.3  2.23
\put{$\scriptstyle{1}$}[t] at -2.35  1.95
\put{$\scriptstyle{2}$}[tl] at -2.25  1.67
\put{$\scriptstyle{1}$}[l] at -1.95  1.65
\put{$\scriptstyle{2}$}[bl] at -1.69  1.73
\put{$\scriptstyle{1}$}[b] at -1.65  2.05
\put{$\scriptstyle{2}$}[r] at -1.75  2.35
\put{$\scriptstyle{1}$}[r] at  0.95 1.35
\put{$\scriptstyle{2}$}[tr] at  0.7 1.23
\put{$\scriptstyle{1}$}[t] at  0.65  0.95
\put{$\scriptstyle{2}$}[tl] at  0.75  0.67
\put{$\scriptstyle{1}$}[l] at 1.05  0.65
\put{$\scriptstyle{2}$}[bl] at 1.31  0.73
\put{$\scriptstyle{1}$}[b] at 1.35 1.05
\put{$\scriptstyle{2}$}[r] at 1.25  1.35
\put{$\scriptstyle{1}$}[r] at -1.05  1.35
\put{$\scriptstyle{2}$}[tr] at -1.3  1.23
\put{$\scriptstyle{1}$}[t] at -1.35  0.95
\put{$\scriptstyle{2}$}[tl] at -1.25  0.67
\put{$\scriptstyle{1}$}[l] at  -0.95  0.65
\put{$\scriptstyle{2}$}[bl] at -0.69  0.73
\put{$\scriptstyle{1}$}[b] at -0.65  1.05
\put{$\scriptstyle{2}$}[r] at -0.75   1.35
\plot -1.2 -1.2   3.2 3.2 /
\plot  1.2 -1.2  -3.2 3.2 /
\plot  0  3.2   0 -1.2 /
\plot  3.2  0  -3.2  0 /
\arrow <5pt> [.2,.67] from 0.175 0.5 to -0.175 0.5   %
\arrow <5pt> [.2,.67] from -0.5 0.175 to -0.175 0.5   
\arrow <5pt> [.2,.67] from -0.5 0.175 to -0.825 0.5   %
\arrow <5pt> [.2,.67] from -0.825 0.5 to -1.175 0.5   %
\plot -1.5 0.175  -1.375 0.315 / %
\plot -1.375 0.316  -1.42 0.35 / %
\arrow <5pt> [.2,.67] from -1.42 0.35 to -1.55 0.2   %
\arrow <5pt> [.2,.67] from -1.2 0.55 to -1.5 0.825   %
\arrow <5pt> [.2,.67] from -1.5 1.125 to -1.5 0.825   %
\arrow <5pt> [.2,.67] from -1.5 1.125 to -1.2 1.5   %
\plot -1.2 1.5  -1.05 1.5 /   %
\plot -1.05 1.5 -1.05 1.55 /  %
\arrow <5pt> [.2,.67] from -1.05 1.55 to -1.22 1.55   %
\arrow <5pt> [.2,.67] from -1.22 1.55 to -1.5 1.82   %
\arrow <5pt> [.2,.67] from -1.5 1.82 to -1.5 2.2   %
\arrow <5pt> [.2,.67] from -1.825 2.5 to -1.5 2.2   %
\setdashes
\plot  -3 3.2   -3 -1.2 /
\plot  -2 3.2   -2 -1.2 /
\plot  -1 3.2   -1 -1.2 /
\plot   1 3.2    1 -1.2 /
\plot   2 3.2    2 -1.2 /
\plot   3 3.2    3 -1.2 /
\plot  3.2 -1   -3.2 -1 /
\plot  3.2  1   -3.2  1 /
\plot  3.2  2   -3.2  2 /
\plot  3.2  3   -3.2  3 /
%
\plot  3.2 -0.8   2.8 -1.2 /
\plot  3.2 0.2    1.8 -1.2 /
\plot  3.2 1.2   0.8 -1.2 /
\plot  3.2 2.2   -0.2 -1.2 /
\plot  2.2 3.2   -2.2 -1.2 /
\plot  1.2 3.2   -3.2 -1.2 /
\plot  0.2 3.2   -3.2 -0.2 /
\plot  -0.8 3.2  -3.2 0.8 /
\plot  -1.8 3.2  -3.2 1.8 /
\plot  -2.2 3.2   2.3 -1.3 /
\plot  -1.2 3.2   3.2 -1.2 /
\plot  -0.2 3.2   3.2 -0.2 /
\plot  0.8 3.2  3.2 0.8 /
\plot  1.8 3.2  3.2 1.8 /
%
\plot  -3.2 -0.8   -2.8 -1.2 /
\plot  -3.2 0.2    -1.8 -1.2 /
\plot  -3.2 1.2   -0.8  -1.2 /
\plot  -3.2 2.2   0.2 -1.2 /
\endpicture
$$
but, by first applying relations $f_i^\mp = - f_i^\pm$ and then
working left to right applying the relations $c_i^\pm = c_i^\mp+f_i^\pm$, gives
\begin{align*}
c_1^-c_2^+c_0^+c_1^-f_0^-c_2^+c_1^-c_2^+f_1^-c_0^+c_1^+c_2^-
&=-(c_1^-c_2^+c_0^+c_1^-f_0^+c_2^+c_1^-c_2^+f_1^-c_0^+c_1^+c_2^-) \\
&=-(c_1^-(c_2^-+f_2^+)c_0^+c_1^-f_0^+c_2^+c_1^-c_2^+f_1^-c_0^+c_1^+c_2^-) \\
&=-(c_1^-(c_2^-+f_2^+)c_0^+c_1^-f_0^+c_2^+(c_1^++f_1^-)c_2^+f_1^-c_0^+c_1^+c_2^-) \\
&=-\big(c_1^-(c_2^-+f_2^+)c_0^+c_1^-f_0^+c_2^+(c_1^++f_1^-)c_2^+f_1^-c_0^+c_1^+
(c_2^++f_2^-)\big) 
\end{align*}
and every term in the expansion of this expression is an alcove walk.

\end{subsection}

\begin{subsection}{The affine Hecke algebra}

Fix an invertible element $q\in \KK$.
The \emph{affine Hecke algebra} $\tilde H$ is the quotient of the alcove walk
algebra by the relations
$$
\beginpicture
\setcoordinatesystem units <1cm,1cm>         
\setplotarea x from -0.7 to 0.7, y from -0.5 to 0.5  
\put{$i$}[b] at 0 0.6
\put{$\scriptstyle{-}$}[b] at -0.4 0.25
\put{$\scriptstyle{+}$}[b] at 0.4 0.25
\plot  0 -0.4  0 0.5 /
\arrow <5pt> [.2,.67] from -0.5 0 to 0.5 0   %
\endpicture
=
\left(\beginpicture
\setcoordinatesystem units <1cm,1cm>         
\setplotarea x from -0.7 to 0.7, y from -0.5 to 0.5  
\put{$i$}[b] at 0 0.6
\put{$\scriptstyle{-}$}[b] at -0.4 0.25
\put{$\scriptstyle{+}$}[b] at 0.4 0.25
\plot  0 -0.4  0 0.5 /
\arrow <5pt> [.2,.67] from 0.5 0 to -0.5 0   %
\endpicture
\right)^{-1},
\qquad
\beginpicture
\setcoordinatesystem units <1cm,1cm>         
\setplotarea x from -0.7 to 0.7, y from -0.5 to 0.5  
\put{$i$}[b] at 0 0.7
\put{$\scriptstyle{-}$}[b] at -0.4 0.35
\put{$\scriptstyle{+}$}[b] at 0.4 0.35
\plot  0 -0.4  0 0.6 /
\plot -0.5 0  -0.05 0 /
\plot -0.05 0 -0.05 0.1 /
\arrow <5pt> [.2,.67] from -0.05 0.1 to -0.5 0.1   %
\endpicture 
=-(q-q^{-1}),
\qquad
\beginpicture
\setcoordinatesystem units <1cm,1cm>         
\setplotarea x from -0.7 to 0.7, y from -0.5 to 0.5  
\put{$i$}[b] at 0 0.7
\put{$\scriptstyle{-}$}[b] at -0.4 0.35
\put{$\scriptstyle{+}$}[b] at 0.4 0.35
\plot  0 -0.4  0 0.6 /
\plot 0.5 0  0.05 0 /
\arrow <5pt> [.2,.67] from 0.05 0.1 to 0.5 0.1   %
\plot 0.05 0 0.05 0.1 /
\endpicture 
=(q-q^{-1}),
$$
and
\begin{equation}
p=p'
\qquad\hbox{if $p$ and $p'$ are nonfolded walks with $\mathrm{end}(p)=\mathrm{end}(p')$,}
\end{equation}
where $\mathrm{end}(p)$ is the final alcove of $p$.  
Conceptually, the affine Hecke algebra only remembers the ending 
alcove of a walk (and some information about the folds)
and forgets how it got to its destination.

For $w\in W$ and $\lambda\in P$ define elements
\begin{align*}
T_{w^{-1}}^{-1} &= (\hbox{image in $\tilde H$ of a minimal length alcove walk from $A$ to $wA$}), \\
X^\lambda &= (\hbox{image in $\tilde H$ of a minimal length alcove walk from $A$ to $t_\lambda A$}).
\end{align*}
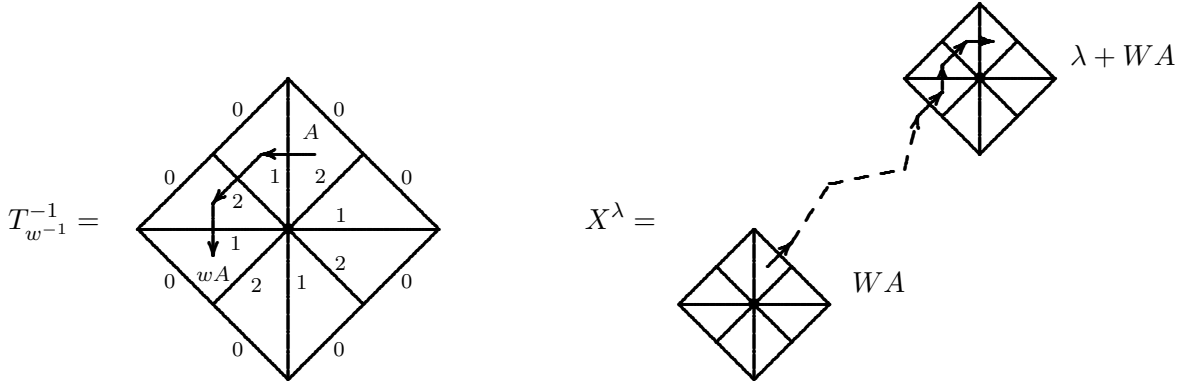
$$
T_{w^{-1}}^{-1} = 
\beginpicture
\setcoordinatesystem units <2cm,2cm>         
\setplotarea x from -1.2 to 1.2, y from -1.2 to 1.2    
\put{$\scriptstyle{0}$}[bl] at 0.3 0.75
\put{$\scriptstyle{0}$}[br] at -0.3 0.75
\put{$\scriptstyle{0}$}[bl] at 0.75  0.3
\put{$\scriptstyle{0}$}[br] at -0.75 0.3
\put{$\scriptstyle{0}$}[tl] at 0.3 -0.75
\put{$\scriptstyle{0}$}[tr] at -0.3 -0.75
\put{$\scriptstyle{0}$}[tl] at 0.75 -0.3
\put{$\scriptstyle{0}$}[tr] at -0.75 -0.3
\put{$\scriptstyle{1}$}[r] at -0.05 0.35
\put{$\scriptstyle{2}$}[tr] at -0.3 0.23
\put{$\scriptstyle{1}$}[t] at -0.35 -0.05
\put{$\scriptstyle{2}$}[tl] at -0.25 -0.33
\put{$\scriptstyle{1}$}[l] at 0.05 -0.35
\put{$\scriptstyle{2}$}[bl] at 0.31 -0.27
\put{$\scriptstyle{1}$}[b] at 0.35 0.05
\put{$\scriptstyle{2}$}[r] at 0.25  0.35
%
\plot  0 -1  0 1 /
\plot  -1 0  1 0 /
\plot  0 -1  1 0 /
\plot   0 -1  -1 0 /
\plot  -1 0  0 1 /
\plot  1 0  0 1 /
\plot  -0.5 0.5  0.5 -0.5 /
\plot  -0.5 -0.5  0.5 0.5 /
\arrow <5pt> [.2,.67] from 0.175 0.5 to  -0.175 0.5   
\arrow <5pt> [.2,.67] from -0.175 0.5 to -0.5 0.175   
\arrow <5pt> [.2,.67] from -0.5 0.175 to -0.5 -0.175   
\put{$\bullet$} at 0  0
\put{$\scriptstyle{A}$}[b] at 0.15  0.6
\put{$\scriptstyle{wA}$}[t] at -0.5 -0.25
\endpicture
\qquad\qquad
X^\lambda = 
\beginpicture
\setcoordinatesystem units <1cm,1cm>         
\setplotarea x from -2.2 to 3.2, y from -2.2 to 3.2    
\plot  -1 -2  -1 0 /
\plot  -2 -1  0 -1 /
\plot  -1 -2  0 -1 /
\plot   -1 -2  -2 -1 /
\plot  -2 -1  -1 0 /
\plot  0 -1  -1 0 /
\plot  -1.5 -0.5  -0.5 -1.5 /
\plot  -1.5 -1.5  -0.5 -0.5 /
\plot  2 1  2 3 /
\plot  2 1  3 2 /
\plot   2 3  3 2 /
\plot  1 2  3 2 /
\plot  1 2  2 3 /
\plot  1 2  2 1 /
\plot  1.5 2.5  2.5 1.5 /
\plot  1.5 1.5  2.5 2.5 /
\arrow <5pt> [.2,.67] from -0.825 -0.5 to  -0.5 -0.175   
\arrow <5pt> [.2,.67] from 1.175 1.5 to  1.5 1.825   
\arrow <5pt> [.2,.67] from 1.825 2.5 to 2.175 2.5  
\arrow <5pt> [.2,.67] from 1.5 2.175 to 1.825 2.5   
\arrow <5pt> [.2,.67] from 1.5 1.825 to 1.5 2.175   
\put{$\bullet$} at -1  -1
\put{$\bullet$} at 2  2
\put{$WA$}[l] at 0.3  -0.7
\put{$\lambda+WA$}[l] at 3.2 2.3
\setdashes
\plot -0.5 -0.175  0 0.6 /
\plot 0 0.6  1 0.8 /
\arrow <5pt> [.2,.67] from 1 0.8 to 1.175 1.5   
\endpicture
$$
The following proposition shows that the alcove walk definition of the affine Hecke algebra
coincides with the standard definition by generators and relations (see [IM] and [Lu]).
A consequence of the proposition is that 
\begin{equation}
\begin{array}{lll}
\hbox{the \emph{finite Hecke algebra}},\qquad
&H = \hbox{span}\{ T_{w^{-1}}^{-1}\ |\ w\in W\},
&\hbox{and} \\
\hbox{the Laurent polynomial ring,}\qquad
&\KK[P] = \hbox{span}\{X^\lambda\ |\ \lambda\in P\},
\end{array}
\end{equation}
are subalgebras of $\tilde H$.

\begin{prop}  Let $g\in \Omega$, $\lambda,\mu\in P$, $w\in W$ and $1\le i\le n$.
Let $\varphi$ be the element of $R^+$ such that $H_{\alpha_0} = H_{\varphi,1}$ is
the wall of $A$ which is not a wall of $C$ and let $s_\varphi$ be the reflection in $H_\varphi$.
Let $w_0$ be the longest element of $W$.  
The following identities hold in $\tilde H$.
\item[(a)]  $X^\lambda X^\mu = X^{\lambda+\mu}=X^\mu X^\lambda$.
\item[(b)]  $\displaystyle{T_{s_i}T_w = \begin{cases}
T_{s_iw}, &\hbox{if $\ell(s_iw)>\ell(w)$,} \\
T_{s_iw} + (q-q^{-1})T_w, &\hbox{if $\ell(s_iw)<\ell(w)$.}
\end{cases}
}$
\item[(c)] If $\langle\lambda,\alpha_i^\vee\rangle = 0$ then $T_{s_i}X^\lambda = X^\lambda T_{s_i}$.
\item[(d)] If $\langle\lambda,\alpha_i^\vee\rangle =1$ then 
$T_{s_i}X^{s_i\lambda}T_{s_i} = X^\lambda$.
\item[(e)] 
$\displaystyle{
T_{s_i}X^\lambda = X^{s_i\lambda}T_{s_i} + (q-q^{-1})\frac{X^\lambda-X^{s_i\lambda}}
{1-X^{-\alpha_i}}. }$
\item[(f)]  $T_{s_0}T_{s_{\varphi}} = X^\varphi$.
\item[(g)] $X^{\omega_i} = gT_{w_0w_i}$, where the action of $g$ on A sends the origin to $\omega_i$
and $w_i$ is the longest element of the stabilizer $W_{\omega_i}$ of $\omega_i$ in $W$.
\end{prop}
\begin{proof}  Use notations for alcove walks as in (3.4).

\smallskip\noindent
(a)  If $p_\lambda$ is a minimal length walk from $A$ to $t_\lambda A$ and 
$p_\mu$ is a minimal length walk from from $A$ to $t_\mu A$ then
$$
\hbox{$p_\lambda p_\mu$
and $p_\mu p_\lambda$ are both nonfolded walks from $A$ to $t_{\lambda+\mu}A$.}
$$
Thus the images of $p_\lambda p_\mu$ and $p_\mu p_\lambda$
are equal in $\tilde H$.

\smallskip\noindent
(b) If $\ell(ws_i)>\ell(w)$ and $p_{w}$ is a minimal length walk from $A$ to $wA$ then
$$\hbox{$p_{ws_i} = p_{w}c_i^-$ is a minimal length walk from $A$ to $ws_iA$.}
$$
and so 
$T_{s_iw^{-1}}^{-1}
=T_{{ws_i}^{-1}}^{-1} = T_{w^{-1}}^{-1}T_{s_i}^{-1}=(T_{s_i}T_{w^{-1}})^{-1}$ 
in $\tilde H$.  Taking inverses gives the first
result, and the second follows by switching $w$ and $ws_i$ and using the relation
$T_{s_i}^{-1} = T_{s_i} -(q-q^{-1})$ which follows from (3.2) and (3.5).

\smallskip\noindent
(c)  Let $p_\lambda$ be a minimal length alcove walk from $A$ to $t_\lambda A$.
If $\langle \lambda,\alpha_i^\vee\rangle = 0$ then $H_{\alpha_i}$ is a 
wall of $t_\lambda A$ and $s_i\lambda = \lambda$ and
$$
\hbox{$c_i^- p_\lambda c_i^+$ is a nonfolded walk from $A$ to $t_\lambda A$.}
$$
Thus $T_{s_i}^{-1} X^\lambda T_{s_i} = X^\lambda = X^{s_i\lambda}$ in $\tilde H$.
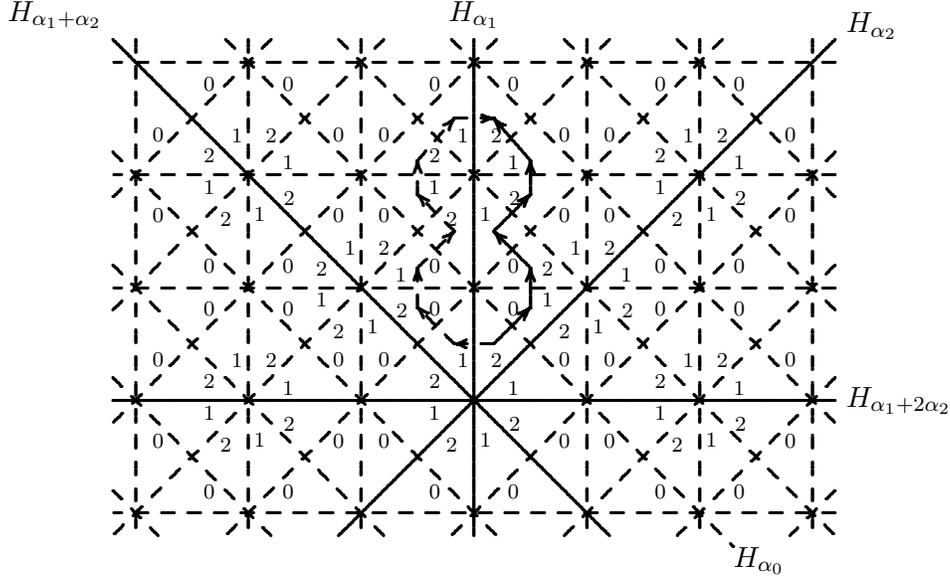
$$
\beginpicture
\setcoordinatesystem units <1.5cm,1.5cm>         
\setplotarea x from -3.5 to 3.5, y from -1.5 to 4    
\put{$H_{\alpha_1}$}[b] at 0 3.3
\put{$H_{\alpha_2}$}[bl] at 3.3 3.2
\put{$H_{\alpha_1+\alpha_2}$}[br] at -3.3 3.3
\put{$H_{\alpha_0}$}[tl] at 2.3 -1.3 
\put{$H_{\alpha_1+2\alpha_2}$}[l] at 3.3 0
\put{$\scriptstyle{0}$}[bl] at 0.3 0.75
\put{$\scriptstyle{0}$}[br] at -0.3 0.75
\put{$\scriptstyle{0}$}[bl] at 0.75  0.3
\put{$\scriptstyle{0}$}[br] at -0.75 0.3
\put{$\scriptstyle{0}$}[tl] at 0.3 -0.75
\put{$\scriptstyle{0}$}[tr] at -0.3 -0.75
\put{$\scriptstyle{0}$}[tl] at 0.75 -0.3
\put{$\scriptstyle{0}$}[tr] at -0.75 -0.3
\put{$\scriptstyle{1}$}[r] at -0.05 0.35
\put{$\scriptstyle{2}$}[tr] at -0.3 0.23
\put{$\scriptstyle{1}$}[t] at -0.35 -0.05
\put{$\scriptstyle{2}$}[tl] at -0.25 -0.33
\put{$\scriptstyle{1}$}[l] at 0.05 -0.35
\put{$\scriptstyle{2}$}[bl] at 0.31 -0.27
\put{$\scriptstyle{1}$}[b] at 0.35 0.05
\put{$\scriptstyle{2}$}[r] at 0.25  0.35
\put{$\scriptstyle{0}$}[bl] at 2.3 0.75
\put{$\scriptstyle{0}$}[br] at 1.7 0.75
\put{$\scriptstyle{0}$}[bl] at 2.75  0.3
\put{$\scriptstyle{0}$}[br] at 1.25 0.3
\put{$\scriptstyle{0}$}[tl] at 2.3 -0.75
\put{$\scriptstyle{0}$}[tr] at 1.7 -0.75
\put{$\scriptstyle{0}$}[tl] at 2.75 -0.3
\put{$\scriptstyle{0}$}[tr] at 1.25 -0.3
\put{$\scriptstyle{1}$}[r] at 1.95 0.35
\put{$\scriptstyle{2}$}[tr] at 1.7 0.23
\put{$\scriptstyle{1}$}[t] at 1.65 -0.05
\put{$\scriptstyle{2}$}[tl] at 1.75 -0.33
\put{$\scriptstyle{1}$}[l] at 2.05 -0.35
\put{$\scriptstyle{2}$}[bl] at 2.31 -0.27
\put{$\scriptstyle{1}$}[b] at 2.35 0.05
\put{$\scriptstyle{2}$}[r] at 2.25  0.35
\put{$\scriptstyle{0}$}[bl] at -1.7 0.75
\put{$\scriptstyle{0}$}[br] at -2.3 0.75
\put{$\scriptstyle{0}$}[bl] at -1.25  0.3
\put{$\scriptstyle{0}$}[br] at -2.75 0.3
\put{$\scriptstyle{0}$}[tl] at -1.7 -0.75
\put{$\scriptstyle{0}$}[tr] at -2.3 -0.75
\put{$\scriptstyle{0}$}[tl] at -1.25 -0.3
\put{$\scriptstyle{0}$}[tr] at -2.75 -0.3
\put{$\scriptstyle{1}$}[r] at -2.05 0.35
\put{$\scriptstyle{2}$}[tr] at -2.3 0.23
\put{$\scriptstyle{1}$}[t] at -2.35 -0.05
\put{$\scriptstyle{2}$}[tl] at -2.25 -0.33
\put{$\scriptstyle{1}$}[l] at -1.95 -0.35
\put{$\scriptstyle{2}$}[bl] at -1.69 -0.27
\put{$\scriptstyle{1}$}[b] at -1.65 0.05
\put{$\scriptstyle{2}$}[r] at -1.75  0.35
\put{$\scriptstyle{0}$}[bl] at 0.3 2.75
\put{$\scriptstyle{0}$}[br] at -0.3 2.75
\put{$\scriptstyle{0}$}[bl] at 0.75  2.3
\put{$\scriptstyle{0}$}[br] at -0.75 2.3
\put{$\scriptstyle{0}$}[tl] at 0.3  1.25
\put{$\scriptstyle{0}$}[tr] at -0.3  1.25
\put{$\scriptstyle{0}$}[tl] at 0.75  1.7
\put{$\scriptstyle{0}$}[tr] at -0.75  1.7
\put{$\scriptstyle{1}$}[r] at -0.05 2.35
\put{$\scriptstyle{2}$}[tr] at -0.3 2.23
\put{$\scriptstyle{1}$}[t] at -0.35  1.95
\put{$\scriptstyle{2}$}[tl] at -0.25  1.67
\put{$\scriptstyle{1}$}[l] at 0.05   1.65
\put{$\scriptstyle{2}$}[bl] at 0.31  1.73
\put{$\scriptstyle{1}$}[b] at 0.35 2.05
\put{$\scriptstyle{2}$}[r] at 0.25  2.35
\put{$\scriptstyle{0}$}[bl] at 2.3 2.75
\put{$\scriptstyle{0}$}[br] at 1.7 2.75
\put{$\scriptstyle{0}$}[bl] at 2.75  2.3
\put{$\scriptstyle{0}$}[br] at 1.25 2.3
\put{$\scriptstyle{0}$}[tl] at 2.3 1.25
\put{$\scriptstyle{0}$}[tr] at 1.7  1.25
\put{$\scriptstyle{0}$}[tl] at 2.75  1.7
\put{$\scriptstyle{0}$}[tr] at 1.25  1.7
\put{$\scriptstyle{1}$}[r] at 1.95  2.35
\put{$\scriptstyle{2}$}[tr] at 1.7  2.23
\put{$\scriptstyle{1}$}[t] at 1.65  1.95
\put{$\scriptstyle{2}$}[tl] at 1.75  1.67
\put{$\scriptstyle{1}$}[l] at 2.05   1.65
\put{$\scriptstyle{2}$}[bl] at 2.31   1.73
\put{$\scriptstyle{1}$}[b] at 2.35  2.05
\put{$\scriptstyle{2}$}[r] at 2.25   2.35
\put{$\scriptstyle{0}$}[bl] at -1.7  2.75
\put{$\scriptstyle{0}$}[br] at -2.3  2.75
\put{$\scriptstyle{0}$}[bl] at -1.25  2.3
\put{$\scriptstyle{0}$}[br] at -2.75  2.3
\put{$\scriptstyle{0}$}[tl] at -1.7  1.25
\put{$\scriptstyle{0}$}[tr] at -2.3  1.25
\put{$\scriptstyle{0}$}[tl] at -1.25  1.7
\put{$\scriptstyle{0}$}[tr] at -2.75  1.7
\put{$\scriptstyle{1}$}[r] at -2.05  2.35
\put{$\scriptstyle{2}$}[tr] at -2.3  2.23
\put{$\scriptstyle{1}$}[t] at -2.35  1.95
\put{$\scriptstyle{2}$}[tl] at -2.25  1.67
\put{$\scriptstyle{1}$}[l] at -1.95  1.65
\put{$\scriptstyle{2}$}[bl] at -1.69  1.73
\put{$\scriptstyle{1}$}[b] at -1.65  2.05
\put{$\scriptstyle{2}$}[r] at -1.75  2.35
\put{$\scriptstyle{1}$}[r] at  0.95 1.35
\put{$\scriptstyle{2}$}[tr] at  0.7 1.23
\put{$\scriptstyle{1}$}[t] at  0.65  0.95
\put{$\scriptstyle{2}$}[tl] at  0.75  0.67
\put{$\scriptstyle{1}$}[l] at 1.05  0.65
\put{$\scriptstyle{2}$}[bl] at 1.31  0.73
\put{$\scriptstyle{1}$}[b] at 1.35 1.05
\put{$\scriptstyle{2}$}[r] at 1.25  1.35
\put{$\scriptstyle{1}$}[r] at -1.05  1.35
\put{$\scriptstyle{2}$}[tr] at -1.3  1.23
\put{$\scriptstyle{1}$}[t] at -1.35  0.95
\put{$\scriptstyle{2}$}[tl] at -1.25  0.67
\put{$\scriptstyle{1}$}[l] at  -0.95  0.65
\put{$\scriptstyle{2}$}[bl] at -0.69  0.73
\put{$\scriptstyle{1}$}[b] at -0.65  1.05
\put{$\scriptstyle{2}$}[r] at -0.75   1.35
\plot -1.2 -1.2   3.2 3.2 /
\plot  1.2 -1.2  -3.2 3.2 /
\plot  0  3.2   0 -1.2 /
\plot  3.2  0  -3.2  0 /
\arrow <5pt> [.2,.67] from 0.175 0.5 to  0.5 0.825   %
\arrow <5pt> [.2,.67] from 0.5 0.825 to 0.5 1.175   %
\arrow <5pt> [.2,.67] from 0.5 1.175 to 0.175 1.5   %
\arrow <5pt> [.2,.67] from 0.175 1.5 to 0.5 1.825   %
\arrow <5pt> [.2,.67] from 0.5 1.825 to 0.5 2.125   %
\arrow <5pt> [.2,.67] from 0.5 2.125 to 0.175 2.5   %
\setdashes
\arrow <5pt> [.2,.67] from 0.175 0.5 to -0.175 0.5   %
\arrow <5pt> [.2,.67] from -0.175 0.5 to -0.5 0.825   %
\arrow <5pt> [.2,.67] from -0.5 0.825 to -0.5 1.175   %
\arrow <5pt> [.2,.67] from -0.5 1.175 to -0.175 1.5   %
\arrow <5pt> [.2,.67] from -0.175 1.5 to -0.5 1.825   %
\arrow <5pt> [.2,.67] from -0.5 1.825 to -0.5 2.125   %
\arrow <5pt> [.2,.67] from -0.5 2.125 to -0.175 2.5   %
\arrow <5pt> [.2,.67] from -0.175 2.5 to 0.175 2.5   %
\plot  -3 3.2   -3 -1.2 /
\plot  -2 3.2   -2 -1.2 /
\plot  -1 3.2   -1 -1.2 /
\plot   1 3.2    1 -1.2 /
\plot   2 3.2    2 -1.2 /
\plot   3 3.2    3 -1.2 /
\plot  3.2 -1   -3.2 -1 /
\plot  3.2  1   -3.2  1 /
\plot  3.2  2   -3.2  2 /
\plot  3.2  3   -3.2  3 /
%
\plot  3.2 -0.8   2.8 -1.2 /
\plot  3.2 0.2    1.8 -1.2 /
\plot  3.2 1.2   0.8 -1.2 /
\plot  3.2 2.2   -0.2 -1.2 /
\plot  2.2 3.2   -2.2 -1.2 /
\plot  1.2 3.2   -3.2 -1.2 /
\plot  0.2 3.2   -3.2 -0.2 /
\plot  -0.8 3.2  -3.2 0.8 /
\plot  -1.8 3.2  -3.2 1.8 /
\plot  -2.2 3.2   2.3 -1.3 /
\plot  -1.2 3.2   3.2 -1.2 /
\plot  -0.2 3.2   3.2 -0.2 /
\plot  0.8 3.2  3.2 0.8 /
\plot  1.8 3.2  3.2 1.8 /
%
\plot  -3.2 -0.8   -2.8 -1.2 /
\plot  -3.2 0.2    -1.8 -1.2 /
\plot  -3.2 1.2   -0.8  -1.2 /
\plot  -3.2 2.2   0.2 -1.2 /
\endpicture
$$

\smallskip\noindent
(d) Let $p_\lambda$ be a minimal length walk from $A$ to $t_\lambda A$.  
If $\langle \lambda,\alpha_i^\vee\rangle = 1$ then
there is a minimal length walk from $A$ to $t_\lambda A$ of the form
$p_\lambda = p_{t_\lambda s_i}c_i^+$ where $p_{t_\lambda s_i}$ is minimal length walk
from $A$ to $t_\lambda s_i A$.  Then 
$$\hbox{$c_i^-p_{t_\lambda s_i}$ is a minimal length walk from $A$ to $t_{s_i\lambda}A$.}$$
Thus  $T_{s_i}^{-1}(X^\lambda T_{s_i^{-1}}) = X^{s_i\lambda}$ in $\tilde H$.
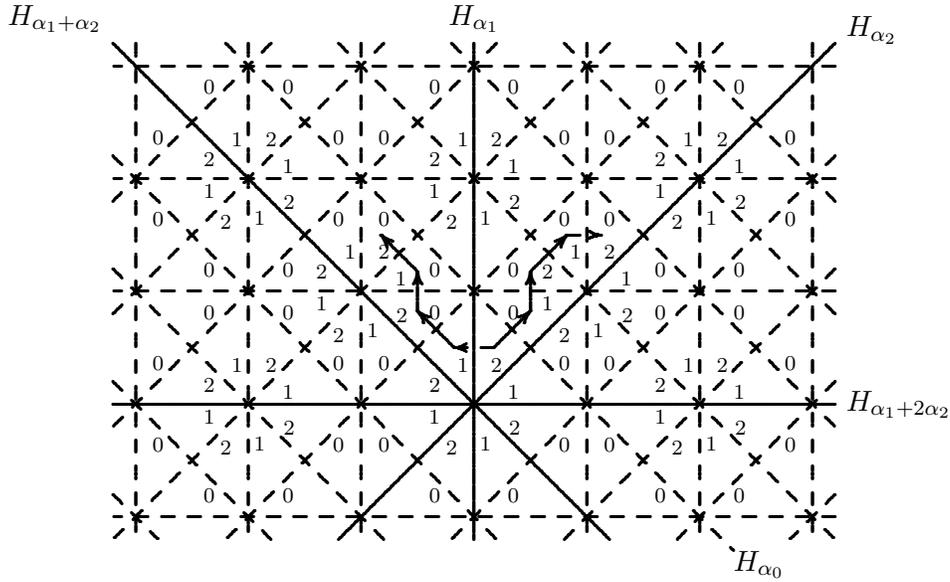
$$
\beginpicture
\setcoordinatesystem units <1.5cm,1.5cm>         
\setplotarea x from -3.5 to 3.5, y from -1.5 to 4    
\put{$H_{\alpha_1}$}[b] at 0 3.3
\put{$H_{\alpha_2}$}[bl] at 3.3 3.2
\put{$H_{\alpha_1+\alpha_2}$}[br] at -3.3 3.3
\put{$H_{\alpha_0}$}[tl] at 2.3 -1.3 
\put{$H_{\alpha_1+2\alpha_2}$}[l] at 3.3 0
\put{$\scriptstyle{0}$}[bl] at 0.3 0.75
\put{$\scriptstyle{0}$}[br] at -0.3 0.75
\put{$\scriptstyle{0}$}[bl] at 0.75  0.3
\put{$\scriptstyle{0}$}[br] at -0.75 0.3
\put{$\scriptstyle{0}$}[tl] at 0.3 -0.75
\put{$\scriptstyle{0}$}[tr] at -0.3 -0.75
\put{$\scriptstyle{0}$}[tl] at 0.75 -0.3
\put{$\scriptstyle{0}$}[tr] at -0.75 -0.3
\put{$\scriptstyle{1}$}[r] at -0.05 0.35
\put{$\scriptstyle{2}$}[tr] at -0.3 0.23
\put{$\scriptstyle{1}$}[t] at -0.35 -0.05
\put{$\scriptstyle{2}$}[tl] at -0.25 -0.33
\put{$\scriptstyle{1}$}[l] at 0.05 -0.35
\put{$\scriptstyle{2}$}[bl] at 0.31 -0.27
\put{$\scriptstyle{1}$}[b] at 0.35 0.05
\put{$\scriptstyle{2}$}[r] at 0.25  0.35
\put{$\scriptstyle{0}$}[bl] at 2.3 0.75
\put{$\scriptstyle{0}$}[br] at 1.7 0.75
\put{$\scriptstyle{0}$}[bl] at 2.75  0.3
\put{$\scriptstyle{0}$}[br] at 1.25 0.3
\put{$\scriptstyle{0}$}[tl] at 2.3 -0.75
\put{$\scriptstyle{0}$}[tr] at 1.7 -0.75
\put{$\scriptstyle{0}$}[tl] at 2.75 -0.3
\put{$\scriptstyle{0}$}[tr] at 1.25 -0.3
\put{$\scriptstyle{1}$}[r] at 1.95 0.35
\put{$\scriptstyle{2}$}[tr] at 1.7 0.23
\put{$\scriptstyle{1}$}[t] at 1.65 -0.05
\put{$\scriptstyle{2}$}[tl] at 1.75 -0.33
\put{$\scriptstyle{1}$}[l] at 2.05 -0.35
\put{$\scriptstyle{2}$}[bl] at 2.31 -0.27
\put{$\scriptstyle{1}$}[b] at 2.35 0.05
\put{$\scriptstyle{2}$}[r] at 2.25  0.35
\put{$\scriptstyle{0}$}[bl] at -1.7 0.75
\put{$\scriptstyle{0}$}[br] at -2.3 0.75
\put{$\scriptstyle{0}$}[bl] at -1.25  0.3
\put{$\scriptstyle{0}$}[br] at -2.75 0.3
\put{$\scriptstyle{0}$}[tl] at -1.7 -0.75
\put{$\scriptstyle{0}$}[tr] at -2.3 -0.75
\put{$\scriptstyle{0}$}[tl] at -1.25 -0.3
\put{$\scriptstyle{0}$}[tr] at -2.75 -0.3
\put{$\scriptstyle{1}$}[r] at -2.05 0.35
\put{$\scriptstyle{2}$}[tr] at -2.3 0.23
\put{$\scriptstyle{1}$}[t] at -2.35 -0.05
\put{$\scriptstyle{2}$}[tl] at -2.25 -0.33
\put{$\scriptstyle{1}$}[l] at -1.95 -0.35
\put{$\scriptstyle{2}$}[bl] at -1.69 -0.27
\put{$\scriptstyle{1}$}[b] at -1.65 0.05
\put{$\scriptstyle{2}$}[r] at -1.75  0.35
\put{$\scriptstyle{0}$}[bl] at 0.3 2.75
\put{$\scriptstyle{0}$}[br] at -0.3 2.75
\put{$\scriptstyle{0}$}[bl] at 0.75  2.3
\put{$\scriptstyle{0}$}[br] at -0.75 2.3
\put{$\scriptstyle{0}$}[tl] at 0.3  1.25
\put{$\scriptstyle{0}$}[tr] at -0.3  1.25
\put{$\scriptstyle{0}$}[tl] at 0.75  1.7
\put{$\scriptstyle{0}$}[tr] at -0.75  1.7
\put{$\scriptstyle{1}$}[r] at -0.05 2.35
\put{$\scriptstyle{2}$}[tr] at -0.3 2.23
\put{$\scriptstyle{1}$}[t] at -0.35  1.95
\put{$\scriptstyle{2}$}[tl] at -0.25  1.67
\put{$\scriptstyle{1}$}[l] at 0.05   1.65
\put{$\scriptstyle{2}$}[bl] at 0.31  1.73
\put{$\scriptstyle{1}$}[b] at 0.35 2.05
\put{$\scriptstyle{2}$}[r] at 0.25  2.35
\put{$\scriptstyle{0}$}[bl] at 2.3 2.75
\put{$\scriptstyle{0}$}[br] at 1.7 2.75
\put{$\scriptstyle{0}$}[bl] at 2.75  2.3
\put{$\scriptstyle{0}$}[br] at 1.25 2.3
\put{$\scriptstyle{0}$}[tl] at 2.3 1.25
\put{$\scriptstyle{0}$}[tr] at 1.7  1.25
\put{$\scriptstyle{0}$}[tl] at 2.75  1.7
\put{$\scriptstyle{0}$}[tr] at 1.25  1.7
\put{$\scriptstyle{1}$}[r] at 1.95  2.35
\put{$\scriptstyle{2}$}[tr] at 1.7  2.23
\put{$\scriptstyle{1}$}[t] at 1.65  1.95
\put{$\scriptstyle{2}$}[tl] at 1.75  1.67
\put{$\scriptstyle{1}$}[l] at 2.05   1.65
\put{$\scriptstyle{2}$}[bl] at 2.31   1.73
\put{$\scriptstyle{1}$}[b] at 2.35  2.05
\put{$\scriptstyle{2}$}[r] at 2.25   2.35
\put{$\scriptstyle{0}$}[bl] at -1.7  2.75
\put{$\scriptstyle{0}$}[br] at -2.3  2.75
\put{$\scriptstyle{0}$}[bl] at -1.25  2.3
\put{$\scriptstyle{0}$}[br] at -2.75  2.3
\put{$\scriptstyle{0}$}[tl] at -1.7  1.25
\put{$\scriptstyle{0}$}[tr] at -2.3  1.25
\put{$\scriptstyle{0}$}[tl] at -1.25  1.7
\put{$\scriptstyle{0}$}[tr] at -2.75  1.7
\put{$\scriptstyle{1}$}[r] at -2.05  2.35
\put{$\scriptstyle{2}$}[tr] at -2.3  2.23
\put{$\scriptstyle{1}$}[t] at -2.35  1.95
\put{$\scriptstyle{2}$}[tl] at -2.25  1.67
\put{$\scriptstyle{1}$}[l] at -1.95  1.65
\put{$\scriptstyle{2}$}[bl] at -1.69  1.73
\put{$\scriptstyle{1}$}[b] at -1.65  2.05
\put{$\scriptstyle{2}$}[r] at -1.75  2.35
\put{$\scriptstyle{1}$}[r] at  0.95 1.35
\put{$\scriptstyle{2}$}[tr] at  0.7 1.23
\put{$\scriptstyle{1}$}[t] at  0.65  0.95
\put{$\scriptstyle{2}$}[tl] at  0.75  0.67
\put{$\scriptstyle{1}$}[l] at 1.05  0.65
\put{$\scriptstyle{2}$}[bl] at 1.31  0.73
\put{$\scriptstyle{1}$}[b] at 1.35 1.05
\put{$\scriptstyle{2}$}[r] at 1.25  1.35
\put{$\scriptstyle{1}$}[r] at -1.05  1.35
\put{$\scriptstyle{2}$}[tr] at -1.3  1.23
\put{$\scriptstyle{1}$}[t] at -1.35  0.95
\put{$\scriptstyle{2}$}[tl] at -1.25  0.67
\put{$\scriptstyle{1}$}[l] at  -0.95  0.65
\put{$\scriptstyle{2}$}[bl] at -0.69  0.73
\put{$\scriptstyle{1}$}[b] at -0.65  1.05
\put{$\scriptstyle{2}$}[r] at -0.75   1.35
\plot -1.2 -1.2   3.2 3.2 /
\plot  1.2 -1.2  -3.2 3.2 /
\plot  0  3.2   0 -1.2 /
\plot  3.2  0  -3.2  0 /
\arrow <5pt> [.2,.67] from 0.175 0.5 to  0.5 0.825   
\arrow <5pt> [.2,.67] from 0.5 0.825 to 0.5 1.175   
\arrow <5pt> [.2,.67] from 0.5 1.175 to 0.825 1.5   
\arrow <5pt> [.2,.67] from -0.175 0.5 to -0.5 0.825   
\arrow <5pt> [.2,.67] from -0.5 0.825 to -0.5 1.175   
\arrow <5pt> [.2,.67] from -0.5 1.175 to -0.825 1.5   
\setdashes
\plot  -3 3.2   -3 -1.2 /
\plot  -2 3.2   -2 -1.2 /
\plot  -1 3.2   -1 -1.2 /
\plot   1 3.2    1 -1.2 /
\plot   2 3.2    2 -1.2 /
\plot   3 3.2    3 -1.2 /
\plot  3.2 -1   -3.2 -1 /
\plot  3.2  1   -3.2  1 /
\plot  3.2  2   -3.2  2 /
\plot  3.2  3   -3.2  3 /
%
\plot  3.2 -0.8   2.8 -1.2 /
\plot  3.2 0.2    1.8 -1.2 /
\plot  3.2 1.2   0.8 -1.2 /
\plot  3.2 2.2   -0.2 -1.2 /
\plot  2.2 3.2   -2.2 -1.2 /
\plot  1.2 3.2   -3.2 -1.2 /
\plot  0.2 3.2   -3.2 -0.2 /
\plot  -0.8 3.2  -3.2 0.8 /
\plot  -1.8 3.2  -3.2 1.8 /
\plot  -2.2 3.2   2.3 -1.3 /
\plot  -1.2 3.2   3.2 -1.2 /
\plot  -0.2 3.2   3.2 -0.2 /
\plot  0.8 3.2  3.2 0.8 /
\plot  1.8 3.2  3.2 1.8 /
%
\plot  -3.2 -0.8   -2.8 -1.2 /
\plot  -3.2 0.2    -1.8 -1.2 /
\plot  -3.2 1.2   -0.8  -1.2 /
\plot  -3.2 2.2   0.2 -1.2 /
%
\arrow <5pt> [.2,.67] from 0.825 1.5 to 1.125 1.5   
\arrow <5pt> [.2,.67] from 0.175 0.5 to -0.175 0.5   

\endpicture
$$

\smallskip\noindent
(e) Note that (c) and (d) are special cases of (e).
If the statement of (e) holds for $\lambda$ then, by multiplying
on the left by $X^{-s_i\lambda}$ and on the right by $X^-\lambda$,
it holds for $-\lambda$,
If the statement (e) holds for $\lambda$ and $\mu$ then it holds for $\lambda+\mu$ since
\begin{align*}
T_{s_i}X^\lambda X^{\mu}
&= \left(X^{s_i\lambda}T_{s_i} + (q-q^{-1})\frac{X^\lambda-X^{s_i\lambda}}
{1-X^{-\alpha_i}} \right)X^{\mu} \\
&= X^{s_i\lambda}\left(X^{s_i\mu}T_{s_i} + (q-q^{-1})\frac{X^\mu-X^{s_i\mu}}
{1-X^{-\alpha_i}} \right)+ (q-q^{-1})\left(\frac{X^\lambda-X^{s_i\lambda}}
{1-X^{-\alpha_i}} \right)X^{\mu} \\
&= X^{s_i(\lambda+\mu)}T_{s_i} + (q-q^{-1})\frac{X^{\lambda+\mu}-X^{s_i(\lambda+\mu)} }
{1-X^{-\alpha_i}}. 
\end{align*}
Thus, to prove (e) it is sufficient to verify (c) and (d), which has already been done.

\smallskip\noindent
(f)    Let $p_{s_\varphi}$
be a minimal length walk from $s_\varphi A$ to $A$, then 
$$\hbox{$p_{\varphi} = c_0^+p_{s_\varphi}$ is a minimal
length walk from $A$ to $t_\varphi A$.}$$
Thus $T_0T_{s_\varphi} = X^{\varphi}$ in $\tilde H$.

\smallskip\noindent
(g)  If $p_{w_0w_i}$ is a minimal length 
walk from  $w_iw_0A$ to $A$ then 
$$\hbox{$p_{\omega_i}=gp_{w_0w_i}$ is a minimal length walk from
$A$ to $t_{\omega_i}A$. }
$$
Thus $X^{\omega_i} = gT_{w_0w_i}$ in $\tilde H$.  For example, in type $C_2$,
$w_0=s_2s_1s_2s_1$ and
there is one element $g$ in $\Omega$ such that $g\ne 1$ for which $g\omega_2=0$ and
$w_2=s_1$ so that $w_0w_2=s_2s_1s_2$.
\end{proof}

The sets
\begin{equation}
\{T_{w^{-1}}^{-1}X^\lambda\ |\ w\in W, \lambda\in P\}
\qquad\hbox{and}\qquad
\{ X^\mu T_{v^{-1}}^{-1}\ |\ \mu\in P, v\in W\}
\end{equation}
are bases of $\tilde H$.
If $p$ is an alcove walk then the \emph{weight} of $p$ and the
\emph{final direction} of $p$ are 
\begin{equation}
\hbox{$\wt(p)\in P$ and
$\varphi(p)\in W$}
\qquad\hbox{such that}\qquad
\hbox{$p$ ends in the alcove $\wt(p)+\varphi(p) A$.}
\end{equation}
Let
\begin{equation}
\begin{array}{rl}
f^-(p) &= \hbox{(number of negative folds of $p$)},
\\
f^+(p) &= \hbox{(number of positive folds of $p$)},
\qquad\hbox{and} \\
f(p) &= \hbox{(total number of folds of $p$)}.
\end{array}
\end{equation}
The following theorem provides a combinatorial formulation of the transition matrix
between the bases in (3.7).  It is a $q$-version of the main result of [LP] and an extension of
Corollary 6.1 of [Sc].

\begin{thm}   Use notations as in (3.4).  Let $\lambda\in P$ and $w\in W$.  Fix
a minimal length walk $p_w=c_{i_1}^-c_{i_2}^-\cdots c_{i_r}^-$ 
from $A$ to $wA$ and a minimal length walk 
$p_\lambda=c_{j_1}^{\epsilon_1}\cdots c_{j_s}^{\epsilon_s}$ 
from $A$ to $t_\lambda A$.
Then, with notations as in (3.8) and (3.9),
$$T_{w^{-1}}^{-1}X^\lambda 
= \sum_p (-1)^{f^-(p)}(q-q^{-1})^{f(p)} X^{\wt(p)}T_{\varphi(p)^{-1}}^{-1},$$
where the sum is over all alcove walks 
$p=c_{i_1}^-\cdots c_{i_r}^- p_{j_1}\cdots p_{j_s}$ such that 
$p_{j_k}$ is either $c_{j_k}^{\epsilon_k}$, $c_{j_k}^{-\epsilon_k}$ or $f_{j_k}^{\epsilon_k}$.
\end{thm}
\begin{proof}
The product
$p_wp_\lambda = c_{i_1}^-c_{i_2}^-\cdots c_{i_r}^-c_{j_1}^{\epsilon_1}\cdots c_{j_s}^{\epsilon_s}$ 
may not necessarily be walk, but its straightening produces a sum of walks, and
this decomposition gives the formula in the statement.
\end{proof}

\begin{remark}  The \emph{initial direction} $\iota(p)$ and the
\emph{final direction} $\varphi(p)$ of an alcove walk $p$ appear naturally in 
Theorem 3.3.  These statistics also appear
in the Pieri-Chevalley formula in the $K$-theory of the flag variety (see [PR], [GR], [Br] and [LP]).
\end{remark}

\begin{remark}  In Theorem 3.3, for certain $\lambda$ the walk $p_\lambda$ may be
chosen so that all the terms in the expansion of $T_{w^{-1}}^{-1}X^\lambda$ have the same
sign.  For example, if $\lambda$ is dominant, then $p_\lambda$ can be taken with
all $\epsilon_k=+$, in which case all folds which appear in the straightening of
$p_wp_\lambda$ will be positive folds and so all terms in the expansion will be positive.
If $\lambda$ is antidominant then $p_\lambda$ can be taken with all $\epsilon_k=-$
and all terms in the expansion will be negative.  This fact gives positivity results
for products in the cohomology and the $K$-theory of the flag variety (see [PR], [Br]).
\end{remark}

\begin{remark}
The affine Hecke algebra $\tilde H$ has basis 
$\{ X^\lambda T_{w^{-1}}^{-1}\ |\ \lambda\in P, w\in W\}$ in bijection with
the alcoves in $\Omega\times \fh_\RR^*$,
where $X^\lambda T_{w^{-1}}^{-1}$ is the image in $\tilde H$ of a minimal length alcove walk
from $A$ to the alcove $\lambda+wA$.   Changing the orientation of the walls of the
alcoves changes the resulting basis in the affine Hecke algebra $\tilde H$.  The orientation
in (3.1) is the one such that 
\begin{equation}
\hbox{the most negative point
is $-\infty\rho$, deep in the chamber $w_0C$}.
\end{equation}
Another standard orientation is where
\begin{equation}
\hbox{the most negative point is
the center of the fundamental alcove $A$.}
\end{equation}
Using the orientation of the walls given by (3.11)
produces the basis commonly denoted $\{ T_w\ |\ w\in \widetilde{W}\}$ by taking
$T_w$ to be the image in $\tilde H$ of a minimal length alcove walk from $A$ to $w^{-1}A$.
Since $T_i^{-1} = T_i - (q-q^{-1})$ the transition matrix between the basis
$\{ X^\lambda T_{w^{-1}}^{-1}\ |\ \lambda\in P, w\in W\}$
and the basis $\{ T_w\ |\ w\in \widetilde{W}\}$ is triangular.  
\end{remark}

\end{subsection}

\end{section}

\begin{section}{Polynomials and symmetric functions}

Recall, from (3.6), that the \emph{finite Hecke algebra} is the subalgebra of the 
affine Hecke algebra $\tilde H$ given by 
$$H = \hbox{span}\{ T_{w^{-1}}^{-1}\ |\ w\in W\}.$$
Let $\mathbf{1}_0$ be the element of $H$ given by
\begin{equation}
\mathbf{1}_0^2 = \mathbf{1}_0
\qquad\hbox{and}\qquad
T_{w^{-1}}^{-1}\mathbf{1}_0 = q^{-\ell(w)}\mathbf{1}_0,
\end{equation}
for $w\in W$.  Two explicit formulas for $\mathbf{1}_0$ are
$$\mathbf{1}_0 = \frac{1}{W_0(q^{-2})}\sum_{w\in W} q^{-\ell(w)}T_{w^{-1}}^{-1}
=\frac{1}{W_0(q^2)}\sum_{w\in W} q^{\ell(w)}T_w,
\qquad\hbox{where}\qquad
W_0(t) = \sum_{w\in W} t^{\ell(w)}$$
is the \emph{Poincar\'e polynomial of $W$}.
 
As observed in (3.6), $\KK[P] = \hbox{span}\{ X^\lambda\ |\ \lambda\in P\}$ is a subalgebra
of $\tilde H$.  The vector space $\KK[P]$ also sits inside $\tilde H$ in a
\emph{different} way.
Since $\{ X^\lambda\mathbf{1}_0\ |\ \lambda\in P\}$ is a basis of $\tilde H\mathbf{1}_0$
there is
\begin{equation}
\hbox{a vector space isomorphism}\qquad
\begin{matrix}
\KK[P] &\longrightarrow &\tilde H \mathbf{1}_0 \\
f &\longmapsto &f\mathbf{1}_0.
\end{matrix}
\end{equation}

The ring of \emph{symmetric functions} is 
\begin{equation}
\KK[P]^W = \{ f\in \KK[P]\ |\ \hbox{$wf = f$ for all $w\in W$}\}.
\end{equation}
By a theorem of Bernstein (see [NR, Theorem 1.4]) this subalgebra of $\tilde H$ is the center of $\tilde H$,
\begin{equation}
\KK[P]^W = Z(\tilde H).
\end{equation}
The \emph{spherical Hecke algebra} is the ring $\mathbf{1}_0X^\lambda\mathbf{1}_0$
and the restriction of the map (4.2) to $Z(\tilde H)$ is the \emph{Satake isomorphism}
of the following theorem.

\begin{thm}   Let $\mathbf{1}_0$ and $\KK[P]^W$ be as in (4.1) and (4.3), respectively.  Then
$$\begin{array}{lcc}
\KK[P]^W = Z(\tilde H) &\longrightarrow &\mathbf{1}_0 \tilde H \mathbf{1}_0 \\
\ \ f &\longmapsto &f \mathbf{1}_0
\end{array}
\qquad\hbox{is a $\KK$-algebra isomorphism.}
$$
\end{thm}
\begin{proof}
The map is a well defined homomorphism since, if $f, f_1, f_2\in Z(\tilde H)$, then
$f\mathbf{1}_0 = f\mathbf{1}_0^2 = \mathbf{1}_0f\mathbf{1}_0$ and
$f_1f_2\mathbf{1}_0 = f_1f_2\mathbf{1}_0^2 = f_1\mathbf{1}_0f_2\mathbf{1}_0$.

Suppose that $\langle \lambda,\alpha_i^\vee\rangle >0$ so that
$\lambda$ is on the positive side of $H_{\alpha_i}$.  Then, by Proposition 3.2e,
\begin{align*}
\mathbf{1}_0X^{s_i\lambda}\mathbf{1}_0
&= q^{-1}\mathbf{1}_0T_{s_i} X^{\lambda}\mathbf{1}_0 \\
&= q^{-1}\mathbf{1}_0X^{\lambda}T_{s_i}\mathbf{1}_0 
-q^{-1}(q-q^{-1})\mathbf{1}_0(X^{s_i\lambda+\alpha_i}+\cdots+X^{\lambda-\alpha_i}
+X^\lambda)\mathbf{1}_0 \\
&= \mathbf{1}_0X^{\lambda}\mathbf{1}_0 
-(1-q^{-2})\mathbf{1}_0(X^{s_i\lambda+\alpha_i}+\cdots+X^{\lambda-\alpha_i}
+X^\lambda)\mathbf{1}_0.
\end{align*}
so that
\begin{equation}
\mathbf{1}_0X^{s_i\lambda}\mathbf{1}_0
= q^{-2}\mathbf{1}_0X^{\lambda}\mathbf{1}_0 
-(1-q^{-2})\mathbf{1}_0(X^{s_i\lambda+\alpha_i}+\cdots
+X^{\lambda-\alpha_i})\mathbf{1}_0
\end{equation}
or, equivalently,
\begin{equation}
\mathbf{1}_0(X^{s_i\lambda}+\cdots+X^{\lambda-\alpha_i})\mathbf{1}_0
=q^{-2}\mathbf{1}_0(X^{s_i\lambda+\alpha_i}+\cdots+X^{\lambda})\mathbf{1}_0.
\end{equation}
From the relation in (4.5),
\begin{align*}
\mathbf{1}_0X^{s_i\lambda}\mathbf{1}_0
-\mathbf{1}_0X^{s_i\lambda+\alpha_i}\mathbf{1}_0
&= q^{-2}\mathbf{1}_0X^{\lambda}\mathbf{1}_0 
-(1-q^{-2})\mathbf{1}_0(X^{s_i\lambda+\alpha_i}+\cdots
+X^{\lambda-\alpha_i})\mathbf{1}_0 \\
&\qquad\quad
- q^{-2}\mathbf{1}_0X^{\lambda-\alpha_i}\mathbf{1}_0 
+(1-q^{-2})\mathbf{1}_0(X^{s_i\lambda+2\alpha_i}+\cdots
+X^{\lambda-2\alpha_i})\mathbf{1}_0,
\end{align*}
so that
\begin{equation}
\mathbf{1}_0X^{s_i\lambda}\mathbf{1}_0
=q^{-2}\mathbf{1}_0X^{\lambda}\mathbf{1}_0
+q^{-2}\mathbf{1}_0X^{s_i\lambda+\alpha_i}\mathbf{1}_0
-\mathbf{1}_0X^{\lambda-\alpha_i}\mathbf{1}_0.
\end{equation}
It follows from these relations that any element of $\mathbf{1}_0 \tilde H \mathbf{1}_0$
can, inductively, be written as a linear combination of the elements $\mathbf{1}_0X^\lambda\mathbf{1}_0$,
$\lambda\in {P^+}$.   Using Theorem 3.3 to expand $\mathbf{1}_0 X^\lambda$
in terms of the basis $\{ X^\mu T_{v^{-1}}^{-1}\ |\ \mu\in P, v\in W\}$ produces 
$$\mathbf{1}_0X^\lambda 
= X^{w_0\lambda} T_{w_0}^{-1} 
+ \sum_{\mu>w_0\lambda} d_{\mu,v} X^\mu T_{v^{-1}}^{-1},$$
and, since these leading terms are all different (as $\lambda$ runs over $P^+$), it follows that  
\begin{equation}
\mathbf{1}_0\tilde H\mathbf{1}_0
\quad\hbox{has basis}\quad
\{ \mathbf{1}_0 X^\lambda \mathbf{1}_0\ |\ \lambda\in P^+\}.
\end{equation}

The \emph{orbit sums}
\begin{equation}
m_\lambda = \sum_{\gamma\in W\lambda} X^\gamma,
\qquad\lambda\in P^+,
\end{equation}
form a basis of $\KK[P]^W$.
The relation in Proposition 3.2e implies that, if $f\in \KK[P]^W$ then $T_w f = f T_w$ for all $w\in W$,
and so the $m_\lambda \mathbf{1}_0 = \mathbf{1}_0 m_\lambda \mathbf{1}_0$
are in $\mathbf{1}_0 \tilde H \mathbf{1}_0$.  Viewing these in terms of the basis
$\{ X^\mu T_{v^{-1}}^{-1}\ |\ \mu\in P, v\in W\}$ of $\tilde H$ one sees that the 
$m_\lambda \mathbf{1}_0$, $\lambda\in P^+$,  are linearly independent  and so
\begin{equation}
\mathbf{1}_0\tilde H\mathbf{1}_0
\quad\hbox{has basis}\quad
\{ m_\lambda \mathbf{1}_0\ |\ \lambda\in P^+\}.
\end{equation}
The point is that the transition matrix between the basis in (4.8) and the basis in (4.10) is
triangular.
\end{proof}

\begin{subsection}{Hall-Littlewood polynomials}

For $\mu\in P$ let $W_\mu = \mathrm{Stab}(\mu)$ be the stabilizer of $\mu$.
The \emph{Poincar\'e polynomial} of $W_\mu$ is 
\begin{equation}
W_\mu(t) = \sum_{w\in W_\mu} t^{\ell(w)}.
\end{equation}
For $\mu\in P$, the \emph{Hall-Littlewood polynomial}  or \emph{Macdonald spherical
function} $P_\mu(X;t)$ is the element of
$\KK[P]^W$ defined by
\begin{equation}
P_\mu(X;q^{-2}) \mathbf{1}_0
=\Big(\sum_{w\in W^\mu} q^{-\ell(w)} T_{w^{-1}}^{-1}\Big)
X^\mu \mathbf{1}_0,
\end{equation}
where $W^\mu$ is a set of minimal length coset representatives for the
cosets in $W/W_\mu$.  
Since every element $w\in W$ has a unique
expression $w=uv$ with $u\in W^\mu$ and $v\in W_\mu$,
\begin{align*}
\sum_{w\in W^\mu} &q^{-\ell(w)} T_{w^{-1}}^{-1} X^\mu \mathbf{1}_0 
=\frac{1}{W_\mu(q^{-2})} \left(\sum_{u\in W^\mu} q^{-\ell(u)}T_{u^{-1}}^{-1}\right)
X^\mu
\left(\sum_{v\in W_\mu} q^{-\ell(v)}T_{v^{-1}}^{-1} \right)\mathbf{1}_0 \\
&=\frac{1}{W_\mu(q^{-2})} \left(\sum_{u\in W^\mu} q^{-\ell(u)}T_{u^{-1}}^{-1}\right)
\left(\sum_{v\in W_\mu} q^{-\ell(v)}T_{v^{-1}}^{-1} \right)
X^\mu
\mathbf{1}_0 
=\frac{W_0(q^{-2})}{W_\mu(q^{-2})}\mathbf{1}_0 X^\mu \mathbf{1}_0,
\end{align*}
and hence 
$$P_\mu(X, q^{-2})\mathbf{1}_0
\quad\hbox{is exactly}\quad
\mathbf{1}_0X^\mu\mathbf{1}_0
\quad\hbox{except normalized}
$$
so that the coefficient of $X^\mu \mathbf{1}_0$ is $1$.

Macdonald's formula for the spherical functions on a $p$-adic group [Mac1, Theorem 4.1.2]
is
\begin{equation}
P_\mu(X;q^{-2}) = \frac{1}{W_\mu(q^{-2})}
\sum_{w\in W} w\left(X^\mu\prod_{\alpha\in R^+}
\frac{1-q^{-2}X^{-\alpha}}{1-X^{-\alpha}}\right).
\end{equation}
See [NR, Theorem 2.9a] for a proof in this context.  

The following theorem gives additional formulas for $P_\lambda(X;q^{-2})$.
\begin{equation}
\hbox{A \emph{positively folded} alcove walk is an alcove walk with no negative folds.}
\end{equation}
In the following theorem we shall consider alcove walks which do not necessarily begin
at $A$.  This is the natural way to account for the sum over $W^\lambda$ which appears in 
the definition of $P_\lambda$ in (4.12).  The \emph{type} of $p$ is the sequence of 
labels of the folds and the wall crossings of $p$.

\begin{thm} 
For $\lambda\in P^+$ let $t_\lambda\in \widetilde{W}$ be the translation in $\lambda$
and let $n_\lambda$ be the maximal length element in the double coset $Wt_\lambda W$.
\item[(a)] {\rm [Sc, Theorem 1.1]}
Let $\lambda\in P^+$ and fix a minimal length walk
$p_\lambda=c_{i_1}^+\cdots c_{i_\ell}^+$ from $A$ to $\lambda+A$.  Let 
$$B_q(p_\lambda) 
= \{ \hbox{positively folded alcove walks of type $(i_1,\ldots, i_\ell)$ which begin at $wA$}
\ |\ w\in W^\lambda\}
$$
Then 
$$P_\lambda(X;q^{-2}) = 
\sum_{p\in B_q(p_\lambda)} q^{-(\ell(\iota(p))+\ell(\varphi(p))-f(p))}
(1-q^{-2})^{f(p)}X^{\wt(p)},$$
where $\iota(p)$ is the alcove where $p$ begins,
$\wt(p)+\varphi(p)A$ is the alcove where $p$ ends, and $f(p)$ is the number of folds in $p$.
\item[(b)]  Let $\lambda\in P^+$.
Then
$$q^{\ell(w_0)}W_0(q^{-2})P_\lambda\mathbf{1}_0
= \sum_{x\in Wt_\lambda W} q^{\ell(x)-\ell(n_\lambda)}T_x.$$
\end{thm}
\begin{proof}
(a)  The proof is accomplished by using Theorem 3.3 to expand
the sum in (4.12).  Since all crossings in the walk $p_\lambda$ are positive crossings
Theorem 3.3 gives
$$\left(\sum_{w\in W^\lambda} q^{-\ell(w)} T_{w^{-1}}^{-1}\right) X^\lambda
= \sum_{w\in W^\lambda} q^{-\ell(w)} 
\sum_{p\in B_q(p_\lambda)\atop \iota(p)=w} (q-q^{-1})^{f(p)}X^{\wt(p)}T_{\varphi(p)}^{-1}.$$
Hence
$$P_\lambda \mathbf{1}_0
= \sum_{p\in B_q(\lambda)} q^{-\ell(\iota(p))}(q-q^{-1})^{f(p)}X^{\wt(p)}q^{-\ell(\varphi(p))}\mathbf{1}_0.
$$
(b) Let $\lambda\in P^+$.  
Let $W_\lambda={\rm Stab}(\lambda)$ and let
$w_0$ and $w_\lambda$ be the maximal length elements in $W$ and 
$W_\lambda$, respectively.  
Let $m_\lambda$ and $n_\lambda$ be the minimal and maximal 
length elements respectively in the double coset $Wt_\lambda W$.
If $\lambda = 2\omega_2$ in type $C_2$, then $W_{\lambda} = \{1, s_1\}$,
$w_\lambda = s_1$, $w_0=s_1s_2s_1s_2$, $\ell(t_\lambda) = 6$,
$\ell(m_\lambda) = 3$, and $\ell(n_\lambda)=10$.  Labeling the alcove
$wA$ by the element $w$, the 32 alcoves $wA$ with $w\in Wt_\lambda W$ 
make up the four shaded diamonds.
$$
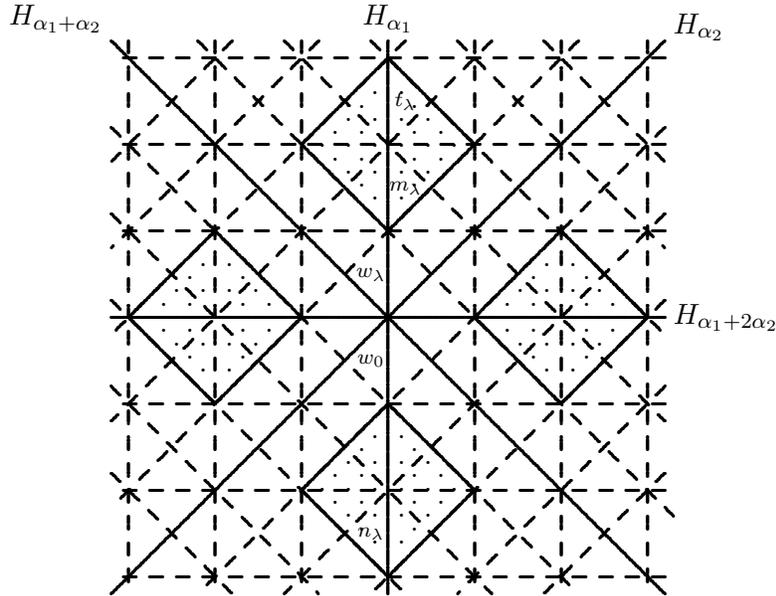
\beginpicture
\setcoordinatesystem units <1.15cm,1.15cm>         
\setplotarea x from -4.5 to 4.5, y from -3.5 to 4    
\put{$H_{\alpha_1}$}[b] at 0 3.3
\put{$H_{\alpha_2}$}[bl] at 3.3 3.2
\put{$H_{\alpha_1+\alpha_2}$}[br] at -3.3 3.3
\put{$H_{\alpha_1+2\alpha_2}$}[l] at 3.3 0
\put{$\scriptstyle{t_\lambda}$} at 0.2 2.5
\put{$\scriptstyle{m_\lambda}$} at 0.2 1.5
\put{$\scriptstyle{n_\lambda}$} at -0.2 -2.5
\put{$\scriptstyle{w_0}$} at -0.2 -0.5
\put{$\scriptstyle{w_\lambda}$} at -0.2 0.5
\plot -3.2 -3.2   3.2 3.2 /
\plot  3.2 -3.2  -3.2 3.2 /
\plot  0  3.2   0 -3.2 /
\plot  3.2  0  -3.2  0 /
\plot -1 2  0 3 /
\plot -1 2  0 1 /
\plot  0 1  1 2 /
\plot  0 3  1 2 /
\plot -1 -2  0 -3 /
\plot -1 -2  0 -1 /
\plot  0 -1  1 -2 /
\plot  0 -3  1 -2 /
\plot  1 0  2 1 /
\plot  2 1  3 0 /
\plot  1 0  2 -1 /
\plot  2 -1  3 0 /
\plot  -1 0  -2 1 /
\plot  -2 1  -3 0 /
\plot  -1 0  -2 -1 /
\plot  -2 -1  -3 0 /
\setdashes
\plot  -3 3.2   -3 -3.2 /
\plot  -2 3.2   -2 -3.2 /
\plot  -1 3.2   -1 -3.2 /
\plot   1 3.2    1 -3.2 /
\plot   2 3.2    2 -3.2 /
\plot   3 3.2    3 -3.2 /
\plot  3.2 -3   -3.2 -3 /
\plot  3.2 -2   -3.2 -2 /
\plot  3.2 -1   -3.2 -1 /
\plot  3.2  1   -3.2  1 /
\plot  3.2  2   -3.2  2 /
\plot  3.2  3   -3.2  3 /
\plot  3.2 -1.8   1.8 -3.2 /
\plot  3.2 -0.8   0.8 -3.2 /
\plot  3.2 0.2    -0.2 -3.2 /
\plot  3.2 1.2   -1.2 -3.2 /
\plot  3.2 2.2   -2.2 -3.2 /
\plot  2.2 3.2   -3.2 -2.2 /
\plot  1.2 3.2   -3.2 -1.2 /
\plot  0.2 3.2   -3.2 -0.2 /
\plot  -0.8 3.2  -3.2 0.8 /
\plot  -1.8 3.2  -3.2 1.8 /
\plot  -2.2 3.2   3.3 -2.3 /
\plot  -1.2 3.2   3.2 -1.2 /
\plot  -0.2 3.2   3.2 -0.2 /
\plot  0.8 3.2  3.2 0.8 /
\plot  1.8 3.2  3.2 1.8 /
\plot  -3.2 -1.8   -1.8 -3.2 /
\plot  -3.2 -0.8   -0.8 -3.2 /
\plot  -3.2 0.2    0.2 -3.2 /
\plot  -3.2 1.2   1.2 -3.2 /
\plot  -3.2 2.2   2.2 -3.2 /
\vshade 1 0 0   2 -1 1 /
\vshade 2 -1 1  3 0 0 /
\vshade -2 -1 1   -1 0 0 /
\vshade -3 0 0  -2 -1 1 /
\vshade -1 2 2   0 1 3 /
\vshade  0 1 3   1 2 2 /
\vshade -1 -2 -2   0 -3 -1 /
\vshade  0 -3 -1   1 -2 -2 /
\endpicture
$$
Then
\begin{align*}
q^{\ell(w_0)}W_0(q^{-2})P_\lambda\mathbf{1}_0 
&=
q^{\ell(w_0)}W_0(q^{-2})q^{-2\ell(w_0)+2\ell(w_\lambda)}
\left(\sum_{u\in W^\lambda} q^{\ell(u)}T_u\right)
X^\lambda \mathbf{1}_0\\
&= 
q^{\ell(w_0)}q^{-2\ell(w_0)}W_0(q^{2})q^{-2\ell(w_0)+2\ell(w_\lambda)}
\left(\sum_{u\in W^\lambda} q^{\ell(u)}T_u\right)
T_{m_\lambda}T_{w_0w_\lambda}
\mathbf{1}_0\\
&= 
q^{-3\ell(w_0+2\ell(w_\lambda)}
\left(\sum_{u\in W^\lambda} q^{\ell(u)}T_u\right)
T_{m_\lambda} q^{\ell(w_0)-\ell(w_\lambda)}
\left(\sum_{w\in W} q^{\ell(w)}T_w\right) \\
&= 
q^{-2\ell(w_0)+\ell(w_\lambda)}
\left(\sum_{u\in W^\lambda\atop w\in W} q^{\ell(u)+\ell(w)}T_{um_\lambda w}\right) \\
&= 
q^{-2\ell(w_0)+\ell(w_\lambda)}
\left(\sum_{x\in Wt_\lambda W} q^{\ell(x)-\ell(m_\lambda)}T_x\right), 
\end{align*}
and the result follows from the identity
$\ell(n_\lambda) = \ell(w_0)+\ell(t_\lambda)
=\ell(w_0)+\ell(w_0)-\ell(w_\lambda)+\ell(m_\lambda)$.
\end{proof}

\begin{remark} The set $B_q(p_\lambda)$ appearing in Theorem 4.2a
is heavily dependent on the choice of $p_\lambda$.  One way to make this 
seem less dependent on this choice is to consider the region 
$$[\lambda] = \{ x = \sum_i x_i\omega_i\ |\ 0\le x_i\le \lambda_i\},
\qquad\hbox{where}\qquad
\lambda = \sum_i \lambda_i\omega_i.$$
This region is a union of alcoves and any minimal length walk $p_\lambda$ from
$A$ to $\lambda+w_0A$ lies in $[\lambda]$.  Foldings of the walk $p_\lambda$
are then produced by folding the region $[\lambda]$ along the ``creases'' formed by the
hyperplanes.  This process produces a bijection between the paths in $B_q(p_\lambda)$ and
the set $B_q([\lambda])$ of ``positively folded foldings'' of the region $[\lambda]$, and the set
$B_q([\lambda])$ does not depend the choice of an initial path.
The moral is that the best way to forget the choice of the initial path $p_\lambda$ is to 
remember all the possible initial paths all at once.
This translation of foldings was explained to me by J.\ Ramagge in Fall 2000.
\end{remark}

\begin{remark} Let $B_q(p_\lambda)$ be as in Theorem 4.2a and let
$p\in B_q(p_\lambda)$.  Suppose that $p$ has $f$ folds.
For $0\le i\le f$, let $p^{(i)}$ be the positively folded walk in $B_q(p_\lambda)$
which coincides with $p$ up
to the $i$th fold and is nonfolded thereafter.  Then
$p^{(0)}, \ldots, p^{(f)}$ is a sequence of positively folded walks such that 
$$p^{(f)}=p,
\quad \iota(p^{(i)}) = \iota(p),
\quad \varphi(p^{(0)}) = \iota(p),
\quad\hbox{and}\quad \varphi(p^{(i)})=s_\alpha\varphi(p^{(i-1)})
$$
if the $i$th fold is on the hyperplane $H_{\alpha,j}$.  Since
$$\varphi(p^{(i-1)})> \varphi(p^{(i)})
\quad\hbox{and}\quad
(-1)^{\ell(\varphi(p^{(i)}))}
=(-1)^{\ell(s_\alpha)}(-1)^{\varphi(p^{(i-1)})}
=(-1)(-1)^{\varphi(p^{(i-1)})},$$
$\ell(\varphi(p^{(i-1)}))-\ell(\varphi(p^{(i)}))-1$ is an even
integer $\ge 0$.  It follows that 
\begin{align*}
\ell(\iota(p))+\ell(\varphi(p))-f(p)
&=\ell(\iota(p))-\ell(\varphi(p^{(1)}))-1 \\
&\qquad+\ell(\varphi(p^{(1)}))-\ell(\varphi(p^{(2)}))-1 \\
&\qquad\qquad+\ell(\varphi(p^{(2)}))-\ell(\varphi(p^{(3)}))-1 \\
&\qquad\qquad\qquad+\cdots
+\ell(\varphi(p^{(f-1)}))-\ell(\varphi(p^{(f)}))-1 \quad+2\ell(\varphi(p))
\end{align*}
is an even integer $\ge 0$.  This proves that $f(p)\le \ell(\iota(p))-\ell(\varphi(p))$ and that
$P_\lambda(X;q^{-2})$ really is a polynomial in the variable $q^{-2}$.
\end{remark}

\end{subsection}

\begin{subsection}{Demazure operators}

The group $W$ acts on $\KK[P]=\hbox{span}\{X^\lambda\ |\ \lambda\in P\}$ by 
\begin{equation}
wX^\lambda = X^{w\lambda},
\qquad\hbox{for $w\in W$, $\lambda\in P$.}
\end{equation}
For each $1\le i\le n$, define \emph{Demazure operators}
$$\Delta_i\colon \KK[P]\longrightarrow \KK[P]
\qquad\hbox{and}\qquad
\tilde \Delta_i\colon \KK[P]\longrightarrow \KK[P]
$$
by
\begin{equation}
\Delta_i f = \frac{1}{1-X^{-\alpha_i}}(1-s_i)f
\qquad\hbox{and}\qquad
\tilde \Delta_i f = \frac{1}{X^{\alpha_i}-1}(X^{\alpha_i}-s_i)f,
\end{equation}
respectively.  

Via the isomorphism in (4.2), the vector space $\KK[P]$ is an $\tilde H$-module.
Let 
\begin{equation}
C_i=q^{-2}+q^{-1}T_{s_i} = 1+q^{-1}T_{s_i}^{-1}=(1+q^{-2})\mathbf{1}_i,
\end{equation}
where
$\mathbf{1}_i$ is the element of
$H$ such that $\mathbf{1}_i^2 = \mathbf{1}_i$ and $T_{s_i}^{-1}\mathbf{1}_i = q^{-1}\mathbf{1}_i$.
The element $\mathbf{1}_i$ is the rank 1 version of the element $\mathbf{1}_0$ in (4.1).

The following proposition shows that, at $q^{-2}=0$, the action of $C_i$ on $\KK[P]$ 
is the Demazure operator $\tilde \Delta_i$.
In geometry, the Demazure operators arise naturally as push-pull operators on the
K-theory of the flag variety (see [PR, Proposition]).  

\begin{prop}  Let $\rho = \omega_1+\cdots +\omega_n$ as in (2.16).
As operators on $\KK[P]$,
\item[(a)] $\tilde\Delta_i = X^{-\rho}\Delta_i X^\rho = \Delta_i+s_i$,
\item[(b)] $\displaystyle{
C_i=(1-q^{-2})\Delta_i + (s_i+q^{-2})
= (1+s_i)\left(\frac{1-q^{-2}X^{-\alpha_i}}{1-X^{-\alpha_i}}\right),}
$
\end{prop}
\begin{proof}  (a)
Let $\lambda\in P$.  Since $s_i\rho = \rho-\langle\rho,\alpha_i^\vee\rangle\alpha_i
= \rho-\alpha_i$,
$$
(X^{-\rho}\Delta_i X^\rho)(X^\lambda)
= X^{-\rho}\frac{ X^{\lambda+\rho} - X^{s_i\lambda+\rho-\alpha_i} }{1-X^{-\alpha_i} }
=\frac{X^\lambda - X^{s_i\lambda-\alpha_i} }{ 1-X^{-\alpha_i}}
=\tilde\Delta_i(X^\lambda),
$$
and, as operators,
$$\Delta_i+s_i = \frac{1}{1-X^{-\alpha_i}}(1-s_i)+s_i
=\frac{1}{1-X^{-\alpha_i}}(1-s_i+s_i-X^{-\alpha_i}s_i)
=\tilde \Delta_i.$$
(b)  Using Proposition 3.2e,
\begin{align*}
q^{-1}T_{s_i} X^\lambda \mathbf{1}_0
&= \left(
q^{-1} X^{s_i\lambda}T_{s_i} 
+ (1-q^{-2})\frac{ X^\lambda-X^{s_i\lambda}}{1-X^{-\alpha_i}}
\right)\mathbf{1}_0 
= \left(
X^{s_i\lambda} 
+ (1-q^{-2})\frac{ X^\lambda-X^{s_i\lambda}}{1-X^{-\alpha_i}}
\right)\mathbf{1}_0 \\
&= \left(
\frac{X^{s_i\lambda}-X^{s_i\lambda-\alpha_i}+X^\lambda-X^{s_i\lambda}
-q^{-2}(X^{\lambda}-X^{s_i\lambda}) } 
{1-X^{-\alpha_i}}
\right)\mathbf{1}_0 \\
&= \big(\tilde \Delta_i - q^{-2}\Delta_i\big)(X^\lambda)\mathbf{1}_0,
\end{align*}
and the first formula in (b) now follows from the second formula in (a).
Then
\begin{eqnarray*}
C_i
&=&(q^{-2}+\tilde\Delta_i-q^{-2}\Delta_i) 
=\frac{1}{1-X^{-\alpha_i}}(q^{-2}-q^{-2}X^{-\alpha_i}+1-X^{-\alpha_i}s_i
-q^{-2}+q^{-2}s_i) \\
&=& \frac{1-q^{-2}X^{-\alpha_i}}{1-X^{-\alpha_i}}
+s_i\frac{q^{-2}-X^{\alpha_i}}{1-X^{\alpha_i}}
=(1+s_i)\left(\frac{1-q^{-2}X^{-\alpha_i}}{1-X^{-\alpha_i}}\right).
\end{eqnarray*}
\end{proof}

\begin{remark}  Slightly renormalizing the generators of the affine Hecke algebra by setting 
$\tilde T_i = q^{-1}T_{s_i}$ allows one to let
$q^{-1}=0$  so that $\tilde T_i$ acts on $\KK[P]$ by $\tilde \Delta_i$.
This is the action of the \emph{nil affine Hecke algebra} on $\KK[P]$.
Since the $\tilde T_i$ satisfy the braid relations so do the $\tilde \Delta_i$.  The first 
formula in Proposition 4.1 shows that $\Delta_i$ is a conjugate of $\tilde \Delta_i$ and
so the $\Delta_i$ also satisfy the braid relations.  Although $C_i$ equals $\tilde \Delta_i$
at $q^{-2}=0$, the operators $C_i$ do \emph{not} satisfy
the braid relations.   Furthermore, if
$w_0=s_{i_1}\cdots s_{i_\ell}$ is a reduced word for the longest
element then 
$$C_{i_1}\cdots C_{i_\ell} 
= W_0(q^{-2})\mathbf{1}_0 + q^{-2}(\hbox{extra terms}).$$
In contrast to the case for Weyl characters (when $q^{-2}=0$), because
of the $q^{-2}(\hbox{extra terms})$ the Hall-Littlewood polynomial cannot be generated by
applying the product $C_{i_1}\cdots C_{i_\ell}$
unless one somehow knows how to throw away the extra terms.
\end{remark}

\begin{remark}
As operators on $\KK[P]$,
$$q-T_{s_i} = \left(\frac{q^{-1}-qX^{-\alpha_i} }{1-X^{-\alpha_i}}\right)
(1-s_i)
\quad\hbox{and}\quad
\mathbf{1}_0 = \frac{1}{W_0(q^{-2})}\sum_{w\in W} w \left( \prod_{\alpha\in R^+} 
\frac{1-q^{-2}X^{-\alpha}}{1-X^{-\alpha}}\right).
$$
The second formula is equivalent to 
Macdonald's spherical function formula (4.13).
%
%
%
%
\end{remark}

\end{subsection}

\begin{subsection}{Root operators}

The idea of root operators is to 
give an alcove walk interpretation of the action of the operator
$C_i$ on $\KK[P]$ by considering the projections of the alcove
walks onto the line perpendicular to $H_{\alpha_i}$.    The main point is 
the identity (4.21) which gives a combinatorial description
of the action of $C_i$ on $\KK[P]$.  The appropriate
combinatorics is more or less forced by the \emph{Leibnitz rule} or
\emph{tensor product rule} for the operator $C_i$ given in (4.27).

Let $p$ be a positively folded alcove walk.
$$
\beginpicture
\setcoordinatesystem units <1.15cm,1.15cm>         
\setplotarea x from -4.5 to 4.5, y from -3.5 to 4    
\put{$H_{\alpha_1}$}[b] at 0 3.3
\put{$H_{\alpha_2}$}[bl] at 3.3 3.2
\put{$H_{\alpha_1+\alpha_2}$}[br] at -3.3 3.3
\put{$H_{\alpha_1+2\alpha_2}$}[l] at 4.3 0
\plot -3.2 -3.2   3.2 3.2 /
\plot  3.2 -3.2  -3.2 3.2 /
\plot  0  3.2   0 -3.2 /
\plot  4.2  0  -3.2  0 /
\arrow <5pt> [.2,.67] from 0.175 0.5 to 0.5 0.175 %
\arrow <5pt> [.2,.67] from 0.5 0.175 to 0.5 -0.175 %
\arrow <5pt> [.2,.67] from 0.5 -0.175 to 0.175 -0.5 %
\plot 0.175 -0.5  0.05 -0.5 /   
\plot 0.05 -0.5  0.05 -0.55 /  %
\arrow <5pt> [.2,.67] from 0.05 -0.55 to 0.175 -0.55   
\arrow <5pt> [.2,.67] from 0.175 -0.55 to 0.5 -0.825 %
\arrow <5pt> [.2,.67] from 0.5 -0.825 to 0.825 -0.5 %
\arrow <5pt> [.2,.67] from 0.825 -0.5 to 1.175 -0.5 %
\arrow <5pt> [.2,.67] from 1.175 -0.5 to 1.5 -0.175 %
\arrow <5pt> [.2,.67] from 1.5 -0.175 to 1.5 0.175 %
\arrow <5pt> [.2,.67] from 1.5 0.175 to 1.825 0.5 %
\arrow <5pt> [.2,.67] from 1.825 0.5 to 2.175 0.5 %
\arrow <5pt> [.2,.67] from 2.175 0.5 to 2.5 0.175 %
\arrow <5pt> [.2,.67] from 2.5 0.175 to 2.5 -0.175 %
\arrow <5pt> [.2,.67] from 2.5 -0.175 to 2.175 -0.5 %
\arrow <5pt> [.2,.67] from 2.175 -0.5 to 2.5 -0.825 %
\arrow <5pt> [.2,.67] from 2.5 -0.825 to 2.5 -1.175 %
\arrow <5pt> [.2,.67] from 2.5 -1.175 to 2.825 -1.5 %
\arrow <5pt> [.2,.67] from 2.825 -1.5 to 2.5 -1.825 %
\arrow <5pt> [.2,.67] from 2.5 -1.825 to 2.175 -1.5 %
\arrow <5pt> [.2,.67] from 2.175 -1.5 to 1.825 -1.5 %
\arrow <5pt> [.2,.67] from 1.825 -1.5 to 1.5 -1.825 %
\arrow <5pt> [.2,.67] from 1.5 -1.825 to 1.5 -2.175 %
\arrow <5pt> [.2,.67] from 1.5 -2.175 to 1.175 -2.5 %
\arrow <5pt> [.2,.67] from 1.175 -2.5 to 0.825 -2.5 %
\arrow <5pt> [.2,.67] from 0.825 -2.5 to 0.5 -2.825 %
\arrow <5pt> [.2,.67] from 0.5 -2.825 to 0.175 -2.5 %
\arrow <5pt> [.2,.67] from 0.175 -2.5 to -0.175 -2.5 %
\arrow <5pt> [.2,.67] from -0.175 -2.5 to -0.5 -2.825 %
\arrow <5pt> [.2,.67] from -0.5 -2.825 to -0.825 -2.5 %
\arrow <5pt> [.2,.67] from -0.825 -2.5 to -1.175 -2.5 %
\arrow <5pt> [.2,.67] from -1.175 -2.5 to -1.5 -2.825 %
\arrow <5pt> [.2,.67] from -1.5 -2.825 to -1.825 -2.5 %
\plot -1.825 -2.5  -1.95 -2.5 /   
\plot -1.95 -2.5  -1.95 -2.45 /  %
\arrow <5pt> [.2,.67] from -1.95 -2.45 to -1.825 -2.45   
\arrow <5pt> [.2,.67] from -1.825 -2.45 to -1.5 -2.175 %
\arrow <5pt> [.2,.67] from -1.5 -2.175 to -1.175 -2.45 %
\arrow <5pt> [.2,.67] from -1.175 -2.45 to -0.825 -2.45 %
\arrow <5pt> [.2,.67] from -0.825 -2.45 to -0.5 -2.175 %
\arrow <5pt> [.2,.67] from -0.5 -2.175 to -0.175 -2.45 %
\arrow <5pt> [.2,.67] from -0.175 -2.45 to 0.175 -2.45 %
\arrow <5pt> [.2,.67] from 0.175 -2.45 to 0.5 -2.175 %
\arrow <5pt> [.2,.67] from 0.5 -2.175 to 0.5 -1.825 %
\arrow <5pt> [.2,.67] from 0.5 -1.825 to 0.175 -1.5 %
\arrow <5pt> [.2,.67] from 0.175 -1.5 to -0.175 -1.5 %
\arrow <5pt> [.2,.67] from -0.175 -1.5 to -0.5 -1.175 %
\arrow <5pt> [.2,.67] from -0.5 -1.175 to -0.825 -1.5 %
\arrow <5pt> [.2,.67] from -0.825 -1.5 to -1.175 -1.5 %
\arrow <5pt> [.2,.67] from -1.175 -1.5 to -1.5 -1.825 %
\arrow <5pt> [.2,.67] from -1.5 -1.825 to -1.825 -1.5 %
\arrow <5pt> [.2,.67] from -1.825 -1.5 to -2.175 -1.5 %
\arrow <5pt> [.2,.67] from -2.175 -1.5 to -2.45 -1.825 %
\arrow <5pt> [.2,.67] from -2.45 -1.825 to -2.45 -2.175 %
\arrow <5pt> [.2,.67] from -2.45 -2.175 to -2.175 -2.5 %
\arrow <5pt> [.2,.67] from -2.175 -2.5 to -2.5 -2.825 %
\arrow <5pt> [.2,.67] from -2.5 -2.825 to -2.825 -2.5 %
\plot -2.825 -2.5  -2.95 -2.5 /   
\plot -2.95 -2.5  -2.95 -2.45 /  %
\arrow <5pt> [.2,.67] from -2.95 -2.45 to -2.825 -2.45   
\arrow <5pt> [.2,.67] from -2.825 -2.45 to -2.55 -2.175 %
\arrow <5pt> [.2,.67] from -2.55 -2.175 to -2.55 -1.825 %
\arrow <5pt> [.2,.67] from -2.55 -1.825 to -2.825 -1.5 %
\arrow <5pt> [.2,.67] from -2.825 -1.5 to -2.5 -1.175 %
\arrow <5pt> [.2,.67] from -2.5 -1.175 to -2.5 -0.825 %
\arrow <5pt> [.2,.67] from -2.5 -0.825 to -2.175 -0.5 %
\arrow <5pt> [.2,.67] from -2.175 -0.5 to -1.825 -0.5 %
\plot -1.825 -0.5   -1.72 -0.67  /  %
\plot -1.72  -0.67  -1.684 -0.625  / %
\arrow <5pt> [.2,.67] from -1.684 -0.625 to -1.8 -0.45   %
\arrow <5pt> [.2,.67] from -1.8 -0.45 to -1.5 -0.175 %
\arrow <5pt> [.2,.67] from -1.5 -0.175 to -1.5 0.175 %
\arrow <5pt> [.2,.67] from -1.5 0.175 to -1.175 0.5 %
\arrow <5pt> [.2,.67] from -1.175 0.5 to -0.825 0.5 %
\arrow <5pt> [.2,.67] from -0.825 0.5 to -0.5 0.825 %
\arrow <5pt> [.2,.67] from -0.5 0.825 to -0.5 1.175 %
\arrow <5pt> [.2,.67] from -0.5 1.175 to -0.175 1.5 %
\arrow <5pt> [.2,.67] from -0.175 1.5 to 0.175 1.5 %
\arrow <5pt> [.2,.67] from 0.175 1.5 to 0.5 1.175 %
\arrow <5pt> [.2,.67] from 0.5 1.175 to 0.5 0.825 %
\arrow <5pt> [.2,.67] from 0.5 0.825 to 0.8 0.525 %
\plot  0.8 0.525   0.69 0.375  /  %
\plot 0.69  0.375  0.709 0.33  / %
\arrow <5pt> [.2,.67] from 0.709 0.33 to 0.85 0.5   %
%
%
%
\arrow <5pt> [.2,.67] from 0.85 0.5 to 1.175 0.5 %
\arrow <5pt> [.2,.67] from 1.175 0.5 to 1.5 0.825 %
\arrow <5pt> [.2,.67] from 1.5 0.825 to 1.5 1.175 %
\arrow <5pt> [.2,.67] from 1.5 1.175 to 1.825 1.5 %
\arrow <5pt> [.2,.67] from 1.825 1.5 to 2.175 1.5 %
\arrow <5pt> [.2,.67] from 2.175 1.5 to 2.475 1.275 %
\plot  2.475 1.275   2.475 1.05  /  %
\plot 2.475  1.05  2.525 1.05  / %
\arrow <5pt> [.2,.67] from 2.525 1.05 to 2.525 1.275   %
\arrow <5pt> [.2,.67] from 2.525 1.275 to 2.825 1.5 %
\arrow <5pt> [.2,.67] from 2.825 1.5 to 3.175 1.5 %
\arrow <5pt> [.2,.67] from 3.175 1.55 to 2.825 1.55 %
\arrow <5pt> [.2,.67] from 2.825 1.55 to 2.5 1.825 %
\arrow <5pt> [.2,.67] from 2.5 1.825 to 2.175 1.55 %
\arrow <5pt> [.2,.67] from 2.175 1.55 to 1.825 1.55 %
\arrow <5pt> [.2,.67] from 1.825 1.55 to 1.5 1.825 %
\arrow <5pt> [.2,.67] from 1.5 1.825 to 1.5 2.175 %
\arrow <5pt> [.2,.67] from 1.5 2.175 to 1.175 2.5 %
\plot 1.175 2.5  1.05 2.5 /   
\plot 1.05 2.5  1.05 2.55 /  %
\arrow <5pt> [.2,.67] from 1.05 2.55 to 1.175 2.55   
\arrow <5pt> [.2,.67] from 1.175 2.55 to 1.5 2.825 %
\arrow <5pt> [.2,.67] from 1.5 2.825 to 1.825 2.5 %
\arrow <5pt> [.2,.67] from 1.825 2.5 to 2.175 2.5 %
\arrow <5pt> [.2,.67] from 2.175 2.5 to 2.5 2.825 %
\arrow <5pt> [.2,.67] from 2.5 2.825 to 2.825 2.5 %
\arrow <5pt> [.2,.67] from 2.825 2.5 to 3.175 2.5 %
\arrow <5pt> [.2,.67] from 3.175 2.5 to 3.5 2.175 %
\arrow <5pt> [.2,.67] from 3.5 2.175 to 3.825 2.5 %
\arrow <5pt> [.2,.67] from 3.825 2.5 to 4.175 2.5 %
\setdashes
\plot  -3 3.2   -3 -3.2 /
\plot  -2 3.2   -2 -3.2 /
\plot  -1 3.2   -1 -3.2 /
\plot   1 3.2    1 -3.2 /
\plot   2 3.2    2 -3.2 /
\plot   3 3.2    3 -3.2 /
\plot   4 3.2    4 -3.2 /
\plot  4.2 -3   -3.2 -3 /
\plot  4.2 -2   -3.2 -2 /
\plot  4.2 -1   -3.2 -1 /
\plot  4.2  1   -3.2  1 /
\plot  4.2  2   -3.2  2 /
\plot  4.2  3   -3.2  3 /
\plot  4.2 -0.8   1.8 -3.2 /
\plot  4.2 0.2   0.8 -3.2 /
\plot  4.2 1.2    -0.2 -3.2 /
\plot  4.2 2.2   -1.2 -3.2 /
\plot  4.2 3.2   -2.2 -3.2 /
\plot  2.2 3.2   -3.2 -2.2 /
\plot  1.2 3.2   -3.2 -1.2 /
\plot  0.2 3.2   -3.2 -0.2 /
\plot  -0.8 3.2  -3.2 0.8 /
\plot  -1.8 3.2  -3.2 1.8 /
\plot  -2.2 3.2   4.2 -3.2 /
\plot  -1.2 3.2   4.2 -2.2 /
\plot  -0.2 3.2   4.2 -1.2 /
\plot  0.8 3.2  4.2 -0.2 /
\plot  1.8 3.2  4.2 0.8 /
\plot  -3.2 -1.8   -1.8 -3.2 /
\plot  -3.2 -0.8   -0.8 -3.2 /
\plot  -3.2 0.2    0.2 -3.2 /
\plot  -3.2 1.2   1.2 -3.2 /
\plot  -3.2 2.2   2.2 -3.2 /
\endpicture
$$
Let $1\le i\le n$.
The projection of $p$ onto a line perpendicular to $H_{\alpha_i}$ is a
positively folded alcove walk $\bar p$ ``with respect to $\alpha_i$'' (the only important
information in the projection is the relative position of the walk to each of the hyperplanes
parallel to $H_{\alpha_i}$).
\begin{equation}
\beginpicture
\setcoordinatesystem units <1.15cm,0.75cm>         
\setplotarea x from -3.5 to 3.5, y from -1 to 1.1  
\put{$\scriptstyle{\vert}\longleftarrow$}[l] at -3 1.1
\put{$\longrightarrow\scriptstyle{\vert}$}[r] at 4.5 1.1
\put{$\scriptstyle{d_i^+(p)=d^+(\bar p)}$} at 0.75 1.1
\put{$\scriptstyle{\vert}\longleftarrow$}[l] at -3 -1.1
\put{$\longrightarrow\scriptstyle{\vert}$}[r] at 0 -1.1
\put{$\scriptstyle{d_i^-(p)=d^-(\bar p)}$} at -1.5 -1.1
\plot 0.5 -0.7  0 -0.7 /
\plot 0 -0.7  0 -0.5 /
\arrow <5pt> [.2,.67] from 0 -0.5 to 0.45 -0.5 %
\arrow <5pt> [.2,.67] from 0.55 -0.5 to 1.45 -0.5 %
\arrow <5pt> [.2,.67] from 1.55 -0.5 to 2.45 -0.5 %
\arrow <5pt> [.2,.67] from 2.45 -0.3 to 1.55 -0.3 %
\arrow <5pt> [.2,.67] from 1.45 -0.3 to 0.55 -0.3 %
\arrow <5pt> [.2,.67] from 0.45 -0.3 to -0.45 -0.3 %
\arrow <5pt> [.2,.67] from -0.55 -0.3 to -1.45 -0.3 %
\plot -1.55 -0.3  -2 -0.3 /
\plot -2 -0.3  -2 -0.1 /
\arrow <5pt> [.2,.67] from -2 -0.1 to -1.55 -0.1 %
\arrow <5pt> [.2,.67] from -1.45 -0.1 to -0.55 -0.1 %
\arrow <5pt> [.2,.67] from -0.45 -0.1 to 0.45 -0.1 %
\arrow <5pt> [.2,.67] from 0.45 0.1 to -0.45 0.1 %
\arrow <5pt> [.2,.67] from -0.55 0.1 to -1.45 0.1 %
\arrow <5pt> [.2,.67] from -1.55 0.1 to -2.45 0.1 %
\plot -2.45 0.1  -3 0.1 /
\plot -3 0.1  -3 0.3 /
\arrow <5pt> [.2,.67] from -3 0.3 to -2.55 0.3 %
\arrow <5pt> [.2,.67] from -2.45 0.3 to -1.55 0.3 %
\arrow <5pt> [.2,.67] from -1.45 0.3 to -0.55 0.3 %
\arrow <5pt> [.2,.67] from -0.45 0.3 to 0.45 0.3 %
\arrow <5pt> [.2,.67] from 0.55 0.3 to 1.45 0.3 %
\arrow <5pt> [.2,.67] from 1.55 0.3 to 2.45 0.3 %
\arrow <5pt> [.2,.67] from 2.55 0.3 to 3.45 0.3 %
\arrow <5pt> [.2,.67] from 3.45 0.5 to 2.55 0.5 %
\arrow <5pt> [.2,.67] from 2.45 0.5 to 1.55 0.5 %
\plot 1.45 0.5  1 0.5 /
\plot 1 0.5  1 0.7 /
\arrow <5pt> [.2,.67] from 1 0.7 to 1.45 0.7 %
\arrow <5pt> [.2,.67] from 1.55 0.7 to 2.45 0.7 %
\arrow <5pt> [.2,.67] from 2.55 0.7 to 3.45 0.7 %
\arrow <5pt> [.2,.67] from 3.55 0.7 to 4.45 0.7 %
\setplotsymbol({\bf .})                   
\plot 0.5 -0.7  0 -0.7 /
\plot 0 -0.3  -0.45 -0.3 /
\plot -0.55 -0.3  -1.45 -0.3 /
\plot -1.55 -0.3  -2 -0.3 /
\plot -2 0.1  -2.45 0.1 /
\plot -2.55 0.1  -3 0.1 /
\plot -3 0.3  -2.55 0.3 /
\plot -2.45 0.3  -1.55 0.3 /
\plot -1.45 0.3  -0.55 0.3 /
\plot -0.45 0.3  0.45 0.3 /
\plot 0.55 0.3  1 0.3 /
\plot 1 0.7  1.45 0.7 /
\plot 1.55 0.7  2.45 0.7 /
\plot 2.55 0.7  3.45 0.7 /
\plot 3.55 0.7  4.45 0.7 /
\endpicture
\qquad\quad
\end{equation}
Because $\bar p$ is positively folded it
is a concatenation of negative-positive sections of the form
$c^-c^-\cdots c^-fc^+c^+\cdots c^+$, where $c^+$ denotes a positive crossing,
$c^-$ a negative crossing, and $f$ a (positive) fold.  
The \emph{outer edge} (bottom most negative traveling portion 
and topmost positive traveling portion) of the walk is a single negative-positive walk
\begin{equation}
\quad \underbrace{c^-c^-\cdots c^-c^-}_{d^-(\bar p)\ \mathrm{factors}}f
\underbrace{c^+c^+\cdots c^+c^+}_{d^+(\bar p)\ \mathrm{factors}}. 
\qquad\qquad
\beginpicture
\setcoordinatesystem units <0.75cm,0.75cm>         
\setplotarea x from -3.5 to 3.5, y from -0.5 to 0.5  
\setplotsymbol({\bf .})                   
\arrow <5pt> [.2,.67] from 0.45 -0.1 to -0.45 -0.1 %
\arrow <5pt> [.2,.67] from -0.55 -0.1 to -1.45 -0.1 %
\arrow <5pt> [.2,.67] from -1.55 -0.1 to -2.45 -0.1 %
\plot -2.55 -0.1  -3 -0.1 /
\plot  -3 -0.1  -3 0.1 / 
\arrow <5pt> [.2,.67] from -3 0.1 to -2.55 0.1 %
\arrow <5pt> [.2,.67] from -2.45 0.1 to -1.55 0.1 %
\arrow <5pt> [.2,.67] from -1.45 0.1 to -0.55 0.1 %
\arrow <5pt> [.2,.67] from -0.45 0.1 to 0.45 0.1  %
\arrow <5pt> [.2,.67] from 0.55 0.1 to 1.45 0.1  %
\arrow <5pt> [.2,.67] from 1.55 0.1 to 2.45 0.1  %
\arrow <5pt> [.2,.67] from 2.55 0.1 to 3.45 0.1  %
\arrow <5pt> [.2,.67] from 3.55 0.1 to 4.45 0.1  %
\endpicture
\end{equation}
If $\bar p$ is the walk in (4.16), the outer edge is the darkened portion of the path,
$d^+(\bar p)=7$ and $d^-(\bar p)=3$.

The \emph{root operators} $\tilde e$ and $\tilde f$
change the outer edge of the path $\bar p$ and leave all other parts of the
walk unchanged.  Define $\tilde e \bar p$ and $\tilde f \bar p$ to be the
positively folded alcove walks which are the same as $\bar p$ except that 
\begin{align*}
\hbox{$\tilde e \bar p$ has outer edge}\quad 
&\underbrace{c^-c^-\cdots c^-c^-c^-}_{d^-(\bar p)+1\ \mathrm{factors}}f
\underbrace{c^+\cdots c^+c^+}_{d^+(\bar p)-1\ \mathrm{factors}}, \ \ \hbox{and} \\
\hbox{$\tilde f \bar p$ has outer edge}\quad 
&\quad \underbrace{c^-c^-\cdots c^-}_{d^-(\bar p)-1\ \mathrm{factors}}f
\underbrace{c^+c^+c^+\cdots c^+c^+}_{d^+(\bar p)+1\ \mathrm{factors}} 
\end{align*}
If $\bar p$ is the walk in (4.16) then
$$
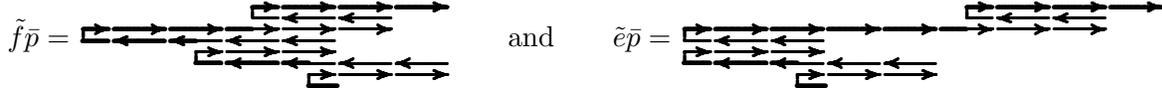
\tilde f \bar p =
\beginpicture
\setcoordinatesystem units <0.75cm,0.75cm>         
\setplotarea x from -4.1 to 2.5, y from -0.8 to 0.8  
\plot 0.5 -0.7  0 -0.7 /
\plot 0 -0.7  0 -0.5 /
\arrow <5pt> [.2,.67] from 0 -0.5 to 0.45 -0.5 %
\arrow <5pt> [.2,.67] from 0.55 -0.5 to 1.45 -0.5 %
\arrow <5pt> [.2,.67] from 1.55 -0.5 to 2.45 -0.5 %
\arrow <5pt> [.2,.67] from 2.45 -0.3 to 1.55 -0.3 %
\arrow <5pt> [.2,.67] from 1.45 -0.3 to 0.55 -0.3 %
\arrow <5pt> [.2,.67] from 0.45 -0.3 to -0.45 -0.3 %
\arrow <5pt> [.2,.67] from -0.55 -0.3 to -1.45 -0.3 %
\plot -1.55 -0.3  -2 -0.3 /
\plot -2 -0.3  -2 -0.1 /
\arrow <5pt> [.2,.67] from -2 -0.1 to -1.55 -0.1 %
\arrow <5pt> [.2,.67] from -1.45 -0.1 to -0.55 -0.1 %
\arrow <5pt> [.2,.67] from -0.45 -0.1 to 0.45 -0.1 %
\arrow <5pt> [.2,.67] from 0.45 0.1 to -0.45 0.1 %
\arrow <5pt> [.2,.67] from -0.55 0.1 to -1.45 0.1 %
\arrow <5pt> [.2,.67] from -1.55 0.1 to -2.45 0.1 %
\arrow <5pt> [.2,.67] from -2.55 0.1 to -3.45 0.1 %
\plot -3.55 0.1  -4 0.1 /
\plot -4 0.1  -4 0.3 /
\arrow <5pt> [.2,.67] from -4 0.3 to -3.55 0.3 %
\arrow <5pt> [.2,.67] from -3.45 0.3 to -2.55 0.3 %
\arrow <5pt> [.2,.67] from -2.45 0.3 to -1.55 0.3 %
\arrow <5pt> [.2,.67] from -1.45 0.3 to -0.55 0.3 %
\arrow <5pt> [.2,.67] from -0.45 0.3 to 0.45 0.3 %
\arrow <5pt> [.2,.67] from 0.55 0.3 to 1.45 0.3 %
\arrow <5pt> [.2,.67] from 1.45 0.5 to 0.55 0.5 %
\arrow <5pt> [.2,.67] from 0.45 0.5 to -0.45 0.5 %
\plot -0.55 0.5  -1 0.5 /
\plot -1 0.5  -1 0.7 /
\arrow <5pt> [.2,.67] from -1 0.7 to -0.55 0.7 %
\arrow <5pt> [.2,.67] from -0.45 0.7 to 0.45 0.7 %
\arrow <5pt> [.2,.67] from 0.55 0.7 to 1.45 0.7 %
\arrow <5pt> [.2,.67] from 1.55 0.7 to 2.45 0.7 %
\setplotsymbol({\bf .})                   
\plot 0.5 -0.7  0 -0.7 /
\plot 0 -0.3  -0.45 -0.3 /
\plot -0.55 -0.3  -1.45 -0.3 /
\plot -1.55 -0.3  -2 -0.3 /
\plot -2 0.1  -2.45 0.1 /
\plot -2.55 0.1  -3.45 0.1 /
\plot -3.55 0.1 -4 0.1 /
\plot -4 0.3  -3.55 0.3 /
\plot -3.45 0.3  -2.55 0.3 /
\plot -2.45 0.3  -1.55 0.3 /
\plot -1.45 0.3  -1 0.3 /
\plot -1 0.7  -0.55 0.7 /
\plot -0.45 0.7  0.45 0.7 /
\plot 0.55 0.7  1.45 0.7 /
\plot 1.55 0.7  2.45 0.7 /
\endpicture
\qquad\hbox{and}\qquad 
\tilde e \bar p =
\beginpicture
\setcoordinatesystem units <0.75cm,0.75cm>         
\setplotarea x from -2.1 to 6.5, y from -0.8 to 0.8  
\plot 0.5 -0.7  0 -0.7 /
\plot 0 -0.7  0 -0.5 /
\arrow <5pt> [.2,.67] from 0 -0.5 to 0.45 -0.5 %
\arrow <5pt> [.2,.67] from 0.55 -0.5 to 1.45 -0.5 %
\arrow <5pt> [.2,.67] from 1.55 -0.5 to 2.45 -0.5 %
\arrow <5pt> [.2,.67] from 2.45 -0.3 to 1.55 -0.3 %
\arrow <5pt> [.2,.67] from 1.45 -0.3 to 0.55 -0.3 %
\arrow <5pt> [.2,.67] from 0.45 -0.3 to -0.45 -0.3 %
\arrow <5pt> [.2,.67] from -0.55 -0.3 to -1.45 -0.3 %
\plot -1.55 -0.3  -2 -0.3 /
\plot -2 -0.3  -2 -0.1 /
\arrow <5pt> [.2,.67] from -2 -0.1 to -1.55 -0.1 %
\arrow <5pt> [.2,.67] from -1.45 -0.1 to -0.55 -0.1 %
\arrow <5pt> [.2,.67] from -0.45 -0.1 to 0.45 -0.1 %
\arrow <5pt> [.2,.67] from 0.45 0.1 to -0.45 0.1 %
\arrow <5pt> [.2,.67] from -0.55 0.1 to -1.45 0.1 %
\plot -1.55 0.1  -2 0.1 /
\plot -2 0.1  -2 0.3 /
\arrow <5pt> [.2,.67] from -2 0.3 to -1.55 0.3 %
\arrow <5pt> [.2,.67] from -1.45 0.3 to -0.55 0.3 %
\arrow <5pt> [.2,.67] from -0.45 0.3 to 0.45 0.3 %
\arrow <5pt> [.2,.67] from 0.55 0.3 to 1.45 0.3 %
\arrow <5pt> [.2,.67] from 1.55 0.3 to 2.45 0.3 %
\arrow <5pt> [.2,.67] from 2.55 0.3 to 3.45 0.3 %
\arrow <5pt> [.2,.67] from 3.55 0.3 to 4.45 0.3 %
\arrow <5pt> [.2,.67] from 4.55 0.3 to 5.45 0.3 %
\arrow <5pt> [.2,.67] from 5.45 0.5 to 4.55 0.5 %
\arrow <5pt> [.2,.67] from 4.45 0.5 to 3.55 0.5 %
\plot 3.45 0.5  3 0.5 /
\plot 3 0.5  3 0.7 /
\arrow <5pt> [.2,.67] from 3 0.7 to 3.45 0.7 %
\arrow <5pt> [.2,.67] from 3.55 0.7 to 4.45 0.7 %
\arrow <5pt> [.2,.67] from 4.55 0.7 to 5.45 0.7 %
\arrow <5pt> [.2,.67] from 5.55 0.7 to 6.45 0.7 %
\setplotsymbol({\bf .})                   
\plot 0.5 -0.7  0 -0.7 /
\plot 0 -0.3  -0.45 -0.3 /
\plot -0.55 -0.3  -1.45 -0.3 /
\plot -1.55 -0.3  -2 -0.3 /
\plot -2 0.3  -1.55 0.3 /
\plot -1.45 0.3  -0.55 0.3 /
\plot -0.45 0.3  0.45 0.3 /
\plot 0.55 0.3  1.45 0.3 /
\plot 1.55 0.3  2.45 0.3 /
\plot 2.55 0.3  3 0.3 /
\plot 3 0.7  3.45 0.7 /
\plot 3.55 0.7  4.45 0.7 /
\plot 4.55 0.7  5.45 0.7 /
\plot 5.55 0.7  6.45 0.7 /
\endpicture
$$
The precise rules for the limiting cases, when $d^+(\bar p)$ or $d^-(\bar p)=0$,
are illustrated by the following example, where the dashed arrows indicate
the action of $\tilde e$ and $\tilde f$.
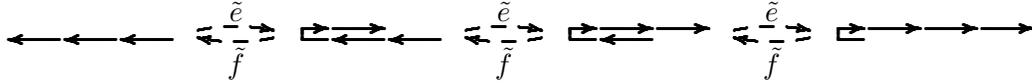
$$
\beginpicture
\setcoordinatesystem units <0.75cm,0.75cm>         
\setplotarea x from -2.5 to 0.5, y from -0.1 to 0.1  
\arrow <5pt> [.2,.67] from 0.45 -0.1 to -0.45 -0.1 %
\arrow <5pt> [.2,.67] from -0.55 -0.1 to -1.45 -0.1 %
\arrow <5pt> [.2,.67] from -1.55 -0.1 to -2.45 -0.1 %
\endpicture
\ \ 
\beginpicture
\setcoordinatesystem units <0.5cm,0.5cm>         
\setplotarea x from -1.1 to 1.1, y from -0.3 to 0.3  
\put{$\tilde e$}[b] at 0 0.35 
\put{$\tilde f$}[t] at 0 -0.35 
\setdashes
\setquadratic \plot 1 -0.1                      
                      0 -0.2                     %
                     -1 -0.1  /                   %
\arrow <5pt> [.2,.67] from -0.95 -0.11  to -1 -0.1   %
\setquadratic \plot -1 0.1                      
                      0 0.2                     %
                     1 0.1  /                   %
\arrow <5pt> [.2,.67] from 0.95 0.11  to 1 0.1   %
\endpicture
\ \ 
\beginpicture
\setcoordinatesystem units <0.75cm,0.75cm>         
\setplotarea x from -2.1 to 0.5, y from -0.2 to 0.2  
\arrow <5pt> [.2,.67] from 0.45 -0.1 to -0.45 -0.1 %
\arrow <5pt> [.2,.67] from -0.55 -0.1 to -1.45 -0.1 %
\plot -1.55 -0.1  -2 -0.1 /
\plot  -2 -0.1  -2 0.1 / 
\arrow <5pt> [.2,.67] from -2 0.1 to -1.55 0.1 %
\arrow <5pt> [.2,.67] from -1.45 0.1 to -0.55 0.1 %
\endpicture
\ \ 
\beginpicture
\setcoordinatesystem units <0.5cm,0.5cm>         
\setplotarea x from -1.1 to 1.1, y from -0.3 to 0.3  
\put{$\tilde e$}[b] at 0 0.35 
\put{$\tilde f$}[t] at 0 -0.35 
\setdashes
\setquadratic \plot 1 -0.1                      
                      0 -0.2                     %
                     -1 -0.1  /                   %
\arrow <5pt> [.2,.67] from -0.95 -0.11  to -1 -0.1   %
\setquadratic \plot -1 0.1                      
                      0 0.2                     %
                     1 0.1  /                   %
\arrow <5pt> [.2,.67] from 0.95 0.11  to 1 0.1   %
\endpicture
\ \ 
\beginpicture
\setcoordinatesystem units <0.75cm,0.75cm>         
\setplotarea x from -1.1 to 1.5, y from -0.2 to 0.2  
\arrow <5pt> [.2,.67] from 0.45 -0.1 to -0.45 -0.1 %
\plot -0.55 -0.1  -1 -0.1 /
\plot  -1 -0.1  -1 0.1 / 
\arrow <5pt> [.2,.67] from -1 0.1 to -0.55 0.1 %
\arrow <5pt> [.2,.67] from -0.45 0.1 to 0.45 0.1 %
\arrow <5pt> [.2,.67] from 0.55 0.1 to 1.45 0.1 %
\endpicture
\ \ 
\beginpicture
\setcoordinatesystem units <0.5cm,0.5cm>         
\setplotarea x from -1.1 to 1.1, y from -0.3 to 0.3  
\put{$\tilde e$}[b] at 0 0.35 
\put{$\tilde f$}[t] at 0 -0.35 
\setdashes
\setquadratic \plot 1 -0.1                      
                      0 -0.2                     %
                     -1 -0.1  /                   %
\arrow <5pt> [.2,.67] from -0.95 -0.11  to -1 -0.1   %
\setquadratic \plot -1 0.1                      
                      0 0.2                     %
                     1 0.1  /                   %
\arrow <5pt> [.2,.67] from 0.95 0.11  to 1 0.1   %
\endpicture
\ \ 
\beginpicture
\setcoordinatesystem units <0.75cm,0.75cm>         
\setplotarea x from -0.1 to 3.5, y from -0.1 to 0.1  
\plot 0.45 -0.1  0 -0.1 /
\plot  0 -0.1  0 0.1 / 
\arrow <5pt> [.2,.67] from 0 0.1 to 0.45 0.1 %
\arrow <5pt> [.2,.67] from 0.55 0.1 to 1.45 0.1 %
\arrow <5pt> [.2,.67] from 1.55 0.1 to 2.45 0.1 %
\arrow <5pt> [.2,.67] from 2.55 0.1 to 3.45 0.1 %
\endpicture
$$
with
$\tilde f(
\beginpicture
\setcoordinatesystem units <0.75cm,0.75cm>         
\setplotarea x from -2.5 to 0.5, y from -0.1 to 0.1  
\arrow <5pt> [.2,.67] from 0.45 0.1 to -0.45 0.1 %
\arrow <5pt> [.2,.67] from -0.55 0.1 to -1.45 0.1 %
\arrow <5pt> [.2,.67] from -1.55 0.1 to -2.45 0.1 %
\endpicture
)=0$ 
and
$\tilde e(
\beginpicture
\setcoordinatesystem units <0.75cm,0.75cm>         
\setplotarea x from -0.1 to 3.5, y from -0.1 to 0.1  
\plot 0.45 -0.1  0 -0.1 /
\plot  0 -0.1  0 0.1 / 
\arrow <5pt> [.2,.67] from 0 0.1 to 0.45 0.1 %
\arrow <5pt> [.2,.67] from 0.55 0.1 to 1.45 0.1 %
\arrow <5pt> [.2,.67] from 1.55 0.1 to 2.45 0.1 %
\arrow <5pt> [.2,.67] from 2.55 0.1 to 3.45 0.1 %
\endpicture
)=0$.

With notations for $d^+(\bar p)$ and $d^-(\bar p)$ as in (4.18) and (4.19), define
\begin{equation}
d_i^+(p) = d^+(\bar p)
\qquad\hbox{and}\qquad
d_i^-(p) = d^-(\bar p),
\end{equation}
where $\bar p$ is the projection of $p$ onto the line perpendicular to $\alpha_i$.
The walks $\tilde e_i p$ and $\tilde f_i p$ are the walks obtained from $p$
by changing the corresponding edges $p$ (so that the projections of
$\tilde e_i p$ and $\tilde f_i p$ are $\tilde e \bar p$ and $\tilde f \bar p$, respectively).
$$
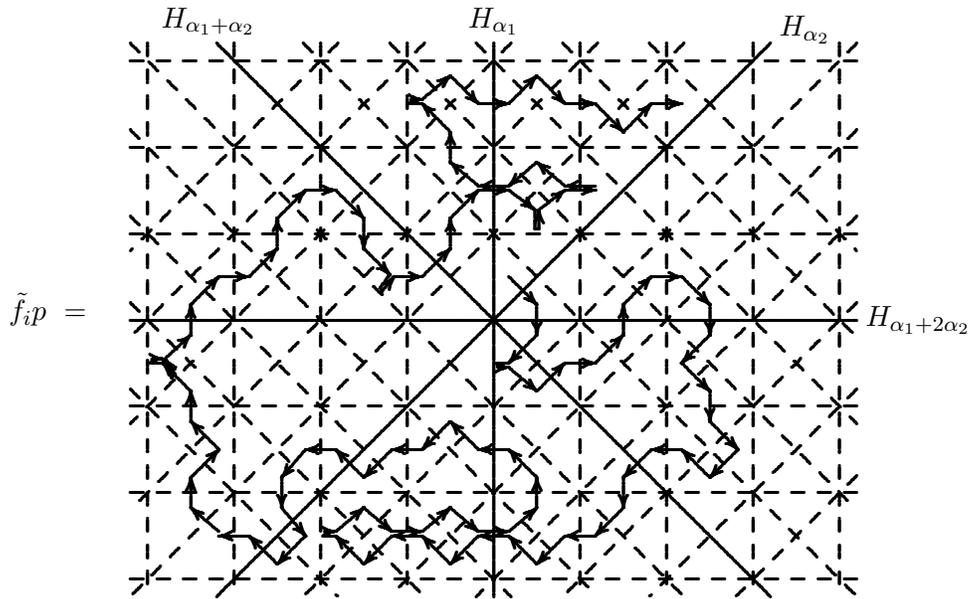
\tilde f_i p ~=~ 
\beginpicture
\setcoordinatesystem units <1.15cm,1.15cm>         
\setplotarea x from -4.5 to 4.5, y from -3.5 to 4    
\put{$H_{\alpha_1}$}[b] at 0 3.3
\put{$H_{\alpha_2}$}[bl] at 3.3 3.2
\put{$H_{\alpha_1+\alpha_2}$}[b] at -3.3 3.3
\put{$H_{\alpha_1+2\alpha_2}$}[l] at 4.3 0
\plot -3.2 -3.2   3.2 3.2 /
\plot  3.2 -3.2  -3.2 3.2 /
\plot  0  3.2   0 -3.2 /
\plot  4.2  0  -4.2  0 /
\arrow <5pt> [.2,.67] from 0.175 0.5 to 0.5 0.175 %
\arrow <5pt> [.2,.67] from 0.5 0.175 to 0.5 -0.175 %
\arrow <5pt> [.2,.67] from 0.5 -0.175 to 0.175 -0.5 %
\plot 0.175 -0.5  0.05 -0.5 /   
\plot 0.05 -0.5  0.05 -0.55 /  %
\arrow <5pt> [.2,.67] from 0.05 -0.55 to 0.175 -0.55   
\arrow <5pt> [.2,.67] from 0.175 -0.55 to 0.5 -0.825 %
\arrow <5pt> [.2,.67] from 0.5 -0.825 to 0.825 -0.5 %
\arrow <5pt> [.2,.67] from 0.825 -0.5 to 1.175 -0.5 %
\arrow <5pt> [.2,.67] from 1.175 -0.5 to 1.5 -0.175 %
\arrow <5pt> [.2,.67] from 1.5 -0.175 to 1.5 0.175 %
\arrow <5pt> [.2,.67] from 1.5 0.175 to 1.825 0.5 %
\arrow <5pt> [.2,.67] from 1.825 0.5 to 2.175 0.5 %
\arrow <5pt> [.2,.67] from 2.175 0.5 to 2.5 0.175 %
\arrow <5pt> [.2,.67] from 2.5 0.175 to 2.5 -0.175 %
\arrow <5pt> [.2,.67] from 2.5 -0.175 to 2.175 -0.5 %
\arrow <5pt> [.2,.67] from 2.175 -0.5 to 2.5 -0.825 %
\arrow <5pt> [.2,.67] from 2.5 -0.825 to 2.5 -1.175 %
\arrow <5pt> [.2,.67] from 2.5 -1.175 to 2.825 -1.5 %
\arrow <5pt> [.2,.67] from 2.825 -1.5 to 2.5 -1.825 %
\arrow <5pt> [.2,.67] from 2.5 -1.825 to 2.175 -1.5 %
\arrow <5pt> [.2,.67] from 2.175 -1.5 to 1.825 -1.5 %
\arrow <5pt> [.2,.67] from 1.825 -1.5 to 1.5 -1.825 %
\arrow <5pt> [.2,.67] from 1.5 -1.825 to 1.5 -2.175 %
\arrow <5pt> [.2,.67] from 1.5 -2.175 to 1.175 -2.5 %
\arrow <5pt> [.2,.67] from 1.175 -2.5 to 0.825 -2.5 %
\arrow <5pt> [.2,.67] from 0.825 -2.5 to 0.5 -2.825 %
\arrow <5pt> [.2,.67] from 0.5 -2.825 to 0.175 -2.5 %
\arrow <5pt> [.2,.67] from 0.175 -2.5 to -0.175 -2.5 %
\arrow <5pt> [.2,.67] from -0.175 -2.5 to -0.5 -2.825 %
\arrow <5pt> [.2,.67] from -0.5 -2.825 to -0.825 -2.5 %
\arrow <5pt> [.2,.67] from -0.825 -2.5 to -1.175 -2.5 %
\arrow <5pt> [.2,.67] from -1.175 -2.5 to -1.5 -2.825 %
\arrow <5pt> [.2,.67] from -1.5 -2.825 to -1.825 -2.5 %
\plot -1.825 -2.5  -1.95 -2.5 /   
\plot -1.95 -2.5  -1.95 -2.45 /  %
\arrow <5pt> [.2,.67] from -1.95 -2.45 to -1.825 -2.45   
\arrow <5pt> [.2,.67] from -1.825 -2.45 to -1.5 -2.175 %
\arrow <5pt> [.2,.67] from -1.5 -2.175 to -1.175 -2.45 %
\arrow <5pt> [.2,.67] from -1.175 -2.45 to -0.825 -2.45 %
\arrow <5pt> [.2,.67] from -0.825 -2.45 to -0.5 -2.175 %
\arrow <5pt> [.2,.67] from -0.5 -2.175 to -0.175 -2.45 %
\arrow <5pt> [.2,.67] from -0.175 -2.45 to 0.175 -2.45 %
\arrow <5pt> [.2,.67] from 0.175 -2.45 to 0.5 -2.175 %
\arrow <5pt> [.2,.67] from 0.5 -2.175 to 0.5 -1.825 %
\arrow <5pt> [.2,.67] from 0.5 -1.825 to 0.175 -1.5 %
\arrow <5pt> [.2,.67] from 0.175 -1.5 to -0.175 -1.5 %
\arrow <5pt> [.2,.67] from -0.175 -1.5 to -0.5 -1.175 %
\arrow <5pt> [.2,.67] from -0.5 -1.175 to -0.825 -1.5 %
\arrow <5pt> [.2,.67] from -0.825 -1.5 to -1.175 -1.5 %
\arrow <5pt> [.2,.67] from -1.175 -1.5 to -1.5 -1.825 %
\arrow <5pt> [.2,.67] from -1.5 -1.825 to -1.825 -1.5 %
\arrow <5pt> [.2,.67] from -1.825 -1.5 to -2.175 -1.5 %
\arrow <5pt> [.2,.67] from -2.175 -1.5 to -2.45 -1.825 %
\arrow <5pt> [.2,.67] from -2.45 -1.825 to -2.45 -2.175 %
\arrow <5pt> [.2,.67] from -2.45 -2.175 to -2.175 -2.5 %
\arrow <5pt> [.2,.67] from -2.175 -2.5 to -2.5 -2.825 %
\arrow <5pt> [.2,.67] from -2.5 -2.825 to -2.825 -2.5 %
%
\arrow <5pt> [.2,.67] from -2.825 -2.5 to -3.175 -2.5 %
\arrow <5pt> [.2,.67] from -3.175 -2.45 to -3.5 -2.175 %
\arrow <5pt> [.2,.67] from -3.5 -2.175 to -3.5 -1.825 %
\arrow <5pt> [.2,.67] from -3.5 -1.825 to -3.175 -1.5 %
\arrow <5pt> [.2,.67] from -3.175 -1.5 to -3.5 -1.175 %
\arrow <5pt> [.2,.67] from -3.5 -1.175 to -3.5 -0.825 %
\arrow <5pt> [.2,.67] from -3.5 -0.825 to -3.825 -0.5 %
\plot -3.825 -0.5  -3.95 -0.5 /   
\plot -3.95 -0.5  -3.95 -0.45 /  %
\arrow <5pt> [.2,.67] from -3.95 -0.45 to -3.825 -0.45   
%
\plot -3.825 -0.45   -3.72 -0.67  /  %
\plot -3.72  -0.67  -3.684 -0.625  / %
\arrow <5pt> [.2,.67] from -3.684 -0.625 to -3.8 -0.45   %
\arrow <5pt> [.2,.67] from -3.8 -0.45 to -3.5 -0.175 %
\arrow <5pt> [.2,.67] from -3.5 -0.175 to -3.5 0.175 %
\arrow <5pt> [.2,.67] from -3.5 0.175 to -3.175 0.5 %
\arrow <5pt> [.2,.67] from -3.175 0.5 to -2.825 0.5 %
\arrow <5pt> [.2,.67] from -2.825 0.5 to -2.5 0.825 %
\arrow <5pt> [.2,.67] from -2.5 0.825 to -2.5 1.175 %
\arrow <5pt> [.2,.67] from -2.5 1.175 to -2.175 1.5 %
\arrow <5pt> [.2,.67] from -2.175 1.5 to -1.825 1.5 %
\arrow <5pt> [.2,.67] from -1.825 1.5 to -1.5 1.175 %
\arrow <5pt> [.2,.67] from -1.5 1.175 to -1.5 0.825 %
\arrow <5pt> [.2,.67] from -1.5 0.825 to -1.2 0.525 %
\plot  -1.2 0.525   -1.31 0.375  /  %
\plot -1.31  0.375  -1.291 0.33  / %
\arrow <5pt> [.2,.67] from -1.291 0.33 to -1.15 0.5   %
%
%
%
\arrow <5pt> [.2,.67] from -1.15 0.5 to -0.825 0.5 %
\arrow <5pt> [.2,.67] from -0.825 0.5 to -0.5 0.825 %
\arrow <5pt> [.2,.67] from -0.5 0.825 to -0.5 1.175 %
\arrow <5pt> [.2,.67] from -0.5 1.175 to -0.175 1.5 %
\arrow <5pt> [.2,.67] from -0.175 1.5 to 0.175 1.5 %
\arrow <5pt> [.2,.67] from 0.175 1.5 to 0.475 1.275 %
\plot  0.475 1.275   0.475 1.05  /  %
\plot 0.475  1.05  0.525 1.05  / %
\arrow <5pt> [.2,.67] from 0.525 1.05 to 0.525 1.275   %
\arrow <5pt> [.2,.67] from 0.525 1.275 to 0.825 1.5 %
\arrow <5pt> [.2,.67] from 0.825 1.5 to 1.175 1.5 %
\arrow <5pt> [.2,.67] from 1.175 1.55 to 0.825 1.55 %
\arrow <5pt> [.2,.67] from 0.825 1.55 to 0.5 1.825 %
\arrow <5pt> [.2,.67] from 0.5 1.825 to 0.175 1.55 %
\arrow <5pt> [.2,.67] from 0.175 1.55 to -0.175 1.55 %
\arrow <5pt> [.2,.67] from -0.175 1.55 to -0.5 1.825 %
\arrow <5pt> [.2,.67] from -0.5 1.825 to -0.5 2.175 %
\arrow <5pt> [.2,.67] from -0.5 2.175 to -0.825 2.5 %
\plot -0.825 2.5  -0.95 2.5 /   
\plot -0.95 2.5  -0.95 2.55 /  %
\arrow <5pt> [.2,.67] from -0.95 2.55 to -0.825 2.55   
\arrow <5pt> [.2,.67] from -0.825 2.55 to -0.5 2.825 %
\arrow <5pt> [.2,.67] from -0.5 2.825 to -0.175 2.5 %
\arrow <5pt> [.2,.67] from -0.175 2.5 to 0.175 2.5 %
\arrow <5pt> [.2,.67] from 0.175 2.5 to 0.5 2.825 %
\arrow <5pt> [.2,.67] from 0.5 2.825 to 0.825 2.5 %
\arrow <5pt> [.2,.67] from 0.825 2.5 to 1.175 2.5 %
\arrow <5pt> [.2,.67] from 1.175 2.5 to 1.5 2.175 %
\arrow <5pt> [.2,.67] from 1.5 2.175 to 1.825 2.5 %
\arrow <5pt> [.2,.67] from 1.825 2.5 to 2.175 2.5 %
\setdashes
\plot  -4 3.2    -4 -3.2 /
\plot  -3 3.2   -3 -3.2 /
\plot  -2 3.2   -2 -3.2 /
\plot  -1 3.2   -1 -3.2 /
\plot   1 3.2    1 -3.2 /
\plot   2 3.2    2 -3.2 /
\plot   3 3.2    3 -3.2 /
\plot   4 3.2    4 -3.2 /
\plot  4.2 -3   -4.2 -3 /
\plot  4.2 -2   -4.2 -2 /
\plot  4.2 -1   -4.2 -1 /
\plot  4.2  1   -4.2  1 /
\plot  4.2  2   -4.2  2 /
\plot  4.2  3   -4.2  3 /
\plot  4.2 -1.8   2.8 -3.2 /
\plot  4.2 -0.8   1.8 -3.2 /
\plot  4.2 0.2   0.8 -3.2 /
\plot  4.2 1.2    -0.2 -3.2 /
\plot  4.2 2.2   -1.2 -3.2 /
\plot  4.2 3.2   -2.2 -3.2 /
\plot  2.2 3.2   -4.2 -3.2 /
\plot  1.2 3.2   -4.2 -2.2 /
\plot  0.2 3.2   -4.2 -1.2 /
\plot  -0.8 3.2  -4.2 -0.2 /
\plot  -1.8 3.2  -4.2 0.8 /
\plot  -2.2 3.2   4.2 -3.2 /
\plot  -1.2 3.2   4.2 -2.2 /
\plot  -0.2 3.2   4.2 -1.2 /
\plot  0.8 3.2  4.2 -0.2 /
\plot  1.8 3.2  4.2 0.8 /
\plot  -4.2 -0.8   -1.8 -3.2 /
\plot  -4.2 0.2   -0.8 -3.2 /
\plot  -4.2 1.2    0.2 -3.2 /
\plot  -4.2 2.2   1.2 -3.2 /
\plot  -4.2 3.2   2.2 -3.2 /
\endpicture
$$
$$
\tilde e_i p ~=~ 
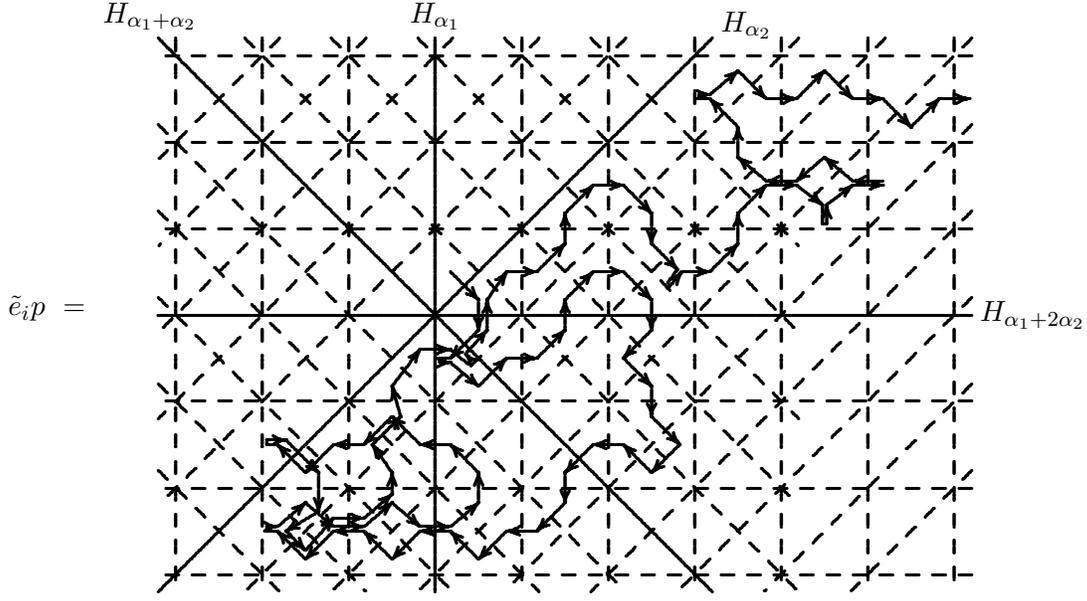
\beginpicture
\setcoordinatesystem units <1.15cm,1.15cm>         
\setplotarea x from -3.0 to 6.5, y from -3.5 to 4    
\put{$H_{\alpha_1}$}[b] at 0 3.3
\put{$H_{\alpha_2}$}[bl] at 3.3 3.2
\put{$H_{\alpha_1+\alpha_2}$}[b] at -3.3 3.3
\put{$H_{\alpha_1+2\alpha_2}$}[l] at 6.3 0
\plot -3.2 -3.2   3.2 3.2 /
\plot  3.2 -3.2  -3.2 3.2 /
\plot  0  3.2   0 -3.2 /
\plot  6.2  0  -3.2  0 /
\arrow <5pt> [.2,.67] from 0.175 0.5 to 0.5 0.175 %
\arrow <5pt> [.2,.67] from 0.5 0.175 to 0.5 -0.175 %
\arrow <5pt> [.2,.67] from 0.5 -0.175 to 0.175 -0.5 %
\plot 0.175 -0.5  0.05 -0.5 /   
\plot 0.05 -0.5  0.05 -0.55 /  %
\arrow <5pt> [.2,.67] from 0.05 -0.55 to 0.175 -0.55   
\arrow <5pt> [.2,.67] from 0.175 -0.55 to 0.5 -0.825 %
\arrow <5pt> [.2,.67] from 0.5 -0.825 to 0.825 -0.5 %
\arrow <5pt> [.2,.67] from 0.825 -0.5 to 1.175 -0.5 %
\arrow <5pt> [.2,.67] from 1.175 -0.5 to 1.5 -0.175 %
\arrow <5pt> [.2,.67] from 1.5 -0.175 to 1.5 0.175 %
\arrow <5pt> [.2,.67] from 1.5 0.175 to 1.825 0.5 %
\arrow <5pt> [.2,.67] from 1.825 0.5 to 2.175 0.5 %
\arrow <5pt> [.2,.67] from 2.175 0.5 to 2.5 0.175 %
\arrow <5pt> [.2,.67] from 2.5 0.175 to 2.5 -0.175 %
\arrow <5pt> [.2,.67] from 2.5 -0.175 to 2.175 -0.5 %
\arrow <5pt> [.2,.67] from 2.175 -0.5 to 2.5 -0.825 %
\arrow <5pt> [.2,.67] from 2.5 -0.825 to 2.5 -1.175 %
\arrow <5pt> [.2,.67] from 2.5 -1.175 to 2.825 -1.5 %
\arrow <5pt> [.2,.67] from 2.825 -1.5 to 2.5 -1.825 %
\arrow <5pt> [.2,.67] from 2.5 -1.825 to 2.175 -1.5 %
\arrow <5pt> [.2,.67] from 2.175 -1.5 to 1.825 -1.5 %
\arrow <5pt> [.2,.67] from 1.825 -1.5 to 1.5 -1.825 %
\arrow <5pt> [.2,.67] from 1.5 -1.825 to 1.5 -2.175 %
\arrow <5pt> [.2,.67] from 1.5 -2.175 to 1.175 -2.5 %
\arrow <5pt> [.2,.67] from 1.175 -2.5 to 0.825 -2.5 %
\arrow <5pt> [.2,.67] from 0.825 -2.5 to 0.5 -2.825 %
\arrow <5pt> [.2,.67] from 0.5 -2.825 to 0.175 -2.5 %
\arrow <5pt> [.2,.67] from 0.175 -2.5 to -0.175 -2.5 %
\arrow <5pt> [.2,.67] from -0.175 -2.5 to -0.5 -2.825 %
\arrow <5pt> [.2,.67] from -0.5 -2.825 to -0.825 -2.5 %
\arrow <5pt> [.2,.67] from -0.825 -2.5 to -1.175 -2.5 %
\arrow <5pt> [.2,.67] from -1.175 -2.5 to -1.5 -2.825 %
\arrow <5pt> [.2,.67] from -1.5 -2.825 to -1.825 -2.5 %
\plot -1.825 -2.5  -1.95 -2.5 /   
\plot -1.95 -2.5  -1.95 -2.45 /  %
\arrow <5pt> [.2,.67] from -1.95 -2.45 to -1.825 -2.45   
\arrow <5pt> [.2,.67] from -1.825 -2.45 to -1.5 -2.175 %
\arrow <5pt> [.2,.67] from -1.5 -2.175 to -1.175 -2.45 %
\arrow <5pt> [.2,.67] from -1.175 -2.45 to -0.825 -2.45 %
\arrow <5pt> [.2,.67] from -0.825 -2.45 to -0.5 -2.175 %
\arrow <5pt> [.2,.67] from -0.5 -2.175 to -0.175 -2.45 %
\arrow <5pt> [.2,.67] from -0.175 -2.45 to 0.175 -2.45 %
\arrow <5pt> [.2,.67] from 0.175 -2.45 to 0.5 -2.175 %
\arrow <5pt> [.2,.67] from 0.5 -2.175 to 0.5 -1.825 %
\arrow <5pt> [.2,.67] from 0.5 -1.825 to 0.175 -1.5 %
\arrow <5pt> [.2,.67] from 0.175 -1.5 to -0.175 -1.5 %
\arrow <5pt> [.2,.67] from -0.175 -1.5 to -0.5 -1.175 %
\arrow <5pt> [.2,.67] from -0.5 -1.175 to -0.825 -1.5 %
\arrow <5pt> [.2,.67] from -0.825 -1.5 to -1.175 -1.5 %
\arrow <5pt> [.2,.67] from -1.175 -1.5 to -1.5 -1.825 %
\arrow <5pt> [.2,.67] from -1.5 -1.825 to -1.825 -1.5 %
\plot -1.825 -1.5  -1.95 -1.5 /   
\plot -1.95 -1.5  -1.95 -1.45 /  %
\arrow <5pt> [.2,.67] from -1.95 -1.45 to -1.725 -1.45   
%
\arrow <5pt> [.2,.67] from -1.725 -1.45 to -1.35 -1.825 %
\arrow <5pt> [.2,.67] from -1.35 -1.825 to -1.35 -2.275 %
\arrow <5pt> [.2,.67] from -1.35 -2.275 to -1.725 -2.5 %
\arrow <5pt> [.2,.67] from -1.725 -2.5 to -1.5 -2.725 %
\arrow <5pt> [.2,.67] from -1.5 -2.725 to -1.175 -2.35 %
\arrow <5pt> [.2,.67] from -1.175 -2.35 to -0.825 -2.35 %
%
\arrow <5pt> [.2,.67] from -0.825 -2.35 to -0.5 -2.075 %
\arrow <5pt> [.2,.67] from -0.5 -2.075 to -0.5 -1.825 %
\arrow <5pt> [.2,.67] from -0.5 -1.825 to -0.725 -1.5 %
\arrow <5pt> [.2,.67] from -0.725 -1.5 to -0.4 -1.175 %
\arrow <5pt> [.2,.67] from -0.4 -1.175 to -0.5 -0.825 %
\arrow <5pt> [.2,.67] from -0.5 -0.825 to -0.175 -0.4 %
\arrow <5pt> [.2,.67] from -0.175 -0.4 to  0.175 -0.4 %
\plot  0.175 -0.4   0.38 -0.57  /  %
\plot  0.38  -0.57   0.416 -0.525  / %
\arrow <5pt> [.2,.67] from  0.416 -0.525 to  0.35 -0.45   %
\arrow <5pt> [.2,.67] from  0.35 -0.45 to  0.6 -0.175 %
\arrow <5pt> [.2,.67] from 0.6 -0.175 to 0.6 0.175 %
\arrow <5pt> [.2,.67] from 0.6 0.175 to  0.825 0.5 %
\arrow <5pt> [.2,.67] from 0.825 0.5 to 1.175 0.5 %
\arrow <5pt> [.2,.67] from 1.175 0.5 to 1.5 0.825 %
\arrow <5pt> [.2,.67] from 1.5 0.825 to 1.5 1.175 %
\arrow <5pt> [.2,.67] from 1.5 1.175 to 1.825 1.5 %
\arrow <5pt> [.2,.67] from 1.825 1.5 to 2.175 1.5 %
\arrow <5pt> [.2,.67] from 2.175 1.5 to 2.5 1.175 %
\arrow <5pt> [.2,.67] from 2.5 1.175 to 2.5 0.825 %
\arrow <5pt> [.2,.67] from 2.5 0.825 to 2.8 0.525 %
\plot  2.8 0.525   2.69 0.375  /  %
\plot 2.69  0.375  2.709 0.33  / %
\arrow <5pt> [.2,.67] from 2.709 0.33 to 2.85 0.5   %
%
%
%
\arrow <5pt> [.2,.67] from 2.85 0.5 to 3.175 0.5 %
\arrow <5pt> [.2,.67] from 3.175 0.5 to 3.5 0.825 %
\arrow <5pt> [.2,.67] from 3.5 0.825 to 3.5 1.175 %
\arrow <5pt> [.2,.67] from 3.5 1.175 to 3.825 1.5 %
\arrow <5pt> [.2,.67] from 3.825 1.5 to 4.175 1.5 %
\arrow <5pt> [.2,.67] from 4.175 1.5 to 4.475 1.275 %
\plot  4.475 1.275   4.475 1.05  /  %
\plot 4.475  1.05  4.525 1.05  / %
\arrow <5pt> [.2,.67] from 4.525 1.05 to 4.525 1.275   %
\arrow <5pt> [.2,.67] from 4.525 1.275 to 4.825 1.5 %
\arrow <5pt> [.2,.67] from 4.825 1.5 to 5.175 1.5 %
\arrow <5pt> [.2,.67] from 5.175 1.55 to 4.825 1.55 %
\arrow <5pt> [.2,.67] from 4.825 1.55 to 4.5 1.825 %
\arrow <5pt> [.2,.67] from 4.5 1.825 to 4.175 1.55 %
\arrow <5pt> [.2,.67] from 4.175 1.55 to 3.825 1.55 %
\arrow <5pt> [.2,.67] from 3.825 1.55 to 3.5 1.825 %
\arrow <5pt> [.2,.67] from 3.5 1.825 to 3.5 2.175 %
\arrow <5pt> [.2,.67] from 3.5 2.175 to 3.175 2.5 %
\plot 3.175 2.5  3.05 2.5 /   
\plot 3.05 2.5  3.05 2.55 /  %
\arrow <5pt> [.2,.67] from 3.05 2.55 to 3.175 2.55   
\arrow <5pt> [.2,.67] from 3.175 2.55 to 3.5 2.825 %
\arrow <5pt> [.2,.67] from 3.5 2.825 to 3.825 2.5 %
\arrow <5pt> [.2,.67] from 3.825 2.5 to 4.175 2.5 %
\arrow <5pt> [.2,.67] from 4.175 2.5 to 4.5 2.825 %
\arrow <5pt> [.2,.67] from 4.5 2.825 to 4.825 2.5 %
\arrow <5pt> [.2,.67] from 4.825 2.5 to 5.175 2.5 %
\arrow <5pt> [.2,.67] from 5.175 2.5 to 5.5 2.175 %
\arrow <5pt> [.2,.67] from 5.5 2.175 to 5.825 2.5 %
\arrow <5pt> [.2,.67] from 5.825 2.5 to 6.175 2.5 %
\setdashes
\plot  -3 3.2   -3 -3.2 /
\plot  -2 3.2   -2 -3.2 /
\plot  -1 3.2   -1 -3.2 /
\plot   1 3.2    1 -3.2 /
\plot   2 3.2    2 -3.2 /
\plot   3 3.2    3 -3.2 /
\plot   4 3.2    4 -3.2 /
\plot   5 3.2    5 -3.2 /
\plot   6 3.2    6 -3.2 /
\plot  6.2 -3   -3.2 -3 /
\plot  6.2 -2   -3.2 -2 /
\plot  6.2 -1   -3.2 -1 /
\plot  6.2  1   -3.2  1 /
\plot  6.2  2   -3.2  2 /
\plot  6.2  3   -3.2  3 /
\plot  6.2 -1.8   4.8 -3.2 /
\plot  6.2 -0.8   3.8 -3.2 /
\plot  6.2 0.2   2.8 -3.2 /
\plot  6.2 1.2   1.8 -3.2 /
\plot  6.2 2.2   0.8 -3.2 /
\plot  6.2 3.2    -0.2 -3.2 /
\plot  5.2 3.2   -1.2 -3.2 /
\plot  4.2 3.2   -2.2 -3.2 /
\plot  2.2 3.2   -3.2 -2.2 /
\plot  1.2 3.2   -3.2 -1.2 /
\plot  0.2 3.2   -3.2 -0.2 /
\plot  -0.8 3.2  -3.2 0.8 /
\plot  -1.8 3.2  -3.2 1.8 /
\plot  -2.2 3.2   4.2 -3.2 /
\plot  -1.2 3.2   4.2 -2.2 /
\plot  -0.2 3.2   4.2 -1.2 /
\plot  0.8 3.2  4.2 -0.2 /
\plot  1.8 3.2  4.2 0.8 /
\plot  -3.2 -1.8   -1.8 -3.2 /
\plot  -3.2 -0.8   -0.8 -3.2 /
\plot  -3.2 0.2    0.2 -3.2 /
\plot  -3.2 1.2   1.2 -3.2 /
\plot  -3.2 2.2   2.2 -3.2 /
\endpicture
$$

The \emph{$i$-string of $p$} $S_i(p)$ is the set of paths generated from $p$ by applying
the root operators $\tilde e_i$ and $\tilde f_i$.  The \emph{head} of the $i$-string is
the path $h$ in $S_i(p)$ which has $d_i^-(h)=0$,  
If $\wt(h)=\lambda$ and $\lambda$ is on the positive side of $H_{\alpha_i}$ then
\begin{equation}
C_i X^\lambda
= q^{-1}X^{s_i\lambda}T_{s_i}^{-1} + 
(1-q^{-2})(X^{s_i\lambda}+X^{s_i\lambda+\alpha_i}+\cdots+X^{\lambda-\alpha_i})+X^\lambda,
\end{equation}
and the terms in this sum correspond to the paths in the $i$-string $S_i(h)$.

 \end{subsection}
 
 \begin{subsection}{$q$-Crystals}
 
 The $q$-crystals provide a combinatorial model for the spherical
 Hecke algebra $\mathbf{1}_0\tilde H\mathbf{1}_0$ in the basis 
 of Hall-Littlewood polynomials.   Three structural properties motivate
 the definition of $q$-crystals:
 \begin{enumerate}
 \item[(a)]
 The Hall-Littlewood polynomials are normalized versions of the 
 basis $\mathbf{1}_0X^\lambda\mathbf{1}_0$.  
 \item[(b)] The element
$\mathbf{1}_0$ is characterized by the property that
$C_i\mathbf{1}_0 = (1+q^{-2})\mathbf{1}_0$.
\item[(c)] The action of $C_i$ on $\tilde H\mathbf{1}_0\cong \KK[P]$
is captured in the combinatorics of $i$-strings.
\end{enumerate}
These properties indicate that the combinatorics of Hall-Littlewood polynomials
can be captured with the root operators.

Let 
\begin{equation}
\hbox{$B_{\mathrm{univ}}$ be the set of positively folded alcove walks which begin in the
$0$-polygon $WA$.}
\end{equation}    
If $B$ is a finite subset of $B_{\mathrm{univ}}$ the \emph{character} of $B$ is 
\begin{equation}
\mathrm{char}(B) = \sum_{p\in B} q^{-(\iota(p)+\varphi(p)-f(p))}(1-q^{-2})^{f(p)-c(p)}X^{\wt(p)}.
\end{equation} 
where $p$ has $f(p)$ folds,
$\iota(p)A$ is the alcove where $p$ begins, $\wt(p)+\varphi(p)A$ is the alcove where $p$
ends and
\begin{equation}
\hbox{$c(p)$ the number of folds of $p$ \emph{touching one of the hyperplanes 
$H_{\alpha_1},\ldots, H_{\alpha_n}$}.}
\end{equation}
A \emph{$q$-crystal} is a finite subset $B$ of $B_{\mathrm{univ}}$ which is 
closed under the action of the root operators.

A positively folded alcove walk is \emph{$i$-dominant} if it never touches 
the hyperplane $H_{\alpha_i,-1}$.  The head $h$ of an $i$-string 
$S_i(p)$ is $i$-dominant and $S_i(h)=S_i(p)$.
A positively folded alcove walk 
\begin{equation}
\hbox{$p$ is \emph{dominant} if $p\subseteq C-\rho$,}
\end{equation}
where 
$$C-\rho = \{ \mu\in \fh_\RR^*\ |\ \hbox{$\langle \mu,\alpha_i^\vee\rangle >-1$ for $1\le i\le n$}\}
$$
In other words, a positively folded alcove walk $p$ is dominant if it is $i$-dominant for all $i$, 
$1\le i\le n$.

\begin{thm}  Let $B$ be a $q$-crystal.  Then, with notations as in (4.23-4.25),
$$\mathrm{char}(B) = \sum_{p\in B\atop p\subseteq C-\rho} 
q^{-(\iota(p)+\varphi(p)-f(p))}(1-q^{-2})^{f(p)-c(p)}P_{\wt(p)}.$$
\end{thm}
\begin{proof}  If $p$ is a dominant walk let
$$
\hbox{$B_q(p)$ be the $q$-crystal generated by $p$}
$$
under the action of the root operators $\tilde e_i$ and $\tilde f_i$.
The point is that the set of all positively folded alcove walks
is partitioned into ``equivalence classes'' given by the sets $B_q(p)$ such that
$p\in B_{\mathrm{univ}}$ is dominant and $p$ is the unique dominant walk in $B_q(p)$.
Because $\mathbf{1}_0$ is characterized by the property that
$C_i\mathbf{1}_0 = (1+q^{-2})\mathbf{1}_0$ and
the action of $C_i$ is modeled by the combinatorics of $i$-strings (4.24),
this equivalence relation is generated by the relations $p\sim \tilde f_i p$ and $p\sim \tilde e_i p$.
\end{proof}

\end{subsection}

\begin{subsection}{Products and restrictions}

The results in this section are generalizations of the Littlewood-Richardson rules.  These
are obtained as corollaries of Theorem 4.8.  

The combinatorial definition of the root operators given above is essentially
a consequence of the \emph{Leibnitz rule} for the Demazure operator,
\begin{equation}
\Delta_i(fg) = \Delta_i(f)g + (s_if)(\Delta_ig),
\qquad\hbox{for $f,g\in \KK[P]$.}
\end{equation} 
The corresponding rule for the operators $C_i$ is 
\begin{equation}
C_i(fg) = (C_if)g +(s_if)((C_i-(1+q^{-2})) g),
\qquad\hbox{for $f,g\in \KK[P]$.}
\end{equation}
This identity is implicit in the additivity in $\lambda$ of the relation
in Proposition 3.2e (the product $T_{s_i}X^{\lambda+\mu} = (T_{s_i} X^\lambda) X^\mu$
can be expanded in two different ways using Proposition 3.2e).

In order to define the product of $p_1\otimes p_2$ of walks 
$p_1,p_2\in B_{\mathrm{univ}}$ the final direction $\varphi(p_1)$ of $p_1$ and
the initial direction $\iota(p_2)$ of $p_2$ need to be taken into account. 
(To properly model the multiplication of Hall-Littlewood polynomials
we must account for the effect of $\mathbf{1}_0$ in the product
$P_\mu\mathbf{1}_0P_\nu\mathbf{1}_0$ and we cannot
just concatenate walks as in the alcove walk algebra).  Define
$p_1\otimes p_2$, recursively, by
$$\begin{array}{lll}
\hbox{If $\varphi(p_1)=\iota(p_2)$}\quad  &\hbox{then}\quad
&p_1\otimes p_2 = p_1p_2,\ \hbox{the concatenation of $p_1$ and $p_2$, and} \\
\hbox{if $\varphi(p_1)\ne \iota(p_2)$}
&\hbox{then}\quad
&p_1\otimes p_2 = p_1'\otimes p_2,
\end{array}
$$
where $p_1'$ is the alcove walk constructed by the following procedure.  
Let $c_{i_1}^{\varepsilon_1}\cdots c_{i_r}^{\varepsilon_r}$
be a minimal length walk from $\varphi(p_1)$ to $\iota(p_2)$.
If $\varepsilon_1=-$ let $p_1' = p_1c_{i_1}^-$.
If $\varepsilon_1=+$ let $H_{\alpha,j}$ be the hyperplane crossed by
the last step of $p_1c_{i_1}^+$ and change the last negative crossing of $H_{\alpha,j}$ 
in $p_1$ to a fold to obtain a new path $p_1'$ with $\varphi(p_1')=\varphi(p_1c_{i_1}^+)$.

In terms of root operators, the Leibnitz rule (4.27) translates to the property
\begin{equation}
\begin{array}{rl}
\tilde e_i(p_1\otimes p_2) &= \begin{cases}
(\tilde e_ip_1)\otimes p_2, &\hbox{if $d_i^+(p_1)\ge d_i^-(p_2)$}, \\
p_1\otimes (\tilde e_i p_2), &\hbox{if $d_i^+(p_1)< d_i^-(p_2)$,}
\end{cases}
\quad\hbox{and} \\ \\
\tilde f_i(p_1\otimes p_2) &= \begin{cases}
(\tilde f_ip_1)\otimes p_2, &\hbox{if $d_i^+(p_1)>d_i^-(p_2)$}
\quad
\big(\overline{p_1\otimes p_2} = 
\beginpicture
\setcoordinatesystem units <1cm,1cm>         
\setplotarea x from -1.1 to 1, y from -0.5 to 0.5  
\put{$\bullet$} at 0 0 
\plot -1 0  0 0 /
\plot -1 0  -1 0.1 /
\plot -1 0.1  0.5 0.1 /
\plot 0.5 0.1  0.5 0.2 /
\plot 0.5 0.2  -0.5 0.2 /
\plot -0.5 0.2  -0.5 0.3 /
\arrow <5pt> [.2,.67] from -0.5 0.3 to 0.8 0.3 %
\endpicture
\big), \\ 
p\otimes (\tilde f_i q), &\hbox{if $d_i^+(p)\le d_i^-(q)$}
\quad
\big(
\overline{p\otimes q} = 
\beginpicture
\setcoordinatesystem units <1cm,1cm>         
\setplotarea x from -1.3 to 0.6, y from -0.5 to 0.5  
\put{$\bullet$} at 0 0 
\plot -0.4 0  0 0 /
\plot -0.4 0  -0.4 0.1 /
\plot -0.4 0.1  0.5 0.1 /
\plot 0.5 0.1  0.5 0.2 /
\plot 0.5 0.2  -1.2 0.2 /
\plot -1.2 0.2  -1.2 0.3 /
\arrow <5pt> [.2,.67] from -1.2 0.3 to -0.3 0.3 %
\endpicture
\big),
\end{cases}
\end{array}
\end{equation}
for $p_1,p_2\in B_{\mathrm{univ}}$.
It follows from this version of the Leibnitz rule that if $B_1$ and $B_2$
are $q$-crystals then the product
\begin{equation}
B_1\otimes B_2 = \{ p_1\otimes p_2\ |\ p_1\in B_1, p_2\in B_2\}
\qquad\hbox{is also a $q$-crystal,}
\end{equation}
and
\begin{equation}
\mathrm{char}(B_1\otimes B_2) = \mathrm{char}(B_1)\mathrm{char}(B_2).
\end{equation}
This last property is not completely trivial.  The general case follows from the
rank one case (projecting onto the line perpendicular to $H_{\alpha}$).
More importantly, the definition of the product $\otimes$ is forced by (4.28-4.30).

The following theorem is a version of [Sc, Theorem 1.2] and results in 
[KM] and [Ha].

\begin{thm}  Recall the notations from Theorem 4.2.
If $\lambda\in P^+$ let $B_q(p_\lambda)$ be the $q$-crystal generated by $p_\lambda$,
where $p_\lambda$ is a fixed minimal length alcove walk from $A$ to $\lambda+A$.
For $\mu, \nu\in P^+$,
$$
P_\mu P_\nu 
= \sum_{p\in B_q(p_\nu)\atop p_\mu\otimes p\subseteq C-\rho}
q^{-(\ell(\iota(p))+\ell(\varphi(p))-f(p))}(1-q^{-2})^{f(p)-c(p)}
P_{\mu+\wt(p)}. 
$$
\end{thm}
\begin{proof}  Using (4.29) and (4.30) and applying Theorem 4.8 to the
$q$-crystal $B_q(p_\mu)\otimes B_q(p_\nu)$ gives
\begin{align*}
P_\mu P_\nu\mathbf{1}_0 
&= \mathrm{char}(B_q(p_\mu))
\mathrm{char}(B_q(p_\nu))\mathbf{1}_0 
= \mathrm{char}((B_q(p_\mu)\otimes B_q(p_\nu))\mathbf{1}_0 \\
&= \sum_{p=p_1\otimes p_2\in B_q(p_\mu)\otimes B_q(p_\nu)\atop p_1\otimes p_2\subseteq C-\rho}
q^{-(\ell(\iota(p))+\ell(\varphi(p))-f(p))}(1-q^{-2})^{f(p)-c(p)}
P_{\mu+\wt(p_2)}\mathbf{1}_0 
\end{align*}
since every path in $p_1\in B_q(p_\mu)$ which is contained in $C-\rho$ has weight $\mu$
(so that $\wt(p_1\otimes p_2) = \mu+\wt(p_2)$.)
\end{proof}

Fix $J\subseteq \{1,2,\ldots, n\}$.  The subgroup of $W$ generated by the reflections in the
hyperplanes $H_{\alpha_j}$, $j\in J$,
$$W_J = \langle s_j\ |\ j\in J\rangle,
\quad\hbox{acts on $\fh_\RR^*$, \qquad with\qquad}
C_J = \{ \mu\in \fh_\RR^*\ |\ \hbox{$\langle \mu,\alpha_j^\vee\rangle >0$ for $j\in J$}\}
$$
as a fundamental chamber.
Let 
$H_J = \hbox{span}\{ T_{w^{-1}}^{-1}\ |\ w\in W_J\}$ and let $\mathbf{1}_J\in H_J$ be given by 
$$\mathbf{1}_J^2 = \mathbf{1}_J
\qquad\hbox{and}\qquad
T_{w^{-1}}^{-1}\mathbf{1}_J = q^{-\ell(w)}\mathbf{1}_J,
\quad\hbox{for $w\in W_J$}.
$$
For $\mu\in P$ let $W^J_\mu$ be the stabilizer of $\mu$ under the $W_J$ action on $P$ and
\begin{equation}
\hbox{define\quad $P^J_\mu(X;q^{-2})\in \KK[P]^{W_J}$\quad by\qquad}
P^J_\mu(X;q^{-2}) \mathbf{1}_0
=\Big(\sum_{w\in W^\mu_J} q^{-\ell(w)} T_{w^{-1}}^{-1}\Big)
X^\mu \mathbf{1}_0,
\end{equation}
where $W^\mu_J$ is a set of minimal length coset representatives for the
cosets in $W_J/W^J_\mu$.  
Then
\begin{equation}
\hbox{up to normalization}\quad
P^J_\lambda(X, q^{-2})\mathbf{1}_0
\quad\hbox{equals}\quad
\mathbf{1}_J X^\lambda\mathbf{1}_0,
\end{equation}
and 
\begin{equation}
\mathbf{1}_J\tilde H\mathbf{1}_0
\qquad\hbox{has basis}\qquad
\{ P_\lambda^J\mathbf{1}_0 \ |\ \lambda\in P_J^+\},
\qquad\hbox{where $P_J^+ = P\cap \overline{C_J}$.}
\end{equation}
with $\overline{C_J} 
= \{ \mu\in \fh_\RR^*\ |\ \hbox{$\langle \mu,\alpha_j^\vee\rangle >0$ for $j\in J$}\}$.

\begin{thm}  Recall the notations from Theorem 4.2.
If $\mu\in P^+$ let 
$\hbox{$B_q(p_\mu)$ be the $q$-crystal generated by $p_\mu$,}$
where $p_\mu$ is a fixed minimal length alcove walk from $A$ to $\mu+A$.
Let $\lambda\in P^+$ and
let  $J\subseteq \{1,2,\ldots, n\}$.  
Then
$$
P_\lambda = \sum_{p\in B(\lambda)\atop p\subseteq C_J-\rho_J} 
q^{-(\ell(\iota(p))+\ell(\varphi(p))-f(p))}(1-q^{-2})^{f(p)-c_J(p)}
P_{\wt(p)}^J,$$
where $C_J-\rho_J = \{ \mu\in \fh_\RR^*\ |\ \hbox{$\langle \mu,\alpha_j^\vee\rangle > -1$ for $j\in J$}\}$
and $c_J(p)$ is the number of folds of $p$ which touch a hyperplane 
$H_{\alpha_j}$ with $j\in J$.
\end{thm}
\begin{proof}  
A \emph{$J$-crystal} is a set of positively folded alcove walks $B$ 
which is closed under the operators $\tilde e_j$,
$\tilde f_j$, for $j\in J$.
Since $P_\lambda\mathbf{1}_0 
= \mathrm{char}(B_q(p_\lambda))\mathbf{1}_0$, 
the statement follows by applying Theorem 4.8 to $B_q(p_\lambda)$ viewed as a $J$-crystal.
\end{proof}

\end{subsection}

\end{section}

\begin{section}{Weyl characters and crystals}

Section 5.1 is an exposition of the theory of Weyl characters
analogous to the theory of Schur functions in [Mac2, Ch.\ 1].  
The element $a_\rho$ in Theorem 5.1 is the \emph{Weyl denominator},  
Lemma 5.2 is a generalization of the \emph{Jacobi-Trudi formula} and the
formulas in Proposition 5.3 are the \emph{quantum dimension formula} and the 
\emph{Weyl dimension formula}, respectively.  The results in Proposition 5.4b and 5.4c 
are the \emph{Kostant partition function formula} and the 
\emph{Brauer-Klimyk formula}, respectively.
Sections 5.2--5.6 give an elementary exposition of the 
theory of crystals and the path model and Section 5.7 explains the relationship
between crystals and column strict tableaux.   
The presentation here is designed to make clear the relationship between the general  
path model and the crystal operators of Lascoux and
Sch\"utzenberger used in the type A case [LS] (see [Ki] for a nice presentation).

The relationship between the path model used in this section and the alcove walks used
in Sections 3 and 4 is as follows.  Let $r\in \RR_{>0}$.  The \emph{dilation} which replaces the 
fundamental alcove $A$ by $\frac{1}{r}A$  induces important maps between the corresponding 
affine Hecke algebras,
$$\tilde H_A\hookrightarrow \tilde H_{\frac{1}{r}A}
\qquad\hbox{and}\qquad
\tilde H_{\frac{1}{r}A}\longrightarrow \tilde H_A,$$
corresponding to stretching the walks by a factor of $r$.
As $r$ gets large the alcove gets small and alcove walks become continuous paths in the limit.
At $q^{-1}=0$ the formula for the Hall-Littlewood polynomial in 
Theorem 4.2a becomes the path model for Weyl characters discovered by P. Littelmann [Li1-3].  
In other words, the $q$-crystals become ``classical'' crystals in the limit.  

\begin{subsection}{Schur functions}

Use notations for the Weyl group $W$ and the lattice $P$ as in Section 2.
The \emph{group algebra of $P$} is the ring
$$\ZZ[P]
\quad\hbox{with basis}\quad
\{ X^\lambda\ |\ \lambda\in P\}
\qquad\hbox{and product}\qquad
X^\lambda X^\mu = X^{\lambda+\mu},
$$
for $\lambda,\mu\in P$.  
The group $W$ acts on $\ZZ[P]$ by
$$wX^\lambda = X^{w\lambda},
\qquad\hbox{for $w\in W$, $\lambda\in P$.}
$$
The ring of \emph{symmetric functions} and \emph{Fock space} are  
\begin{equation}
\begin{array}{rl}
\ZZ[P]^W  &= \{ f\in \ZZ[P]\ |\ \hbox{$wf =f$ for all $w\in W$}\}\qquad\hbox{and} \\
\ZZ[P]^{\det} &= \{ f\in \ZZ[P]\ |\ \hbox{$wf = \det(w)f$ for all $w\in W$}\},
\end{array}
\end{equation}
respectively.
For $\lambda\in P$ define
\begin{equation}
m_\lambda = \sum_{\gamma\in W\lambda} X^{\gamma}
\qquad\hbox{and}\qquad
a_\lambda = \sum_{w\in W} \det(w) X^{w\lambda}.
\end{equation}
The straightening laws for these elements are
\begin{equation}
m_{w\lambda} = m_\lambda
\qquad\hbox{and}\qquad
a_{w\lambda} = \det(w)a_\lambda,
\qquad
\hbox{for $w\in W$ and $\lambda\in P$.}
\end{equation}
The second relation implies that $a_\lambda=0$ if there 
exists $w\in W_\lambda$ with $\det(w)\ne 1$, and it follows from the straightening
laws that
\begin{equation}
\begin{array}{lll}
\ZZ[P]^W\quad &\hbox{has basis}\quad &\{m_\lambda\ |\ \lambda\in P^+\},
\qquad\hbox{and} \\
\ZZ[P]^{\det} &\hbox{has basis} &\{ a_{\lambda+\rho}\ |\ \lambda\in P^+\},
\end{array}
\end{equation}
where $P^+$ and $\rho$ are as in (2.14) and (2.16), respectively.

The \emph{Weyl characters} or \emph{Schur functions} are defined by
\begin{equation}
s_\lambda = \frac{a_{\lambda+\rho}}{a_\rho},
\qquad\hbox{for $\lambda\in P$.}
\end{equation}
The following theorem shows that the $s_\lambda$ are elements of $\ZZ[P]$ and that 
\begin{equation}
\ZZ[P]^W\qquad\hbox{has basis}\qquad
\{ s_\lambda\ |\ \lambda\in P^+\}.
\end{equation}

\begin{thm}  Fock space $\ZZ[P]^{\det}$ is a free $\ZZ[P]^W$ module with generator 
$$
a_\rho
=x^\rho\prod_{\alpha\in R^+} (1-x^{-\alpha})
\qquad\hbox{and the map}\qquad
\begin{matrix}
\ZZ[P]^W &\longrightarrow &\ZZ[P]^{\det} \\
f &\longmapsto &a_\rho f \\
s_\lambda &\longmapsto &a_{\lambda+\rho}
\end{matrix}
$$
is a $\ZZ[P]^W$ module isomorphism.
\end{thm}
\begin{proof}
Let  $f\in \ZZ[P]^{\det}$ and let $\alpha\in R^+$.
If $f_\gamma$ is the coefficient of $x^\gamma$ in $f$ then
$$\sum_{\gamma\in P} f_\gamma x^\gamma = f 
= -s_\alpha f = \sum_{\gamma\in P} -f_\gamma x^{s_\alpha\gamma},
\qquad\hbox{and so}\qquad
f = \sum_{\gamma\in P\atop \langle \gamma,\alpha^\vee\rangle\ge 0} f_\gamma(x^\gamma-x^{s_\alpha\gamma}),
$$
since $f_{s_\alpha\gamma} = -f_\gamma$.
Since  each term $x^\gamma-x^{s_\alpha\gamma}$ is divisible $1-x^{-\alpha}$,
$f$ is divisible by $1-x^{-\alpha}$,
and thus
\begin{equation}
\hbox{each}\quad f\in \ZZ[P]^{\det} 
\quad\hbox{is divisible by}\quad
x^\rho\prod_{\alpha\in R^+} (1-x^{-\alpha})
\end{equation}
since the polynomials $1-x^{-\alpha}$, $\alpha\in R^+$
are coprime in $\ZZ[P]$ and $x^\rho$ is a unit in $\ZZ[P]$.
Comparing coefficients of the maximal terms in $a_\rho$
and $x^\rho\prod_{\alpha\in R^+} (1-x^{-\alpha})$  shows that 
$$a_\rho =x^\rho\prod_{\alpha\in R^+} (1-x^{-\alpha}).$$
Thus each $f\in \ZZ[P]^{\det}$ is divisible by $a_\rho$
and so the inverse of multiplication by $a_\rho$ 
is well defined. 
\end{proof}

The {\it dot action} of $S_n$ on $P$ is given by
\begin{equation}
w\circ \mu = w(\mu+\rho)-\rho,
\qquad\hbox{for $w\in S_n$, $\mu\in P$.}
\end{equation}
The straightening law 
\begin{equation}
s_{w\circ\mu} = \det(w) s_\mu,
\qquad\hbox{for $\mu\in P$, $w\in W$.} 
\end{equation}
for the Schur functions follows from the straightening law for the 
$a_\mu$ in (5.3).

\begin{lemma}   Let  $f\in \ZZ[P]^W$ and write
$\displaystyle{ f = \sum_{\gamma} f_\gamma x^\gamma }$ so that
$f_\gamma$ is the coefficient of $x^\gamma$ in $f$.  Then
$$f = \sum_{\mu\in P^+} f_\mu m_\mu = \sum_{\lambda\in P^+} \eta^\lambda s_\lambda,
\qquad\hbox{where}\qquad
\eta^\lambda = \sum_{w\in W} \det(w) f_{\lambda+\rho-w\rho}.$$
\end{lemma}
\begin{proof}
The first equality is immediate from the definition of $m_\mu$.
Since $f\in \ZZ[P]^W$ and the $s_\lambda$, $\lambda\in P^+$, are a basis of $\ZZ[P]^W$,
the element $f$ can be written as a linear combination of $s_\lambda$.
Then, since $e^{\lambda+rho}$ is the unique dominant term in $a_{\lambda+\rho}$,
\begin{align*}
\eta_\lambda 
&= \hbox{(coefficient of $s_\lambda$ in $f$)} 
= \hbox{(coefficient of $a_{\lambda+\rho}$ in $fa_\rho$)} \\
&= \left(\hbox{coefficient of $e^{\lambda+\rho}$ in}
~~ \sum_{\mu\in P} \sum_{w\in W} \det(w) f_\mu e^{\mu+w\rho}
\right).
\end{align*}
\end{proof}

\begin{prop}
If $\nu\in \fh^*_\RR$ and 
$\displaystyle{
f=\sum_{\mu\in P} f_\mu e^\mu\in \ZZ[P] }$ define
$\displaystyle{
f(e^\nu) = \sum_{\mu\in P} f_\mu e^{\langle \mu,\nu\rangle}. }$
\smallskip\noindent
Let $\lambda\in P^+$, $t\in \RR_{>0}$, $q = e^t$ and $\rho^\vee = \frac12\sum_{\alpha\in R^+} \alpha^\vee$.  Then
$$s_\lambda(q^{\rho^\vee}) 
= \prod_{\alpha\in R^+} {[\langle \lambda+\rho,\alpha^\vee\rangle]
\over [\langle \rho,\alpha^\vee\rangle] }
\qquad\hbox{and}\qquad
s_\lambda(1) 
= \prod_{\alpha\in R^+} {\langle \lambda+\rho,\alpha^\vee\rangle
\over \langle \rho,\alpha^\vee\rangle }
$$
where $[k] = (q^k-1)/(q-1)$ for an integer $k\ne 0$.
\end{prop}
\begin{proof}
\begin{align*}
a_{\lambda+\rho}(e^{t{\rho^\vee}}) 
&=
\sum_{w\in W} \det(w) e^{\langle w(\lambda+\rho),t{\rho^\vee}\rangle} 
=
\sum_{w\in W} \det(w) e^{\langle w\rho^\vee,t(\lambda+\rho)\rangle} \\
&=
a_{\rho^\vee}(e^{t(\lambda+\rho)}) 
= e^{\langle \rho^\vee,t(\lambda+\rho)\rangle} 
\prod_{\alpha\in R^+} 
(1-e^{\langle -\alpha^\vee,t(\lambda+\rho)\rangle} ).
\end{align*}
Thus
$$s_\lambda(e^{t{\rho^\vee}})
=\frac{ a_{\lambda+\rho}(e^{t{\rho^\vee}}) }{a_\rho(e^{t{\rho^\vee}}) }
= \frac{e^{\langle \rho^\vee,t(\lambda+\rho)\rangle} }{ e^{\langle \rho^\vee, t\rho\rangle} }
\prod_{\alpha\in R^+} \frac{ 
1-e^{\langle -\alpha^\vee,t(\lambda+\rho)\rangle}
}{ 1-e^{\langle -\alpha^\vee,t\rho\rangle}
}  
= q^{-\langle \lambda,\rho^\vee\rangle}
\prod_{\alpha\in R^+} \frac{ q^{\langle \lambda+\rho,\alpha^\vee\rangle}-1 }
{ q^{\langle \rho,\alpha^\vee\rangle}-1 }
$$
and
$$
s_\lambda(1) = \lim_{q\to 1} s_\lambda(q^{\rho^\vee})
= \prod_{\alpha\in R^+} \frac{ \langle\lambda+\rho,\alpha^\vee\rangle }
{ \langle \rho,\alpha^\vee\rangle } .
$$
\end{proof}

The \emph{weight multiplicities} are the integers
$K_{\lambda\gamma}$, $\lambda\in P^+$, $\gamma\in P$,
defined by the equations
\begin{equation}
s_\lambda
= \sum_{\gamma\in P} K_{\lambda\gamma} x^\gamma
= \sum_{\mu\in P^+} K_{\lambda\mu} m_\mu.
\end{equation}
The \emph{tensor product multiplicities}  are
the integers $c_{\mu\nu}^\lambda$, $\mu,\nu,\lambda\in P^+$, defined by the 
equations
\begin{equation}
s_\mu s_\nu = \sum_{\lambda\in P^+} c_{\mu\nu}^\lambda s_\lambda.
\end{equation}
The {\it partition function} is the function $p\colon P\to \ZZ_{\ge0}$ defined
by the equation
\begin{equation}
\prod_{\alpha\in R^+} {1\over 1-x^{-\alpha}} 
= \sum_{\gamma\in P} p(\gamma)x^{-\gamma}.
\end{equation}

\begin{prop}   Let $\lambda,\mu,\nu\in P^+$. 
\item[(a)] $K_{\lambda\lambda}=1$,\quad $K_{\lambda,w\mu}=K_{\lambda\mu}$,
for $w\in W$, \quad and \quad
$K_{\lambda\mu}=0$ unless $\mu\le\lambda$.
\smallskip\noindent
\item[(b)] $\displaystyle{
K_{\lambda\mu} = \sum_{w\in W} \det(w)p(w(\lambda+\rho)-(\mu+\rho)).}$
\smallskip\noindent
\item[(c)]
$\displaystyle{
c_{\mu\nu}^\lambda = \sum_{v,w\in W} \det(vw) p (v(\mu+\rho)+w(\nu+\rho)-(\lambda+\rho)-\rho).}$
\end{prop}
\begin{proof}
(a) The equality $K_{\lambda,w\mu} = K_{\lambda\mu}$ follows from the definition and the fact that
$s_\lambda\in \ZZ[P]^W$.
If $w\ne 1$ then $w(\lambda+\rho)<\lambda+\rho$ so that $w(\lambda+\rho)-\rho < \lambda$
and
$$s_\lambda = \left(\sum_{w\in W} \det(w)x^{w(\lambda+\rho)-\rho}\right)\cdot
\prod_{\alpha\in R^+} {1\over 1-x^{-\alpha}}
= x^\lambda+(\hbox{lower terms in dominance order}).$$
Thus $K_{\lambda\lambda}=1$ and $K_{\lambda\mu}=0$ unless $\mu\le \lambda$.

\smallskip\noindent
(b) The coefficient of $x^\mu$ in 
$$s_\lambda = \left(\sum_{w\in W} \det(w)x^{w(\lambda+\rho)-\rho}\right)
\prod_{\alpha\in R^+} {1\over 1-x^{-\alpha}}
=\sum_{w\in W\atop \gamma\in Q^+} \det(w)p(\gamma)x^{w(\lambda+\rho)-\rho-\gamma},$$
has a contribution $\det(w)p(\gamma)$ when $w(\lambda+\rho)-\rho-\gamma=\mu$ so that 
$\gamma=w(\lambda+\rho)-(\mu+\rho)$.  

\smallskip\noindent
(c)  Let $\varepsilon =\sum_{w\in W} \det(w)w$.  Since $c_{\mu\nu}^\lambda$ is the coefficient of $x^{\nu+\rho}$ in 
\begin{align*}
s_\mu s_\nu a_\rho 
&= {\varepsilon(x^{\mu+\rho})\varepsilon(x^{\nu+\rho})\over a_\rho}
=\left(\sum_{v,w\in W} \det(vw)x^{v(\mu+\rho)+w(\nu+\rho)-\rho}\right)
\left(\prod_{\alpha\in R^+} {1\over 1-x^{-\alpha}}\right)\cr
&=\sum_{v,w\in W\atop \gamma\in Q^+} \det(vw) p(\gamma)x^{v(\mu+\rho)+w(\nu+\rho)-\gamma-\rho},
\end{align*}
there is a contribution $\det(vw)p(\gamma)$ to the coefficient $c_{\mu\nu}^\lambda$ when 
$\lambda+\rho = v(\mu+\rho)+w(\nu+\rho)-\gamma-\rho$ so that 
$\gamma = v(\mu+\rho)+w(\mu+\rho) - (\lambda+\rho)-\rho.$
\end{proof}

Fix $J\subseteq \{1,2,\ldots, n\}$.  The subgroup of $W$ generated by the reflections in the
hyperplanes $H_{\alpha_j}$, $j\in J$,
\begin{equation}
W_J = \langle s_j\ |\ j\in J\rangle, \quad\hbox{acts on $\fh_\RR^*$,}
\qquad\hbox{with}\qquad
C_J = \{ \mu\in \fh_\RR^*\ |\ \hbox{$\langle \mu,\alpha_j^\vee\rangle > 0$ for $j\in J$}\}
\end{equation}
as a fundamental chamber.
The group $W_J$ acts on $P$ and 
\begin{equation}
\ZZ[P]^{W_J} = \{ f\in \ZZ[P]\ |\ \hbox{$wf=f$ for $w\in W_J$}\}
\end{equation}
is a subalgebra of $\ZZ[P]$ which contains $\ZZ[P]^W$.
If $\overline{C_J}
=\{ \mu\in \fh_\RR^*\ |\ \hbox{$\langle \mu,\alpha_j^\vee\rangle \ge 0$ for $j\in J$}\}$,
\begin{equation}
P_J^+ = P\cap \overline{C_J}, \qquad
\qquad \rho_J = \sum_{j\in J} \omega_j,
\end{equation}
\begin{equation}
a_\mu^J = \sum_{w\in W_J} \det(w) w X^\mu, \quad \hbox{for $\mu\in P$,} 
\qquad \hbox{and} \qquad
s_\lambda^J = \frac{a_{\lambda+\rho_J}^J}{a_{\rho_J}^J}, \quad \hbox{for $\lambda\in P$,}
\end{equation}
then
$$\hbox{$\{ s_\lambda^J\ |\ \lambda\in P_J^+\}$ is a basis of $\ZZ[P]^{W_J}$.}$$
The \emph{restriction multiplicities} are the integers $c^\lambda_{J,\nu}$ given by
\begin{equation}
s_\lambda = \sum_{\nu\in P_J^+} c^\lambda_{J,\nu} s_\nu^J.
\end{equation}

\end{subsection}

\begin{subsection}{Paths}

Let $\lambda\in P$.  The \emph{straight line path} to $\lambda$ is the map
\begin{equation}
p_\lambda\colon [0,1] \to\fh_\RR^*
\quad\hbox{given by}\quad p_\lambda(t)=\lambda t.
\end{equation}
Let $\ell_1,\ell_2\in \RR_{\ge 0}$.  The \emph{concatenation}
of maps $p_1\colon [0,\ell_1]\to \fh_\RR^*$ and $p_2\colon [0,\ell_2]\to \fh_\RR^*$ is 
the map
$p_1\otimes p_2\colon [0,\ell_1+\ell_2]\to \fh_\RR^*$ given by
\begin{equation}
(p_1\otimes p_2)(t) = \begin{cases}
p_1(t), &\hbox{for $t\in [0,\ell_1]$}, \\
p_1(\ell_1)+p_2(t-\ell_1), &\hbox{for $t\in [\ell_1,\ell_1+\ell_2]$}.
\end{cases}
\end{equation}
Let $r,\ell\in \RR_{\ge 0}$.  The \emph{$r$-stretch} of a map
$p\colon [0,\ell]\to \fh_\RR^*$ is the map $rp\colon [0,r\ell]\to \fh_\RR^*$ given by
\begin{equation}
(rp)(t) = r\cdot p(t/r).
\end{equation}
The \emph{reverse} of a map $p\colon [0,\ell]\to \fh_\RR^*$ is the map
$p^*\colon [0,\ell]\to \fh_\RR^*$ given by
\begin{equation}
p^*(t) = p(\ell-t)-p(\ell).
\end{equation}
The \emph{weight} of a map $p\colon [0,\ell]\to\fh_\RR^*$ is the endpoint of $p$,
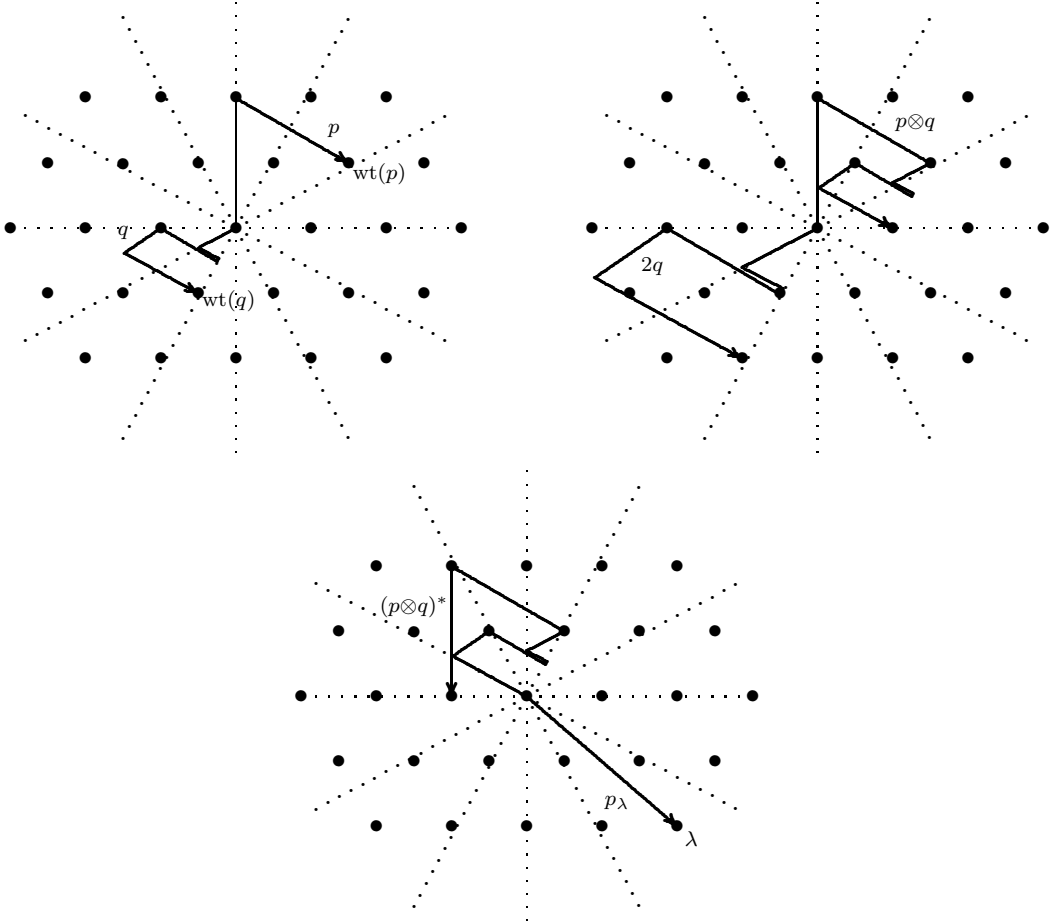
\begin{equation}
\wt(p) = p(\ell).
\end{equation}
$$
\beginpicture
\setcoordinatesystem units <1cm,1cm>         
\setplotarea x from -3 to 3, y from -3 to 3    
\put{$\scriptstyle{p}$} at 1.3 1.3
\put{$\scriptstyle{q}$} at -1.5 -0.05
\put{$\scriptstyle{\wt(p)}$}[tl] at 1.55 0.873
\put{$\scriptstyle{\wt(q)}$}[tl] at -0.45  -0.816
\put{$\bullet$} at 0 0
\put{$\bullet$} at 1 0
\put{$\bullet$} at 2 0
\put{$\bullet$} at 3 0
\put{$\bullet$} at -1 0
\put{$\bullet$} at -2 0
\put{$\bullet$} at -3 0
\put{$\bullet$} at 0.5  0.8660
\put{$\bullet$} at 1.5  0.8660
\put{$\bullet$} at 2.5  0.8660
\put{$\bullet$} at -0.5  0.8660
\put{$\bullet$} at -1.5  0.8460
\put{$\bullet$} at -2.5  0.8660
\put{$\bullet$} at 0  1.732
\put{$\bullet$} at 1  1.732
\put{$\bullet$} at 2  1.732
\put{$\bullet$} at -1  1.732
\put{$\bullet$} at -2  1.732
\put{$\bullet$} at 0.5  -0.8660
\put{$\bullet$} at 1.5  -0.8660
\put{$\bullet$} at 2.5  -0.8660
\put{$\bullet$} at -0.5  -0.8660
\put{$\bullet$} at -1.5  -0.8660
\put{$\bullet$} at -2.5  -0.8660
\put{$\bullet$} at 0  -1.732
\put{$\bullet$} at 1  -1.732
\put{$\bullet$} at 2  -1.732
\put{$\bullet$} at -1  -1.732
\put{$\bullet$} at -2  -1.732
\linethickness=0.5pt                          
\putrule from 0 1.732 to 0 0          
\plot 0 1.732 1.5  0.8660 /    
\plot 0 0 -0.50 -0.26  /  
\plot -0.50 -0.26 -0.23 -0.40 /  %
\plot -0.23 -0.40 -0.25 -0.44 /  %
\plot -0.25 -0.44 -1 0 /  %
\plot -1 0. -1.48  -0.3287 /    
\plot  -1.48  -0.3287  -0.5 -0.8660 /    %
\arrow <5pt> [.2,.67] from 1.35 0.96  to 1.5 0.8660   %
\arrow <5pt> [.2,.67] from -0.6 -0.81  to -0.5 -0.8660   %
\setdots
\putrule from 0 3 to 0 -3        
\plot -1.5 -2.798   1.5 2.798 /  
\putrule from 3 0  to -3 0       
\plot -2.798 -1.5  2.798 1.5 /   
\plot -2.798 1.5  2.798 -1.5 /   
\plot -1.5 2.798  1.5 -2.798 /   
\endpicture
\qquad\qquad
\beginpicture
\setcoordinatesystem units <1cm,1cm>         
\setplotarea x from -3 to 3, y from -3 to 3    
\put{$\scriptstyle{p\otimes q}$} at 1.3 1.4
\put{$\scriptstyle{2q}$} at -2.2 -0.5
\put{$\bullet$} at 0 0
\put{$\bullet$} at 1 0
\put{$\bullet$} at 2 0
\put{$\bullet$} at 3 0
\put{$\bullet$} at -1 0
\put{$\bullet$} at -2 0
\put{$\bullet$} at -3 0
\put{$\bullet$} at 0.5  0.8660
\put{$\bullet$} at 1.5  0.8660
\put{$\bullet$} at 2.5  0.8660
\put{$\bullet$} at -0.5  0.8660
\put{$\bullet$} at -1.5  0.8460
\put{$\bullet$} at -2.5  0.8660
\put{$\bullet$} at 0  1.732
\put{$\bullet$} at 1  1.732
\put{$\bullet$} at 2  1.732
\put{$\bullet$} at -1  1.732
\put{$\bullet$} at -2  1.732
\put{$\bullet$} at 0.5  -0.8660
\put{$\bullet$} at 1.5  -0.8660
\put{$\bullet$} at 2.5  -0.8660
\put{$\bullet$} at -0.5  -0.8660
\put{$\bullet$} at -1.5  -0.8660
\put{$\bullet$} at -2.5  -0.8660
\put{$\bullet$} at 0  -1.732
\put{$\bullet$} at 1  -1.732
\put{$\bullet$} at 2  -1.732
\put{$\bullet$} at -1  -1.732
\put{$\bullet$} at -2  -1.732
\linethickness=0.5pt                          
\putrule from 0 1.732 to 0 0          
\plot 0 1.732 1.5  0.8660 /    
\plot  1.5 0.866   1 0.606  /  
\plot  1 0.606    1.27  0.466 /  %
\plot  1.27 0.466  1.25 0.426 /  %
\plot  1.25 0.426  0.5 0.866 /  %
\plot  0.5 0.866   0.02  0.5373 /    
\plot  0.02  0.5373   1 0 /    %
\plot  0 0  -1 -0.52  /  
\plot -1 -0.52  -0.46 -0.8 /  %
\plot -0.46 -0.8  -0.5 -0.88 /  %
\plot -0.5 -0.88  -2 0 /  %
\plot -2 0   -2.96  -0.6574 /    
\plot  -2.96  -0.6574   -1 -1.732 /    %
\arrow <5pt> [.2,.67] from 0.87 0.08  to 1 0   %
\arrow <5pt> [.2,.67] from -1.13 -1.66  to -1 -1.732   %
\setdots
\putrule from 0 3 to 0 -3        
\plot -1.5 -2.798   1.5 2.798 /  
\putrule from 3 0  to -3 0       
\plot -2.798 -1.5  2.798 1.5 /   
\plot -2.798 1.5  2.798 -1.5 /   
\plot -1.5 2.798  1.5 -2.798 /   
\endpicture
$$
$$
\beginpicture
\setcoordinatesystem units <1cm,1cm>         
\setplotarea x from -3 to 3, y from -3 to 3    
\put{$\scriptstyle{\lambda}$}[tl] at 2.1 -1.8
\put{$\scriptstyle{(p\otimes q)^*}$} at -1.5 1.2
\put{$\scriptstyle{p_\lambda}$}[b] at 1.2 -1.5
\put{$\bullet$} at 0 0
\put{$\bullet$} at 1 0
\put{$\bullet$} at 2 0
\put{$\bullet$} at 3 0
\put{$\bullet$} at -1 0
\put{$\bullet$} at -2 0
\put{$\bullet$} at -3 0
\put{$\bullet$} at 0.5  0.8660
\put{$\bullet$} at 1.5  0.8660
\put{$\bullet$} at 2.5  0.8660
\put{$\bullet$} at -0.5  0.8660
\put{$\bullet$} at -1.5  0.8460
\put{$\bullet$} at -2.5  0.8660
\put{$\bullet$} at 0  1.732
\put{$\bullet$} at 1  1.732
\put{$\bullet$} at 2  1.732
\put{$\bullet$} at -1  1.732
\put{$\bullet$} at -2  1.732
\put{$\bullet$} at 0.5  -0.8660
\put{$\bullet$} at 1.5  -0.8660
\put{$\bullet$} at 2.5  -0.8660
\put{$\bullet$} at -0.5  -0.8660
\put{$\bullet$} at -1.5  -0.8660
\put{$\bullet$} at -2.5  -0.8660
\put{$\bullet$} at 0  -1.732
\put{$\bullet$} at 1  -1.732
\put{$\bullet$} at 2  -1.732
\put{$\bullet$} at -1  -1.732
\put{$\bullet$} at -2  -1.732
\linethickness=0.5pt                          
\putrule from -1 1.732 to -1 0          
\plot -1 1.732  0.5  0.8660 /    
\plot  0.5 0.866   0  0.606  /  
\plot  0 0.606    0.27  0.466 /  %
\plot  0.27 0.466  0.25 0.426 /  %
\plot  0.25 0.426  -0.5 0.866 /  %
\plot  -0.5 0.866   -0.98  0.5373 /    
\plot  -0.98  0.5373   0 0 /    %
\plot 0 0  2 -1.732 /  
\arrow <5pt> [.2,.67] from -1  .1  to -1 0   %
\arrow <5pt> [.2,.67] from 1.9  -1.632  to 2 -1.732   %
\setdots
\putrule from 0 3 to 0 -3        
\plot -1.5 -2.798   1.5 2.798 /  
\putrule from 3 0  to -3 0       
\plot -2.798 -1.5  2.798 1.5 /   
\plot -2.798 1.5  2.798 -1.5 /   
\plot -1.5 2.798  1.5 -2.798 /   
\endpicture
$$

Let 
\smallskip
\begin{equation}
\begin{array}{rl}
B_{\mathrm{univ}}\ \ &\hbox{be the set of maps generated by the straight line paths} \\
&\hbox{by operations of concatenation, reversing and stretching.}
\end{array}
\end{equation}
A \emph{path} is an element $p\colon [0,\ell]\to \fh_\RR^*$ in $B_{\mathrm{univ}}$.
Let $B$ be a set of paths (a subset of $B_{\mathrm{univ}}$).
The \emph{character} of $B$ is the element of $\ZZ[P]$ given by
\begin{equation}
\mathrm{char}(B) = \sum_{p\in B} X^{\wt(p)}.
\end{equation}

A \emph{crystal} is a set of paths $B$ that is closed under the action of
the \emph{root operators}
\begin{equation}
\tilde e_i\colon B_{\mathrm{univ}} \longrightarrow B_{\mathrm{univ}}\cup\{0\} 
\qquad\hbox{and}\qquad
\tilde f_i\colon B_{\mathrm{univ}} \longrightarrow B_{\mathrm{univ}} \cup\{0\},
\qquad 1\le i\le n,
\end{equation}
which are defined and constructed below, in Proposition 5.7 and Theorem 5.8.  
The \emph{crystal graph} of $B$ is the graph with 
\begin{equation}
\hbox{vertices $B$}\qquad\hbox{and}\qquad
\hbox{labeled edges\quad $p'\mapleft{i} p$\ \ \  if $p'=\tilde f_i p$.}
\end{equation}

\end{subsection}

\begin{subsection}{$i$-strings}

Let $B$ be a crystal.  Let $p\in B$ and fix $i$ ($1\le i\le n$).  The \emph{$i$-string}
of $p$ is the set of paths $S_i(p)$ generated from $p$ by applications of the
operators $\tilde e_i$ and $\tilde f_i$.  
\begin{enumerate}
\item[]  The \emph{head} of $S_i(p)$ is $h\in S_i(p)$
such that $\tilde e_i h = 0$.  
\item[]  The \emph{tail} of $S_i(p)$ is $t\in S_i(p)$ such that $\tilde f_i t = 0$.
\end{enumerate}
The weights of the paths in $S_i(p)$ are 
$$\wt(t)=s_i\wt(h) 
=\wt(h)-\langle \wt(h),\alpha_i^\vee\rangle\alpha_i,\ \ 
\ldots,\ \  \wt(h)-2\alpha_i,\ \  \wt(h)-\alpha_i,\ \  \wt(h),$$
and the crystal graph of $S_i(p)$ is 
$$\begin{array}{lcccr}
\scriptstyle{\vert}\longleftarrow &\scriptstyle{d_i^-(p)} 
&\longrightarrow\scriptstyle{\vert}\qquad \\ 
t\  \mapleft{i}\  \tilde e_i t\  \mapleft{i}\  &\cdots\  &\mapleft{i}\  \tilde f_ip\  \mapleft{i}\  
p\ \mapleft{i}\  \tilde e_ip\  \mapleft{i}\  &\cdots &\  \mapleft{i}\  \tilde f_i h\  \mapleft{i}\  h \\ 
&&\qquad \scriptstyle{\vert}\longleftarrow 
&\scriptstyle{d_i^+(p)} &\longrightarrow\scriptstyle{\vert}
\end{array}
$$
where
$$
d_i^+(p) = \hbox{(distance from $h$ to $p$)}
\qquad\hbox{and}\qquad
d_i^-(p) = \hbox{(distance from $p$ to $t$)},
$$
so that $\tilde e_i^{d_i^+(p)}p = h$ and $\tilde f_i^{d_i^-(p)}p=t,$

\end{subsection}

\begin{subsection}{Highest weight paths}

A \emph{highest weight path} is a path $p$ such that
\begin{equation}
\tilde e_i p =0,
\qquad\hbox{for all $1\le i\le n$}.
\end{equation}
A highest weight path is a path $p$ such that, for each $1\le i\le n$, $p$ is the head
of the $i$-string $S_i(p)$.  Thus $\langle p(t),\alpha_i^\vee\rangle> -1$ for all $t$
and all $1\le i\le n$.  So a path $p$ is a highest weight
path if and only if
\begin{equation}
p\subseteq C-\rho,
\qquad\hbox{where}\quad C-\rho = \{\mu-\rho\ |\ \mu\in C\}.
\end{equation}
Following the example at the end of Section 2, for the root system of type $C_2$ the picture is
$$
\begin{array}{c}
\beginpicture
\setcoordinatesystem units <1cm,1cm>         
\setplotarea x from -2.5 to 2.5, y from -2.5 to 2.5  
\put{$H_{\alpha_1}$}[b] at 0 2.1
\put{$H_{\alpha_2}$}[l] at 2.1 2.1
\put{$H_{\alpha_1+\alpha_2}$}[r] at -2.1 2.1
\put{$H_{\alpha_1+2\alpha_2}$}[l] at 2.1 0
\put{$C-\rho$} at 1 2
\put{$\scriptstyle{-\rho}$}[tr] at -0.5 -1.1
\put{$\bullet$} at  -0.5 -1
\put{$\scriptstyle{0}$}[bl] at 0.1 0.1
\plot -0.5 -1   2 1.5 /
\plot  -0.5 2  -0.5 -1 /
\vshade -0.5 -1 2   2 1.5 1.5 /
\setdashes
\plot -2 -2   2 2 /
\plot  2 -2  -2 2 /
\plot  0  2   0 -2 /
\plot  2  0  -2  0 /
\endpicture
\\
\hbox{the region $C-\rho$}
\end{array}
$$

If $p$ is a highest weight path with $\wt(p)\in P$ then, necessarily, $\wt(p)\in P^+$.
The following theorem gives an expression for the character of a crystal in terms of the
basis $\{s_\lambda\ |\ \lambda\in P^+\}$ of $\ZZ[P]^W$.

\begin{thm}  Let $B$ be a crystal.  Let $mathrm{char}(B)$ be as in (5.24) and
$s_\lambda$ as in (5.5).  Then
$$\mathrm{char}(B) = \sum_{p\in B\atop p\subseteq C-\rho} s_{\wt(p)},$$
where the sum is over highest weight paths $p\in B$.
\end{thm}
\begin{proof}
Fix $i$, $1\le i\le n$.  If $p\in B$ let $s_ip$ be the element of the $i$-string of $p$ which
satisfies
$$\wt(s_ip)= s_i\wt(p).
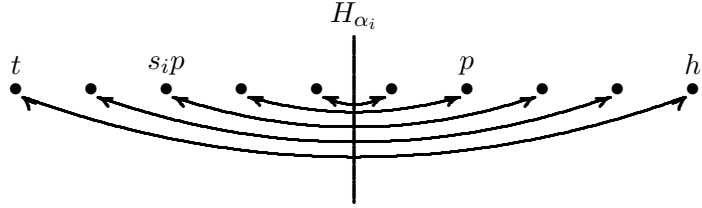
$$
$$\beginpicture
\setcoordinatesystem units <1cm,1cm>         
\setplotarea x from -2.5 to 2.5, y from -2.5 to 1.2  
\put{$H_{\alpha_i}$}[b] at 0 0.8
\put{$t$}[b] at -4.5 0.2
\put{$s_ip$}[b] at -2.5 0.2
\put{$p$}[b] at 1.5 0.2
\put{$h$}[b] at 4.5 0.2
\put{$\bullet$} at  -4.5 0
\put{$\bullet$} at  -3.5 0
\put{$\bullet$} at  -2.5 0
\put{$\bullet$} at  -1.5 0
\put{$\bullet$} at  -0.5 0
\put{$\bullet$} at  0.5 0
\put{$\bullet$} at  1.5 0
\put{$\bullet$} at  2.5 0
\put{$\bullet$} at  3.5 0
\put{$\bullet$} at  4.5 0
\plot 0 -1.5  0 0.7 /
\setquadratic \plot 4.4 -0.1                      
                      0 -0.9                     %
                     -4.4 -0.1  /                   %
\arrow <5pt> [.2,.67] from 4.35 -0.15  to 4.4 -0.1   %
\arrow <5pt> [.2,.67] from -4.35 -0.15  to -4.4 -0.1   %
\setquadratic \plot 3.4 -0.1                      
                      0 -0.7                     %
                     -3.4 -0.1  /                   %
\arrow <5pt> [.2,.67] from 3.3 -0.15  to 3.4 -0.1   %
\arrow <5pt> [.2,.67] from -3.3 -0.15  to -3.4 -0.1   %
\setquadratic \plot 2.4 -0.1                      
                     0 -0.5                     %
                     -2.4 -0.1  /                   %
\arrow <5pt> [.2,.67] from 2.3 -0.15  to 2.4 -0.1   %
\arrow <5pt> [.2,.67] from -2.3 -0.15  to -2.4 -0.1   %
\setquadratic \plot 1.4 -0.1                      
                    0 -0.3                     %
                     -1.4 -0.1  /                   %
\arrow <5pt> [.2,.67] from 1.25 -0.15  to 1.4 -0.1   %
\arrow <5pt> [.2,.67] from -1.25 -0.15  to -1.4 -0.1   %
\setquadratic \plot 0.4 -0.1
                     0 -0.2                      
                     -0.4 -0.1         /                   %
\arrow <5pt> [.2,.67] from 0.25 -0.15  to 0.4 -0.1   %
\arrow <5pt> [.2,.67] from -0.25 -0.15  to -0.4 -0.1   %
\endpicture
$$
Then $s_i(s_ip)=p$ and 
$$s_i\mathrm{char}(B) = \sum_{p\in B} X^{s_i\wt(p)} = \sum_{p\in B} X^{\wt(s_ip)}
=\mathrm{char}(B).$$
Hence
$\mathrm{char}(B)\in \ZZ[P]^W.$

Let 
$$\varepsilon = \sum_{w\in W} \det(w) w
\quad\hbox{so that}\quad
a_\mu = \varepsilon(X^\mu), \quad\hbox{for $\mu\in P$.}
$$
Since $\mathrm{char}(B)\in \ZZ[P]^W$,
$$\mathrm{char}(B)a_\rho = \mathrm{char}(B)\varepsilon(X^\rho)
=\varepsilon(\mathrm{char}(B)X^\rho)$$
and
\begin{equation}
\begin{array}{rcl}
\mathrm{char}(B) 
&=& \displaystyle{
\frac{1}{a_\rho}\mathrm{char}(B)a_\rho 
= \frac{\varepsilon(\mathrm{char}(B)X^\rho)}{a_\rho}  }  \\ \\
&=& \displaystyle{
\sum_{p\in B} \frac{\varepsilon(X^{\wt(p)+\rho})}{a_\rho}
=\sum_{p\in B} \frac{a_{\wt(p)+\rho}}{a_\rho}  
=\sum_{p\in B} s_{\wt(p)}. }
\end{array}
\end{equation}
There is some cancellation which can occur in this sum.  Assume $p\in B$ such that 
$p\not\subseteq C-\rho$ let $t$ be the first time that $p$ leaves the cone $C-\rho$.
In other words let $t\in \RR_{>0}$ be minimal such that there exists an $i$ with
$$p(t)\in H_{\alpha_i,-1}\qquad 
\hbox{where}\qquad
H_{\alpha_i,-1} = \{ \lambda\in \fh_\RR^*\ |\ \langle \lambda,\alpha_i^\vee \rangle = -1\}.
$$
Let $i$ be the minimal index such that the point $p(t)\in H_{\alpha_i,-1}$
and define $s_i\circ p$ to be the element of the $i$-string of $p$ such that 
$$\wt(s_i\circ p) = s_i\circ p.
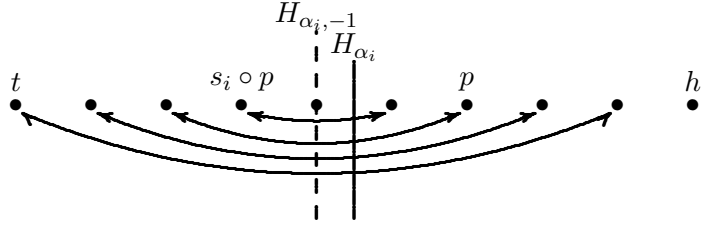
$$
$$\beginpicture
\setcoordinatesystem units <1cm,1cm>         
\setplotarea x from -2.5 to 2.5, y from -2.5 to 1.5  
\put{$H_{\alpha_i}$}[b] at 0 0.6
\put{$H_{\alpha_i,-1}$}[b] at -0.5 1
\put{$t$}[b] at -4.5 0.2
\put{$s_i\circ p$}[b] at -1.5 0.2
\put{$p$}[b] at 1.5 0.2
\put{$h$}[b] at 4.5 0.2
\put{$\bullet$} at  -4.5 0
\put{$\bullet$} at  -3.5 0
\put{$\bullet$} at  -2.5 0
\put{$\bullet$} at  -1.5 0
\put{$\bullet$} at  -0.5 0
\put{$\bullet$} at  0.5 0
\put{$\bullet$} at  1.5 0
\put{$\bullet$} at  2.5 0
\put{$\bullet$} at  3.5 0
\put{$\bullet$} at  4.5 0
\plot 0 -1.5  0 0.6 /
\setdashes
\plot -0.5 -1.5 -0.5 1 /
\setsolid
\setquadratic \plot 3.4 -0.1                      
                     -0.5 -0.9                     %
                     -4.4 -0.1  /                   %
\arrow <5pt> [.2,.67] from 3.35 -0.15  to 3.4 -0.1   %
\arrow <5pt> [.2,.67] from -4.34 -0.15  to -4.4 -0.1   %
\setquadratic \plot 2.4 -0.1                      
                     -0.5 -0.7                     %
                     -3.4 -0.1  /                   %
\arrow <5pt> [.2,.67] from 2.25 -0.15  to 2.4 -0.1   %
\arrow <5pt> [.2,.67] from -3.25 -0.15  to -3.4 -0.1   %
\setquadratic \plot 1.4 -0.1                      
                     -0.5 -0.5                     %
                     -2.4 -0.1  /                   %
\arrow <5pt> [.2,.67] from 1.25 -0.15  to 1.4 -0.1   %
\arrow <5pt> [.2,.67] from -2.25 -0.15  to -2.4 -0.1   %
\setquadratic \plot 0.4 -0.1
                     -0.5 -0.2                      
                     -1.4 -0.1         /                   %
\arrow <5pt> [.2,.67] from 0.2 -0.15  to 0.4 -0.1   %
\arrow <5pt> [.2,.67] from -1.2 -0.15  to -1.4 -0.1   %
\endpicture
$$
Note that since $\langle p_i(t),\alpha_i^\vee\rangle = -1$, $p$ is not the head of
its $i$-string and $s_i\circ p$ is well defined.  If $q=s_i\circ p$ then the first time $t$ that
$q$ leaves the cone $C-\rho$ is the same as the first time that $p$ leaves the cone $C-\rho$
and $p(t)=q(t)$.   Thus $s_i\circ q = p$ and $s_i\circ(s_i\circ p) = p$.  Since
$$s_{\wt(s_i\circ p)} = s_{s_i\circ \wt(p)} = -s_{\wt(p)}$$
the terms $s_{\wt(s_i\circ p)}$ and $s_{wt(p)}$ cancel in the sum in (5.29).  Thus
$$\mathrm{char}(B) = \sum_{p\in B\atop p\subseteq C-\rho} s_{\wt(p)}.$$
\end{proof}

\begin{thm}  Recall the notations for Weyl characters, 
tensor product multiplicities, restriction multiplicities
and paths from (5.5), (5.11), (5.17) and (5.22).
For each $\lambda\in P^+$ fix a highest weight path $p_\lambda^+$ with endpoint
$\lambda$ and let 
$$\hbox{$B(\lambda)$ be the crystal generated by $p_\lambda^+$.}$$
Let $\lambda,\mu, \nu\in P^+$ and
let  $J\subseteq \{1,2,\ldots, n\}$.  Then
$$s_\lambda = \sum_{p\in B(\lambda)} X^{\wt(p)},
\qquad
s_\mu s_\nu = \sum_{q\in B(\nu)\atop p_\mu^+\otimes q\subseteq C-\rho}
s_{\mu+\wt(q)}, 
\qquad\hbox{and}\qquad
s_\lambda = \sum_{p\in B(\lambda)\atop p\subseteq C_J-\rho_J} s_{\wt(p)}^J.$$
\end{thm}
\begin{proof}  (a) The path $p_\lambda^+$ is the unique highest weight path in $B(\lambda)$.
Thus, by Theorem 5.5, $\mathrm{char}(B(\lambda)) = s_\lambda$.

\smallskip\noindent
(b)  By the ``Leibnitz formula'' for the root operators in Theorem 5.8c the set
$$B(\mu)\otimes B(\nu) = \{ p\otimes q\ |\ p\in B(\mu), q\in B(\nu)\}$$
is a crystal.  Since $\wt(p\otimes q) = \wt(p)+\wt(q)$,
\begin{align*}
s_\mu s_\nu &= \mathrm{char}(B(\mu))\mathrm{char}(B(\nu)) 
= \mathrm{char}((B(\mu)\otimes B(\nu)) \\
&= \sum_{p\otimes q\in B(\mu)\otimes B(\nu)\atop p\otimes q\subseteq C-\rho}
s_{\wt(p)+\wt(q)} 
=\sum_{q\in B(\nu)\atop p_\mu^+\otimes q\subseteq C-\rho} s_{\mu+\wt(q)},
\end{align*}
where the third equality is from Theorem 5.5 and the last equality is because the path $p_\mu^+$
has $\wt(p_\mu^+) = \mu$ and is the only highest weight path in $B(\mu)$.

\smallskip\noindent
(c)  A \emph{$J$-crystal} is a set of paths $B$ which is closed under the operators $\tilde e_j$,
$\tilde f_j$, for $j\in J$.
Since $s_\lambda = \mathrm{char}(B(\lambda))$ the statement follows by applying Theorem
5.5 to $B(\lambda)$ viewed as a $J$-crystal.
\end{proof}

\end{subsection}

\begin{subsection}{Root operators for the rank 1 case}

Let 
$$B^{\otimes k} = \{ b_1\otimes \cdots \otimes b_k\ |\ b_i\in B\},
\qquad\hbox{where}\qquad
B = \{+1,-1,0\}.$$
Define
$$\tilde f\colon B^{\otimes k} \to B^{\otimes k}\cup\{0\}
\quad\hbox{and}\quad
\tilde e\colon B^{\otimes k}\to B^{\otimes k}\cup \{0\}$$
as follows.  Let $b=b_1\otimes \cdots \otimes b_k \in B^{\otimes k}$.
Ignoring $0$s successively pair adjacent unpaired $(+1,-1)$ pairs to
obtain a sequence of unpaired $-1$s and $+1$s
$$-1\ -1\ -1\ -1\ -1\ -1\ -1\ +1\ +1\ +1\ +1$$
(after pairing and ignoring $0$s).  Then
\begin{equation}
\begin{array}{rc}
\tilde f b &= \hbox{same as $b$ except the leftmost unpaired $+1$ is changed to $-1$}, \\
\tilde e b &= \hbox{same as $b$ except the righttmost unpaired $-1$ is changed to $+1$}.
\end{array}
\end{equation}
If there is no unpaired $+1$ after pairing then $\tilde f b =0$.  \hfil\break
If there is no unpaired $-1$ after pairing then $\tilde e b=0$.  

\smallskip\noindent
These operators coincide with the operators used in the type A case by Lascoux
and Sch\"utzenberger [LS] (see the nice exposition in [Ki]).  The $(+1,-1)$
pairing procedure is equivalent to the process of taking the ``outer edge'' of the path
(4.18-4.19).
In the context of Section 4 this is natural since only the outer edge of the path contributes 
nontrivially to the image of the path in the affine Hecke algebra.

\smallskip

Let $\fh_\RR^* = \RR$.  By identifying $+1$, $-1$, $0$ with the straight line paths
$$
\begin{matrix}
\beginpicture
\setcoordinatesystem units <1cm,1cm>         
\setplotarea x from -1.5 to 1.5, y from -0.5 to 0.5  
\put{$p_1$}[b] at 0 0.2
\arrow <5pt> [.2,.67] from -0.5 0  to 0.5 0   %
\endpicture
&\beginpicture
\setcoordinatesystem units <1cm,1cm>         
\setplotarea x from -1.5 to 1.5, y from -0.5 to 0.5  
\put{$p_{-1}$}[b] at 0 0.2
\arrow <5pt> [.2,.67] from 0.5 0  to -0.5 0   %
\endpicture
&
\beginpicture
\setcoordinatesystem units <1cm,1cm>         
\setplotarea x from -1.5 to 0.5, y from -0.5 to 0.5  
\put{$p_0$}[b] at 0 0.2
\put{$\bullet$}[r] at  0 0
\endpicture 
\\
\begin{matrix}
p_1\colon &[0,1] &\to &\fh_\RR^* \\
&t &\mapsto &t,
\end{matrix}
\qquad
&\begin{matrix}
p_{-1}\colon &[0,1] &\to &\fh_\RR^* \\
&t &\mapsto &-t,
\end{matrix}
\qquad
&\begin{matrix}
p_0\colon &[0,1] &\to &\fh_\RR^* \\
&t &\mapsto &0,
\end{matrix}
\end{matrix}
$$
respectively, the set
$B^{\otimes k}$ is viewed as a set of maps
$p\colon [0,k]\to \fh_\RR^*$.
Let $B^{\otimes 0} = \{\phi\}$ with $\tilde f\phi = 0$ and $\tilde e\phi =0$.  Then
\begin{equation}
T_\ZZ(B) = \bigsqcup_{k\in \ZZ_{\ge 0}} B^{\otimes k}
\end{equation}
is a set of paths closed under products, reverses and $r$-stretches ($r\in \ZZ_{\ge 0}$)
and the action of $\tilde e$ and $\tilde f$.  For $p\in B$ let
\begin{equation}
\begin{array}{rc}
d^+(p) &= \hbox{(number of unpaired $+1$s after pairing)}, \\
d^-(p) &= \hbox{(number of unpaired $-1$s after pairing)}.
\end{array}
\end{equation}
See the picture in (5.33).
These are the nonnegative integers such that
$$
\tilde f^{d^+(p)}p\ne 0\quad\hbox{and}\quad \tilde f^{d^+(p)+1}p=0,
\qquad\quad\hbox{and}\qquad\quad
\tilde e^{d^-(p)}p\ne 0\quad\hbox{and}\quad \tilde e^{d^-(p)+1}p=0.
$$

\begin{prop}  Use notations for $T_{\ZZ}(B)$ as in (5.30-5.32).
\item[(a)] If $p\in T_\ZZ(B)$ and $\tilde fp\ne 0$ then $\tilde e\tilde f p = p$.
\item[\phantom{(a)}] If $p\in T_{\ZZ}(B)$ and $\tilde ep\ne 0$ then $\tilde f\tilde e p = p$.
\item[(b)] If $p\in T_{\ZZ}(B)$ and $r\in \ZZ_{\ge 0}$ then
$$\tilde f^r(rp) = r(\tilde fp)
\qquad\hbox{and}\qquad
\tilde e^r(rp) = r(\tilde ep).$$
\item[(c)] If $p,q\in T_\ZZ(B)$ then
$$
\tilde f(p\otimes q) = \begin{cases}
\tilde fp\otimes q, &\hbox{if $d^+(p)>d^-(q)$}, \\
p\otimes \tilde f q, &\hbox{if $d^+(p)\le d^-(q)$},
\end{cases}
\qquad\hbox{and}\qquad
\tilde e(p\otimes q) = \begin{cases}
\tilde ep\otimes q, &\hbox{if $d^+(p)\ge d^-(q)$}, \\
p\otimes \tilde e q, &\hbox{if $d^+(p)< d^-(q)$}.
\end{cases}
$$
\item[(d)] If $p\in T_\ZZ(B)$ then
$$\tilde f(p^*)=(\tilde ep)^*\qquad\hbox{and}\qquad \tilde e(p^*)=(\tilde fp)^*.$$
\end{prop}
\begin{proof}
(a), (b) and (d) are direct consequences of the definition of the
operators $\tilde e$ and $\tilde f$ and
the definitions of $r$-stretching and reversing.

\smallskip\noindent
(c)  View $p$ and $q$ as paths.  After pairing, $p$ and $q$ have the form
\begin{equation}
p = 
\beginpicture
\setcoordinatesystem units <1cm,1cm>         
\setplotarea x from -1.5 to 1.5, y from -0.5 to 0.5  
\put{$\bullet$} at 0 0 
\plot -1 0  0 0 /
\plot -1 0  -1 0.1 /
\arrow <5pt> [.2,.67] from -1 0.1  to 1.2 0.1   %
\put{$\vert$}[b] at  -1 0.15
\put{$\vert$}[b] at  1.2 0.15
\put{$\rightarrow$}[r] at  1.15 0.4
\put{$\leftarrow$}[l] at  -1 0.4
\put{$\scriptstyle{d^+(p)}$}[b] at  0.1 0.2
\put{$\vert$}[t] at  -1 -0.05
\put{$\vert$}[t] at  0 -0.05
\put{$\rightarrow$}[r] at  -0.05 -0.2
\put{$\leftarrow$}[l] at  -0.95 -0.2
\put{$\scriptstyle{d^-(p)}$}[t] at  -0.5 -0.5
\endpicture
\qquad\hbox{and}\qquad
q = 
\beginpicture
\setcoordinatesystem units <1cm,1cm>         
\setplotarea x from -1.5 to 1.5, y from -0.5 to 0.5  
\put{$\bullet$} at 1.2 0 
\plot -1 0  -1 0.1 /
\plot -1 0  1.2 0 /
\arrow <5pt> [.2,.67] from -1 0.1  to 0 0.1   %
\put{$\vert$}[t] at  -1 -0.05
\put{$\vert$}[t] at  1.2 -0.05
\put{$\rightarrow$}[r] at  1.15 -0.2
\put{$\leftarrow$}[l] at  -1 -0.2
\put{$\scriptstyle{d^-(q)}$}[t] at  0.1 -0.2
\put{$\vert$}[b] at  -1 0.15
\put{$\vert$}[b] at  0 0.15
\put{$\rightarrow$}[r] at  -0.05 0.4
\put{$\leftarrow$}[l] at  -0.95 0.4
\put{$\scriptstyle{d^+(q)}$}[b] at  -0.5 0.6
\endpicture
\end{equation}
where the left traveling portion of the path corresponds to $-1$s and the right
traveling portion of the path corresponds to $+1$s.  Then
$$\tilde f(p\otimes q) = \begin{cases}
\tilde fp\otimes q, &\hbox{if $\displaystyle{
p\otimes q = 
\beginpicture
\setcoordinatesystem units <1cm,1cm>         
\setplotarea x from -1.5 to 1, y from -0.5 to 0.5  
\put{$\bullet$} at 0 0 
\plot -1 0  0 0 /
\plot -1 0  -1 0.1 /
\plot -1 0.1  0.5 0.1 /
\plot 0.5 0.1  0.5 0.2 /
\plot 0.5 0.2  -0.5 0.2 /
\plot -0.5 0.2  -0.5 0.3 /
\arrow <5pt> [.2,.67] from -0.5 0.3  to 0.8 0.3   %
\endpicture
}$,\quad i.e. $d^+(p)>d^-(q)$}, \\
p\otimes \tilde fq, &\hbox{if $\displaystyle{
p\otimes q = 
\beginpicture
\setcoordinatesystem units <1cm,1cm>         
\setplotarea x from -1.5 to 0.7, y from -0.5 to 0.5  
\put{$\bullet$} at 0 0 
\plot -0.4 0  0 0 /
\plot -0.4 0  -0.4 0.1 /
\plot -0.4 0.1  0.5 0.1 /
\plot 0.5 0.1  0.5 0.2 /
\plot 0.5 0.2  -1.2 0.2 /
\plot -1.2 0.2  -1.2 0.3 /
\arrow <5pt> [.2,.67] from -1.2 0.3  to -0.3 0.3   %
\endpicture
}$,\quad i.e. $d^+(p)\le d^-(q)$}, 
\end{cases}
$$ since, in the first case, the leftmost unpaired $+1$ is from $p$ and,
in the second case, it is from $q$.
\end{proof}

Use property (b) in Proposition 5.7 to extend the operators $\tilde e$ and 
$\tilde f$ to operators on $T_\QQ(B)$, the set of maps $p\colon [0,\ell]\to \RR$
generated by $B$ under the operations of concatentation, reversing and $r$-stretching
($r\in \QQ_{\ge 0}$).  Then, by completion,  the operators $\tilde e $ and $\tilde f$ extend
to operators on 
\begin{equation}
\begin{array}{rl}
T_\RR(B),\ \ &\hbox{the set of maps $p\colon [0,\ell]\to \RR$
generated by $B$} \\
&\hbox{by operations of concatenation, reversing and $r$-stretching ($r\in \RR_{\ge 0}$).}
\end{array}
\end{equation}
A \emph{rank 1 path} is an element of $T_\RR(B)$.

\end{subsection}

\begin{subsection}{The root operators in the general case}

Recall from (5.23) that
\smallskip
$$
\begin{array}{rl}
B_{\mathrm{univ}}\ \ &\hbox{is the set of maps generated by the straight line paths} \\
&\hbox{by operations of concatenation, reversing and stretching.}
\end{array}
$$
and a \emph{path} is an element $p\colon [0,\ell]\to \fh_\RR^*$ in $B_{\mathrm{univ}}$
(see (5.23)).

Let $p\colon [0,\ell]\to \RR$ be a path and let $\alpha\in R^+$ be a positive root.  The
map 
\begin{equation}
p_\alpha\colon [0,\ell]\to \RR
\quad\hbox{given by}\quad
p_\alpha(t) = \langle p(t),\alpha^\vee\rangle
\end{equation}
is a rank 1 path (an element of $T_\RR(B)$).  The rank 1 path
$p_\alpha$ is the projection of $p$ onto the line perpendicular to the hyperplane
$H_{\alpha}$.  Define operators
\begin{equation}
\tilde e_\alpha\colon B_{\mathrm{univ}}\to B_{\mathrm{univ}}\cup \{0\}
\quad\hbox{and}\quad
\tilde f_\alpha\colon B_{\mathrm{univ}}\to B_{\mathrm{univ}}\cup \{0\}
\end{equation}
by
\begin{equation}
\tilde e_\alpha p = p+\hbox{$\frac{1}{2}$}(\tilde e p_\alpha-p_\alpha)\alpha
\qquad\hbox{and}\qquad
\tilde f_\alpha p = p-\hbox{$\frac{1}{2}$}(p_\alpha-\tilde f p_\alpha)\alpha,
\end{equation}
and set
\begin{equation}
\tilde e_i = \tilde e_{\alpha_i}
\quad\hbox{and}\quad
\tilde f_i = \tilde f_{\alpha_i},
\qquad\hbox{for $1\le i\le n$}.
\end{equation}
The operators $\tilde e_i$ and $\tilde f_i$ are designed so that 
after projection onto the line perpendicular to
$H_{\alpha_i}$ they are the operators $\tilde e$ and
$\tilde f$.

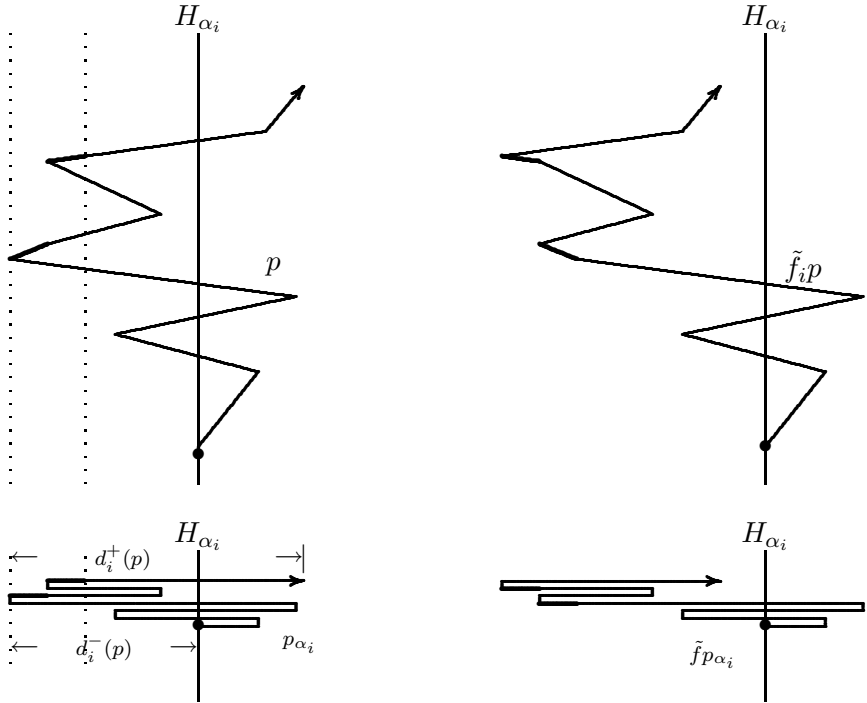
$$
\beginpicture
\setcoordinatesystem units <1cm,1cm>         
\setplotarea x from -3 to 3, y from -3 to 3    
\linethickness=0.5pt                          
\putrule from 0 3 to 0 -3        
\put{$H_{\alpha_i}$} at 0 3.2          
\put{$\bullet$} at 0 -2.6
\plot  0 -2.5   .8 -1.5 /    
\plot .8 -1.5  -1.1 -1 /    
\plot -1.1 -1  1.3 -.5 /    
\plot 1.3 -.5  -2.5 0  /    
\plot -2 .2  -.5 .6 /    
\plot  -.5 .6  -2 1.3 /    
\plot   -1.5 1.36896  .9  1.7 /    
\plot   .9 1.7  1.4  2.3 /    
\arrow <5pt> [.2,.67] from 1.32 2.2  to 1.4 2.3   %
\put{$p$} at 1 -0.1
\setplotsymbol({\bf .})                   
\plot -2.5 0  -2 .2 /    
\plot   -2 1.3  -1.5  1.36896 /    
\setdots
\putrule from -1.5 3 to -1.5 -3        
\putrule from -2.5 3 to -2.5 -3        
\endpicture
\qquad\qquad
\beginpicture
\setcoordinatesystem units <1cm,1cm>         
\setplotarea x from -3 to 3, y from -3 to 3    
\linethickness=0.5pt                          
\putrule from 0 3 to 0 -3        
\put{$H_{\alpha_i}$} at 0 3.2          
\put{$\bullet$} at 0 -2.5
\plot  0 -2.5   .8 -1.5 /    
\plot .8 -1.5  -1.1 -1 /    
\plot -1.1 -1  1.3 -.5 /    
\plot 1.3 -.5  -2.5 0  /    
\plot -3 .2  -1.5 .6 /    
\plot  -1.5 .6  -3 1.3 /    
\plot   -3.5 1.36896  -1.1  1.7 /    
\plot   -1.1 1.7  -0.6  2.3 /    
\arrow <5pt> [.2,.67] from -0.68 2.2  to -0.6 2.3   %
\put{$\tilde f_ip$} at 0.5 -0.1
\setplotsymbol({\bf .})                   
\plot -2.5 0  -3 .2 /    
\plot   -3 1.3  -3.5  1.36896 /    
\setdots
\endpicture
$$
$$
\beginpicture
\setcoordinatesystem units <1cm,1cm>         
\setplotarea x from -3 to 3, y from -0.5 to 1.5    
\linethickness=0.5pt                          
\putrule from 0 1 to 0 -1        
\put{$H_{\alpha_i}$} at 0 1.2          
\put{$\scriptstyle{d_i^-(p)}$} at -1.25 -0.3
\put{$\leftarrow$}[l] at -2.5 -0.3
\put{$\rightarrow$}[r] at 0 -0.3
\put{$\scriptstyle{d_i^+(p)}$} at -1 0.9
\put{$\leftarrow$}[l] at -2.5 0.9
\put{$\rightarrow$}[r] at 1.4 0.9
\put{$\vert$} at 1.4 0.9
\put{$\bullet$} at 0 0
\plot  0 0   .8 0 /    
\plot  0.8 0   .8 0.1 /    
\plot .8 0.1  -1.1 0.1 /    
\plot -1.1 0.1  -1.1 0.2 /    
\plot -1.1 0.2  1.3 0.2 /    
\plot  1.3 0.2  1.3 0.3 /    
\plot 1.3 0.3  -2.5 0.3  /    
\plot  -2.5 0.3  -2.5 0.4  /    
\plot  -2.2 0.4  -.5 0.4 /    
\plot  -0.5 0.4  -.5 0.5 /    
\plot  -.5 0.5  -2 0.5 /    
\plot  -2 0.5  -2 0.6 /    
\plot   -1.5 0.6  .9  0.6 /    
\plot   .9 0.6  1.4  0.6 /    
\arrow <5pt> [.2,.67] from 1.3 0.6  to 1.4 0.6   %
\put{$\scriptstyle{p_{\alpha_i}}$} at 1.35 -0.25
\setplotsymbol({\bf .})                   
\plot -2.5 0.4  -2 0.4 /    
\plot   -2 0.6  -1.5  0.6 /    
\setdots
\putrule from -1.5 -0.5 to -1.5 1        
\putrule from -2.5 -0.5 to -2.5 1        
\endpicture
\qquad\qquad
\beginpicture
\setcoordinatesystem units <1cm,1cm>         
\setplotarea x from -3 to 3, y from -1.5 to 1.5    
\linethickness=0.5pt                          
\putrule from 0 1 to 0 -1        
\put{$H_{\alpha_i}$} at 0 1.2          
\put{$\bullet$} at 0 0
\plot  0 0   .8 0 /    
\plot  0.8 0   .8 0.1 /    
\plot .8 0.1  -1.1 0.1 /    
\plot -1.1 0.1  -1.1 0.2 /    
\plot -1.1 0.2  1.3 0.2 /    
\plot  1.3 0.2   1.3  0.3 /    
\plot 1.3  0.3  -2.5 0.3  /    
\plot  -3  0.3  -3 0.4  /    
\plot -3 0.4  -1.5 0.4 /    
\plot  -1.5 0.4  -1.5 0.5 /    
\plot  -1.5 0.5  -3 0.5 /    
\plot  -3.5 0.5  -3.5 0.6 /    
\plot   -3.5 0.6  -1.1  0.6 /    
\plot   -1.1 0.6  -0.6  0.6 /    
\arrow <5pt> [.2,.67] from -0.7 0.6  to -0.6 0.6   %
\put{$\scriptstyle{\tilde f p_{\alpha_i}}$}[t] at -0.7 -0.2
\setplotsymbol({\bf .})                   
\plot -2.5 0.3  -3 0.3 /    
\plot   -3 0.5  -3.5  0.5 /    
\setdots
\endpicture
$$
The dark parts of the path $p$ are reflected (in a mirror parallel
to $H_{\alpha_i}$) to form the path $f_ip$.  The left dotted line is
the affine hyperplane parallel to $H_{\alpha_i}$ which intersects the path $p$
at its leftmost (most negative) point (relative to
$H_{\alpha_i}$) and the distance between the dotted lines is exactly the
distance between lines of lattice points in $P$ parallel to $H_{\alpha_i}$.

The following theorem is a consequence of Proposition 5.7 and the definition
in (5.34).  The uniqueness of the operators $\tilde f_i$ and $\tilde e_i$ is forced
by the properties (b), (c) and (d) in Theorem 5.8.

\begin{thm}  The operators $\tilde e_i$ and $\tilde f_i$ defined in (5.38)
are the unique operators such that
\begin{enumerate}
\item[(a)] If $p\in B_{\mathrm{univ}}$ and $\tilde f_ip\ne 0$ then $\tilde e_i\tilde f_i p = p$.
\item[] If $p\in B_{\mathrm{univ}}$ and $\tilde e_ip\ne 0$ then $\tilde f_i\tilde e_i p = p$.
\item[(b)] If $\lambda\in P$ and $\langle \lambda,\alpha_i^\vee\rangle\in \ZZ_{>0}$ then
$$\tilde f_i^{\langle \lambda,\alpha_i^\vee\rangle} p_\lambda = p_{s_i\lambda}.$$
\item[(c)] If $p,q\in B_{\mathrm{univ}}$ then
\begin{align*}
\tilde f_i(p\otimes q) &= \begin{cases}
\tilde f_ip\otimes q, &\hbox{if $d_i^+(p)>d_i^-(q)$}, \\
p\otimes \tilde f_i q, &\hbox{if $d_i^+(p)\le d_i^-(q)$},
\end{cases}
\qquad\hbox{and}\\
\\
\tilde e_i(p\otimes q) &= \begin{cases}
\tilde e_ip\otimes q, &\hbox{if $d_i^+(p)\ge d_i^-(q)$}, \\
p\otimes \tilde e_i q, &\hbox{if $d_i^+(p)< d_i^-(q)$}.
\end{cases}
\end{align*}
where $d_i^{\pm}(p) = d^{\pm}(p_{\alpha_i})$ with $d^\pm$ as in 
(5.33) and $p_{\alpha_i}$ as in (5.35).
\item[(d)] If $p\in B_{\mathrm{univ}}$ and $r\in \ZZ_{\ge 0}$ then
$$\tilde f_i^r(rp) = r(\tilde f_ip)
\qquad\hbox{and}\qquad
\tilde e_i^r(rp) = r(\tilde e_ip).$$
\item[(e)] If $p\in B_{\mathrm{univ}}$ then
$$\tilde f_ip^*=(\tilde e_ip)^*\qquad\hbox{and}\qquad \tilde e_ip^*=(\tilde f_ip)^*.$$
\item[(f)] If $p\in B_{\mathrm{univ}}$ and $\tilde f_i p \ne 0$ then $\wt(\tilde f_i p) = \wt(p)-\alpha_i$.
\item[] If $p\in B_{\mathrm{univ}}$ and $\tilde e_i p \ne 0$ then $\wt(\tilde e_i p) = \wt(p)+\alpha_i$.
\end{enumerate}
\end{thm}

\end{subsection}

\begin{subsection}{Column strict tableaux}

A \emph{letter} is an element of $B(\varepsilon_1) = \{\varepsilon_1,\ldots, \varepsilon_n\}$
and a \emph{word of length k} is an element of 
$$B(\varepsilon_1)^{\otimes k} = \{ \varepsilon_{i_1}\otimes \cdots\otimes \varepsilon_{i_k}
\ |\ 1\le i_1,\ldots, i_k\le n\}.$$
For $1\le i\le n-1$ define
$$
\tilde f_i \colon B(\varepsilon_1)^{\otimes k} \longrightarrow B(\varepsilon_1)^{\otimes k}\cup\{0\} 
\quad\hbox{and} \quad
\tilde e_i \colon B(\varepsilon_1)^{\otimes k} \longrightarrow B(\varepsilon_1)^{\otimes k}\cup\{0\} 
$$
as follows.  For $b\in B(\varepsilon_1)^{\otimes k}$,
\begin{equation}
\begin{array}{l}
\hbox{place $+1$ under each $\varepsilon_i$ in $b$,} \\
\hbox{place $-1$ under each $\varepsilon_{i+1}$ in $b$, and} \\
\hbox{place $0$ under each $\varepsilon_j$, $j\ne i, i+1$.}
\end{array}
\end{equation}
Ignoring $0$s, successively pair adjacent $(+1,-1)$ pairs to obtain a sequence of unpaired
$-1$s and $+1$s
\begin{equation}
-1\ -1\ -1\ -1\ -1\ -1\ -1\ +1\ +1\ +1\ +1
\end{equation}
(after pairing and ignoring $0$s).  Then
\begin{align*}
\tilde f_i b &= \hbox{same as $b$ except the letter corresponding to the leftmost unpaired $+1$ is changed to $\varepsilon_{i+1}$,} \\
\tilde e_i b &= \hbox{same as $b$ except the letter corresponding to the rightmost unpaired $-1$ is 
changed to $\varepsilon_i$}.
\end{align*}
If there is no unpaired $+1$ after pairing then $\tilde f_i b =0$.  \hfil\break
If there is no unpaired $-1$ after pairing then $\tilde e_i b=0$.  

A \emph{partition} is a collection $\mu$ of boxes in a corner
where the convention is that gravity goes
up and to the left.  As for matrices,
the rows and columns of $\mu$ are indexed from
top to bottom and left to right, respectively.
\begin{equation}
\begin{array}{ll}
\hbox{The \emph{parts} of $\mu$ are }
&\mu_i = \hbox{(the number of boxes in row $i$ of $\mu$)}, \\
\hbox{the \emph{length} of $\mu$ is} 
&\ell(\mu) = \hbox{(the number of rows of $\mu$)}, \\
\hbox{the \emph{size} of $\mu$ is }
&|\mu| = \mu_1+\cdots+\mu_{\ell(\mu)} 
= \hbox{(the number of boxes of $\mu$)}.
\end{array}
\end{equation}
Then $\mu$ is determined by (and identified with)
the sequence $\mu=(\mu_1,\ldots,\mu_\ell)$ 
of positive integers such that $\mu_1\ge \mu_2\ge \cdots \ge \mu_\ell >0$,
where $\ell = \ell(\mu)$.
For example,
$$
\beginpicture
\setcoordinatesystem units <0.25cm,0.25cm>         
\setplotarea x from -7 to 7, y from 0 to 6    
\linethickness=0.5pt                          
\putrule from 0 6 to 5 6          %
\putrule from 0 5 to 5 5          
\putrule from 0 4 to 5 4          %
\putrule from 0 3 to 3 3          %
\putrule from 0 2 to 3 2          %
\putrule from 0 1 to 1 1          %
\putrule from 0 0 to 1 0          %
\putrule from 0 0 to 0 6        %
\putrule from 1 0 to 1 6        %
\putrule from 2 2 to 2 6        %
\putrule from 3 2 to 3 6        
\putrule from 4 4 to 4 6        %
\putrule from 5 4 to 5 6        %
\put{$(5,5,3,3,1,1)=$} at -6 3.5
\put{.} at 6 3
\endpicture
$$

Let $\lambda$ be a partition and let 
$\mu = (\mu_1,\ldots, \mu_n)\in \ZZ_{\ge 0}^n$ 
be a sequence of nonnegative integers.
A \emph{column strict tableau of shape $\lambda$ and
weight $\mu$} is a filling of the boxes of $\lambda$
with $\mu_1$ ~1s, $\mu_2$ ~2s, $\ldots$,
$\mu_n$ ~$n$s, such that
\begin{enumerate}
\item[(a)]  the rows are weakly increasing from left to right,
\item[(b)]  the columns are strictly increasing from top to bottom.
\end{enumerate}
If $p$ is a column strict tableau write
$\mathrm{shp}(p)$ and $\mathrm{wt}(p)$ for the shape and the
weight of $p$ so that 
\begin{align*}
\mathrm{shp}(p) &= (\lambda_1,\ldots, \lambda_n),
\qquad\hbox{where}\quad
\lambda_i = \hbox{number of boxes in row $i$ of $p$,\quad and} \\
\mathrm{wt}(p) &= (\mu_1,\ldots,\mu_n),
\qquad\hbox{where}\quad
\mu_i = \hbox{number of $i\,$s in $p$}. 
\end{align*}
For example,
$$
\beginpicture
\setcoordinatesystem units <0.5cm,0.5cm>         
\setplotarea x from -2 to 15, y from 0 to 3    
\linethickness=0.5pt                          
\put{$p=$} at -1.5 4.5
\put{has\quad ${\rm shp}(p) = (9,7,7,4,2,1,0)$\quad and}[l] at 11 3.5
\put{\phantom{has}\quad ${\rm wt}(p) = (7,6,5,5,3,2,2)$.}[l] at  11 2.5
\put{7} at 0.5 0.5
\put{6} at 0.5 1.5
\put{4} at 0.5 2.5
\put{3} at 0.5 3.5
\put{2} at 0.5 4.5
\put{1} at 0.5 5.5
\put{7} at 1.5 1.5
\put{5} at 1.5 2.5
\put{3} at 1.5 3.5
\put{2} at 1.5 4.5
\put{1} at 1.5 5.5
\put{5} at 2.5 2.5
\put{3} at 2.5 3.5
\put{2} at 2.5 4.5
\put{1} at 2.5 5.5
\put{6}  at 3.5 2.5
\put{4}  at 3.5 3.5
\put{2}  at 3.5 4.5
\put{1}  at 3.5 5.5
\put{4} at 4.5 3.5
\put{3} at 4.5 4.5
\put{1} at 4.5 5.5
\put{4}  at 5.5 3.5
\put{3}  at 5.5 4.5
\put{1}  at 5.5 5.5
\put{5}  at 6.5 3.5
\put{4}  at 6.5 4.5
\put{1}  at 6.5 5.5
\put{2}  at 7.5 5.5
\put{2}  at 8.5 5.5
\putrule from 0 6 to 9 6          %
\putrule from 0 5 to 9 5          
\putrule from 0 4 to 7 4          %
\putrule from 0 3 to 7 3          %
\putrule from 0 2 to 4 2          %
\putrule from 0 2 to 2 2          %
\putrule from 0 1 to 2 1          %
\putrule from 0 0 to 1 0          %
\putrule from 0 0 to 0 6        %
\putrule from 1 0 to 1 6        %
\putrule from 2 1 to 2 6        %
\putrule from 3 2 to 3 6        
\putrule from 4 2 to 4 6        %
\putrule from 5 3 to 5 6        %
\putrule from 6 3 to 6 6        %
\putrule from 7 3 to 7 6        %
\putrule from 8 5 to 8 6        %
\putrule from 9 5 to 9 6        %
\endpicture
$$
For a partition $\lambda$ and a 
sequence $\mu = (\mu_1,\ldots,\mu_n)\in\ZZ^n_{\ge 0}$
of nonnegative integers write
\begin{equation}
\begin{array}{rl}
B(\lambda) &= \{ 
\hbox{column strict tableaux $p$}\ |\ 
\mathrm{shp}(p) = \lambda\}, \\
B(\lambda)_\mu &= \{ 
\hbox{column strict tableaux $p$}\ |\ 
\mathrm{shp}(p) = \lambda\ \hbox{and}\ {\rm wt}(p)=\mu\}.
\end{array}
\end{equation}

Let $\lambda$ be a partition with $k$ boxes and let 
$$B(\lambda) = \{ \hbox{column strict tableaux of shape $\lambda$}\}.$$
The set $B(\lambda)$ is a subset of $B(\varepsilon_1)^{\otimes k}$ via the injection
\begin{equation}
\begin{matrix}
B(\lambda) &\hookrightarrow &B(\varepsilon_1)^{\otimes k} \\
p&\longmapsto &(\hbox{the arabic reading of $p$}) \\
\\
\beginpicture
\setcoordinatesystem units <0.5cm,0.5cm>         
\setplotarea x from 0 to 8, y from -3 to 3    
\linethickness=0.5pt                          
\put{$\scriptstyle{i_k}$} at 0.5 -2.5
\put{$\scriptstyle{i_{\lambda_1+\lambda_2}}$}[l] at 0.3 1.5
\put{$\scriptstyle{i_{\lambda_1}}$}[l] at 0.3 2.5
\put{$\cdots$} at 2 -0.2
\put{$\cdots$}  at 3.5 1.5
\put{$\cdots$}  at 3.5 2.5
\put{$\scriptstyle{i_{\lambda_1+1}}$}[r]  at 6.7 1.5
\put{$\scriptstyle{i_2}$}  at 7.5 2.5
\put{$\scriptstyle{i_1}$}  at 8.5 2.5
\putrule from 0 3 to 9 3          %
\putrule from 0 2 to 9 2          
\putrule from 0 1 to 7 1          %
\putrule from 4 0 to 7 0          %
\putrule from 2 -1 to 4 -1          %
\putrule from 1 -2 to 2 -2          %
\putrule from 0 -3 to 1 -3          %
\putrule from 0 -3 to 0 3        %
\putrule from 1 -3 to 1 -2        %
\putrule from 1 1 to 1 1        %
\putrule from 2 -2 to 2 -1        %
\putrule from 4 -1 to 4 0        %
\putrule from 7 0 to 7 3        %
\putrule from 8 2 to 8 3        %
\putrule from 9 2 to 9 3        %
\endpicture
&\longmapsto 
&\varepsilon_{i_1}\otimes \varepsilon_{i_2}\otimes \cdots\otimes \varepsilon_{i_k}
\end{matrix}
\end{equation}
where the arabic reading of $p$ is 
$\varepsilon_{i_1}\otimes\varepsilon_{i_2}\otimes\cdots\otimes \varepsilon_{i_k}$
if the entries of $p$ are $i_1, i_2,\ldots, i_k$ read right to left by rows with the rows read in sequence
beginning with the first row.

\begin{prop}  Let $\lambda = (\lambda_1,\ldots, \lambda_n)$ be a partition with $k$ boxes.
Then $B(\lambda)$ is the subset of $B(\varepsilon_1)^{\otimes k}$ generated by
$$p_\lambda = 
\underbrace{
\varepsilon_1\otimes \varepsilon_1\otimes\cdots\otimes 
\varepsilon_1}_{\lambda_1\ \mathrm{factors}}
\otimes\underbrace{
\varepsilon_2\otimes \varepsilon_2\otimes\cdots\otimes 
\varepsilon_2}_{\lambda_2\ \mathrm{factors}}
\otimes\cdots\otimes
\underbrace{
\varepsilon_n\otimes \varepsilon_n\otimes\cdots\otimes 
\varepsilon_n}_{\lambda_n\ \mathrm{factors}}
$$
under the action of the operators $\tilde e_i$, $\tilde f_i$, $1\le i\le n$.
\end{prop}
\begin{proof}
If $P = P(b)$ is a filling of the shape $\lambda$ then
$b_{i_1}\otimes \cdots\otimes b_{i_k} = b$ is obtained from $P$ by reading
the entries of $P$ in arabic reading order (right to left across rows and from top to
bottom down the page).  The tableau
$$P_\lambda = P(p_\lambda) =
\beginpicture
\setcoordinatesystem units <0.5cm,0.5cm>         
\setplotarea x from 0 to 8, y from -3 to 3    
\linethickness=0.5pt                          
\put{$\scriptstyle{n}$} at 0.5 -1.5
\put{$\scriptstyle{n}$} at 2.5 -1.5
\put{$\cdots$} at 1.5 -1.5
\put{$\scriptstyle{2}$}[l] at 0.3 1.5
\put{$\scriptstyle{2}$}[l] at 1 1.5
\put{$\scriptstyle{1}$}[l] at 0.3 2.5
\put{$\scriptstyle{1}$}[l] at 1 2.5
\put{$\cdots$} at 2 -0.2
\put{$\cdots$}  at 3.5 1.5
\put{$\cdots$}  at 4.5 2.5
\put{$\scriptstyle{2}$}[r]  at 5.7 1.5
\put{$\scriptstyle{2}$}[r]  at 6.7 1.5
\put{$\scriptstyle{1}$}  at 6.7 2.5
\put{$\scriptstyle{1}$}  at 7.7 2.5
\put{$\scriptstyle{1}$}  at 8.5 2.5
\putrule from 0 3 to 9 3          %
\putrule from 0 2 to 9 2          
\putrule from 0 1 to 7 1          %
\putrule from 4 0 to 7 0          %
\putrule from 3 -1 to 4 -1          %
\putrule from 0 -2 to 3 -2          %
\putrule from 0 -2 to 0 3        %
\putrule from 1 1 to 1 1        %
\putrule from 3 -2 to 3 -1        %
\putrule from 4 -1 to 4 0        %
\putrule from 7 0 to 7 2        %
\putrule from 9 2 to 9 3        %
\endpicture
$$
is the column strict tableau of shape $\lambda$ with $1$s in the first row, 
$2$s in the second row, and so on.  Define operators $\tilde e_i$ and $\tilde f_i$ on a
filling of $\lambda$ by
$$\tilde e_i P = P(\tilde e_i p)
\quad\hbox{and}\quad
\tilde f_i P = P(\tilde f_i b),
\qquad
\hbox{if $P = P(b)$.}$$
To prove the proposition we shall show that if $P$ is a column strict
tableau of shape $\lambda$ then
\begin{enumerate}
\item[(a)] $\tilde e_i P$ and $\tilde f_i P$ are column strict tableaux,
\item[(b)] $P$ can be obtained by applying a sequence of $\tilde f_i$ to $P_\lambda$.
\end{enumerate}
Let $P^{(j)}$ be the column strict tableau formed by the entries of $P$ which are $\le j$ and let
$\lambda^{(j)} = \mathrm{shp}(P^{(j)})$.  This conversion identifies $P$ with the sequence
$$\begin{array}{l}
P = \big( \emptyset = \lambda^{(0)}\subseteq \lambda^{(1)}\subseteq
\cdots \subseteq \lambda^{(n)} = \lambda \big),
\qquad\hbox{where} \\
\lambda^{(i)}/\lambda^{(i-1)} \hbox{is a \emph{horizontal strip} for each $1\le i\le n$.}
\end{array}
$$

\smallskip\noindent
(a) Let us analyze the action of $\tilde e_i$ and $\tilde f_i$ on $P$.  The
sequence of $+1, -1, 0$ constructed via (5.39) is given by
$$\begin{array}{l}
\hbox{placing $+1$ in each box of $\lambda^{(i)}/\lambda^{(i-1)}$,} \\
\hbox{placing $-1$ in each box of $\lambda^{(i+1)}/\lambda^{(i)}$,} \\
\hbox{placing $0$ in each box of $\lambda^{(j)}/\lambda^{(j-1)}$, for $j\ne i, i+1$,}
\end{array}
$$
and reading the resulting filling in Arabic reading order.  
The process of removing $(+1,-1)$
pairs can be executed on the horizontal strips $\lambda^{(i+1)}/\lambda^{(i)}$ and
$\lambda^{(i)}/\lambda^{(i-1)}$, 
$$\lambda^{(i+1)} = 
\beginpicture
\setcoordinatesystem units <0.5cm,0.25cm>         
\setplotarea x from -2 to 20, y from 0 to 6    
\linethickness=0.5pt                          
\putrule from 0 6 to 17 6          %
\putrule from 12 5 to 17 5          
\putrule from 12 4 to 14 4          %
\putrule from 9 3 to 12 3          %
\putrule from 9 2 to 10.5 2          %
\putrule from 4 -1 to 9 -1         %
\putrule from 4 -2 to 8 -2          %
\putrule from 4 -3 to 6 -3          %
\putrule from 0.5 -5 to 4 -5          %
\putrule from 0.5 -6 to 2 -6          %
\putrule from 0 -8 to 0.5 -8          %
\putrule from 17 5 to 17 6        %
\putrule from 15 5 to 15 6        %
\putrule from 14 4 to 14 6        %
\putrule from 12 3 to 12 6        %
\putrule from 10.5 2 to 10.5 3        %
\putrule from 9 -1 to 9 3        %
\putrule from 8 -2 to 8 -1    %
\putrule from 6 -3 to 6 -1        %
\putrule from 4 -5 to 4 -1    %
\putrule from 2 -6 to 2 -5        %
\putrule from 0.5 -8 to 0.5 -5    %
\putrule from 0 -8 to 0 6        
\put{$\lambda^{(i-1)}$}[tl] at 1.5 3.5
\put{$\scriptstyle{+}$} at 0.75 -5.5
\put{$\scriptstyle{+}$} at 1.25 -5.5
\put{$\scriptstyle{+}$} at 1.75 -5.5
\put{$\scriptstyle{-}$} at 4.25 -2.5
\put{$\scriptstyle{-}$} at 4.75 -2.5
\put{$\scriptstyle{-}$} at 5.25 -2.5
\put{$\scriptstyle{-}$} at 5.75 -2.5
\put{$\scriptstyle{+}$} at 4.25 -1.5
\put{$\scriptstyle{+}$} at 4.75 -1.5
\put{$\scriptstyle{+}$} at 5.25 -1.5
\put{$\scriptstyle{+}$} at 5.75 -1.5
\put{$\scriptstyle{+}$} at 6.25 -1.5
\put{$\scriptstyle{+}$} at 6.75 -1.5
\put{$\scriptstyle{+}$} at 7.25 -1.5
\put{$\scriptstyle{+}$} at 7.75 -1.5
\put{$\scriptstyle{-}$} at 9.25 2.5
\put{$\scriptstyle{-}$} at 9.75 2.5
\put{$\scriptstyle{-}$} at 10.25 2.5
\put{$\scriptstyle{-}$} at 12.25 4.5
\put{$\scriptstyle{-}$} at 12.75 4.5
\put{$\scriptstyle{-}$} at 13.25 4.5
\put{$\scriptstyle{-}$} at 13.75 4.5
\put{$\scriptstyle{+}$} at 12.25 5.5
\put{$\scriptstyle{+}$} at 12.75 5.5
\put{$\scriptstyle{+}$} at 13.25 5.5
\put{$\scriptstyle{+}$} at 13.75 5.5
\put{$\scriptstyle{+}$} at 14.25 5.5
\put{$\scriptstyle{+}$} at 14.75 5.5
\put{$\scriptstyle{-}$} at 15.25 5.5
\put{$\scriptstyle{-}$} at 15.75 5.5
\put{$\scriptstyle{-}$} at 16.25 5.5
\put{$\scriptstyle{-}$} at 16.75 5.5
\endpicture
$$
with the effect that the entries in any configuration of boxes of the form
$$\beginpicture
\setcoordinatesystem units <0.5cm,0.5cm>         
\setplotarea x from 0 to 6, y from -1 to 1    
\linethickness=0.5pt                          
\put{$\scriptstyle{-1}$}  at 0.5 -0.5
\put{$\scriptstyle{-1}$}  at 1.5 -0.5
\put{$\cdots$}  at 3.5 -0.5
\put{$\scriptstyle{-1}$}  at 5.5 -0.5
\put{$\scriptstyle{+1}$}  at 0.5 0.5
\put{$\scriptstyle{+1}$}  at 1.5 0.5
\put{$\cdots$}  at 3.5 0.5
\put{$\scriptstyle{+1}$}  at 5.5 0.5
\putrule from 0 1 to 6 1          %
\putrule from 0 0 to 6 0          
\putrule from 0 -1 to 6 -1          %
\putrule from 0 -1 to 0 1        %
\putrule from 1 -1 to 1 1        %
\putrule from 2 -1 to 2 1        %
\putrule from 5 -1 to 5 1        %
\putrule from 6 -1 to 6 1        %
\endpicture
$$
will be removed.  Additional $+1,-1$ pairs will also be removed and the 
final sequence
\begin{equation}
\underbrace{-1\ -1\ \cdots\ -1}_{d_i^-(p)}
\underbrace{+1 +1\ \cdots\ +1}_{d_i^+(p)}
\end{equation}
will correspond to a configuration of the form
$$\lambda^{(i+1)} = 
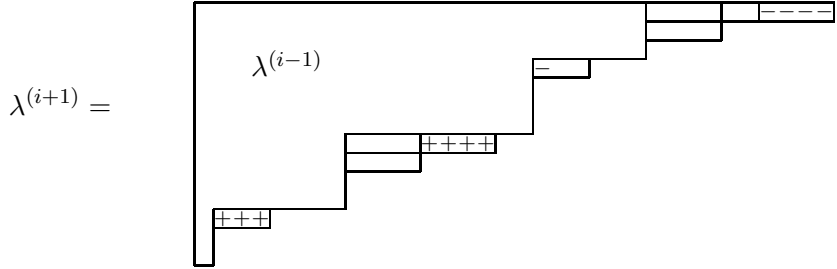
\beginpicture
\setcoordinatesystem units <0.5cm,0.25cm>         
\setplotarea x from -2 to 20, y from 0 to 6    
\linethickness=0.5pt                          
\putrule from 0 6 to 17 6          %
\putrule from 12 5 to 17 5          
\putrule from 12 4 to 14 4          %
\putrule from 9 3 to 12 3          %
\putrule from 9 2 to 10.5 2          %
\putrule from 4 -1 to 9 -1         %
\putrule from 4 -2 to 8 -2          %
\putrule from 4 -3 to 6 -3          %
\putrule from 0.5 -5 to 4 -5          %
\putrule from 0.5 -6 to 2 -6          %
\putrule from 0 -8 to 0.5 -8          %
\putrule from 17 5 to 17 6        %
\putrule from 15 5 to 15 6        %
\putrule from 14 4 to 14 6        %
\putrule from 12 3 to 12 6        %
\putrule from 10.5 2 to 10.5 3        %
\putrule from 9 -1 to 9 3        %
\putrule from 8 -2 to 8 -1    %
\putrule from 6 -3 to 6 -1        %
\putrule from 4 -5 to 4 -1    %
\putrule from 2 -6 to 2 -5        %
\putrule from 0.5 -8 to 0.5 -5    %
\putrule from 0 -8 to 0 6        
\put{$\lambda^{(i-1)}$}[tl] at 1.5 3.5
\put{$\scriptstyle{+}$} at 0.75 -5.5
\put{$\scriptstyle{+}$} at 1.25 -5.5
\put{$\scriptstyle{+}$} at 1.75 -5.5
\put{$\scriptstyle{}$} at 4.25 -2.5
\put{$\scriptstyle{}$} at 4.75 -2.5
\put{$\scriptstyle{}$} at 5.25 -2.5
\put{$\scriptstyle{}$} at 5.75 -2.5
\put{$\scriptstyle{}$} at 4.25 -1.5
\put{$\scriptstyle{}$} at 4.75 -1.5
\put{$\scriptstyle{}$} at 5.25 -1.5
\put{$\scriptstyle{}$} at 5.75 -1.5
\put{$\scriptstyle{+}$} at 6.25 -1.5
\put{$\scriptstyle{+}$} at 6.75 -1.5
\put{$\scriptstyle{+}$} at 7.25 -1.5
\put{$\scriptstyle{+}$} at 7.75 -1.5
\put{$\scriptstyle{-}$} at 9.25 2.5
\put{$\scriptstyle{}$} at 9.75 2.5
\put{$\scriptstyle{}$} at 10.25 2.5
\put{$\scriptstyle{}$} at 12.25 4.5
\put{$\scriptstyle{}$} at 12.75 4.5
\put{$\scriptstyle{}$} at 13.25 4.5
\put{$\scriptstyle{}$} at 13.75 4.5
\put{$\scriptstyle{}$} at 12.25 5.5
\put{$\scriptstyle{}$} at 12.75 5.5
\put{$\scriptstyle{}$} at 13.25 5.5
\put{$\scriptstyle{}$} at 13.75 5.5
\put{$\scriptstyle{}$} at 14.25 5.5
\put{$\scriptstyle{}$} at 14.75 5.5
\put{$\scriptstyle{-}$} at 15.25 5.5
\put{$\scriptstyle{-}$} at 15.75 5.5
\put{$\scriptstyle{-}$} at 16.25 5.5
\put{$\scriptstyle{-}$} at 16.75 5.5
\endpicture
$$
The rightmost $-1$ in the sequence (5.40) corresponds to a box in 
$\lambda^{(i+1)}/\lambda^{(i)}$ which is leftmost in its row and which does not cover a box
of $\lambda^{(i)}/\lambda^{(i-1)}$.  Similarly the leftmost $+1$ in the sequence (5.40) corrresponds
to a box in $\lambda^{(i)}/\lambda^{(i-1)}$ which is rightmost in its row and which does not
have a box of $\lambda^{(i+1)}/\lambda^{(i)}$ covering it.  These conditions guarantee that
$\tilde e_iP$ and $\tilde f_iP$ are column strict tableaux.

\smallskip\noindent
(b)  
The tableau $P$ is obtained from $P_\lambda$ by applying a sequence of $\tilde f_i$ in the 
following way.  Applying the operator 
$$\tilde f_{n,i} = \tilde f_{n-1}\cdots \tilde f_{i+1} \tilde f_i
\quad\hbox{to}\quad P_\lambda$$
will change the rightmost $i$ in row $i$ to $n$.  A sequence of applications of 
$$\tilde f_{n,i},\quad\hbox{as $i$ decreases (weakly) from $n-1$ to $1$,}$$
can be used to produce a column strict tableau $P_n$ in which
\begin{enumerate}
\item[(1)]  the entries equal to $n$ match the entries equal to $n$ in $P$, and
\item[(2)] the subtableau of $P_n$ containing the entries $\le n-1$ is 
$P_{\lambda^{(n-1)}}$.
\end{enumerate}
Iterating the process and applying a sequence of operators
$$\tilde f_{n-1,i},\quad\hbox{as $i$ decreases (weakly) from $n-2$ to $1$,}$$
to the tableau $P_n$ can be used to produce a tableau $P_{n-1}$
in which
\begin{enumerate}
\item[(1)]  the entries equal to $n$ and $n-1$ match the entries equal to $n$ and $n-1$ in $P$, and
\item[(2)] the subtableau of $P_{n-1}$ containing the entries $\le n-2$ is 
$P_{\lambda^{(n-2)}}$.
\end{enumerate}
The tableau $P$ is obtained after a total of $n$ iterations of this process.
\end{proof}

\end{subsection}

\end{section}


\begin{thebibliography}{99}

\bibitem[BD]{BD} Y.\ Billig and M.\ Dyer, \emph{Decompositions of Bruhat type for the Kac-Moody groups},  Nova J.\ Algebra Geom.\ \textbf{3} no.\ 1 (1994), 11--31.

\bibitem[Br]{Br} M.\ Brion, \emph{Positivity in the Grothendieck group of complex flag varieties}, 
J.\ Algebra \textbf{258} no.\ 1 (2002), 137--159. 

\bibitem[GL]{GL} S.\ Gaussent and  P.\ Littelmann, \emph{LS galleries, the path model, and MV cycles},
Duke Math.\ J.\ \textbf{127} no.\ 1 (2005), 35--88.

\bibitem[GR]{GR}  S.\ Griffeth and A.\ Ram, \emph{Affine Hecke algebras and the Schubert calculus},
European J.\ Combin.\ \textbf{25} no.\ 8 (2004), 1263--1283. 

\bibitem[Ha]{Ha} T.\ Haines, \emph{Structure constants for Hecke and representation rings},
Int.\ Math.\ Res.\ Not.\ \textbf{39} (2003), 2103-2119.

\bibitem[IM]{IM}  N.\ Iwahori and H.\  Matsumoto, \emph{On some Bruhat decomposition 
and the structure of the Hecke rings of $\fp$-adic Chevalley groups},
Inst.\ Hautes \'Etudes Sci.\ Publ.\ Math.\ \textbf{25} (1965), 5--48. 


\bibitem[Ki]{Ki} K.\ Killpatrick, \emph{A Combinatorial Proof of a Recursion for the q-Kostka Polynomials},
J.\  Comb.\ Th.\  Ser.\ A \textbf{92} (2000), 29--53.

\bibitem[KM]{KM} M.\ Kapovich and J.J.\ Millson,
\emph{A path model for geodesics in Euclidean buildings and its applications to representation 
theory},  arXiv: math.RT/0411182.

\bibitem[LS]{LS} A.\ Lascoux and M.P.\ Sch\"utzenberger, \emph{Le mono\"ide plaxique}, Quad.\ Ricerce Sci.\ \textbf{109} (1981), 129--156.
 
\bibitem[LP]{LP}  C.\ Lenart and A.\ Postnikov,
\emph{Affine Weyl groups in K-theory and representation theory},
arXiv: math.RT/0309207. 

\bibitem[Li1]{Li1} P.\ Littelmann, 
\emph{A Littlewood-Richardson rule for symmetrizable Kac-Moody algebras}, 
Invent.\ Math.\ \textbf{116} (1994), 329--346. 

\bibitem[Li2]{Li2} P.\ Littelmann, 
\emph{Paths and root operators in representation theory}, 
Ann.\ Math.\ \textbf{142} (1995), 499--525. 

\bibitem[Li3]{Li3} P.\ Littelmann, 
\emph{Characters of representations and paths in $\fh_\RR^*$} , 
Proc.\ Symp.\ Pure Math.\ \textbf{61} (1997), 29--49. 

\bibitem[Lu]{Lu}  G.\ Lusztig, \emph{Affine Hecke algebras and their graded version}, 
J.\ Amer.\ Math.\ Soc.\ \textbf{2} (1989), 599--635. 

\bibitem[Mac1]{Mac1} I.G.\ Macdonald, {\sl Spherical functions on a group of p-adic type}, 
Publ. Ramanujan Institute No. 2, Madras (1971). 

\bibitem[Mac2]{Mac2} I.G.\ Macdonald, {\sl Symmetric functions and Hall polynomials}, 
Oxford Mathematical Monographs, Oxford Univ. Press, New York (second edition, 1995)
 
\bibitem[PR]{PR}  H.\ Pittie and A.\ Ram, \emph{A Pieri-Chevalley formula in the K-theory of 
a G/B-bundle},  Electronic Research Announcements \textbf{5} (1999), 102--107, and
\emph{A Pieri-Chevalley formula for K(G/B)},  arXiv: math.RT/0401332. 

\bibitem[Sc]{Sc} C.\ Schwer, \emph{Galleries, Hall-Littlewood polynomials and structure constants of the spherical Hecke algebra}, arXiv: math.CO/0506287.

\end{thebibliography}
\end{document}